%% file: st_article.tex
\documentclass[onefignum,onetabnum]{siamart190516}

\pdfoutput=1

\input{st_shared}

\ifpdf
\hypersetup{
  pdftitle={Geometrically Higher Order Unfitted Space-Time Methods for PDEs on Moving Domains},
  pdfauthor={F. Heimann, C. Lehrenfeld, J. Preu\ss{}}
}
\fi
 
\begin{document}

\maketitle

\begin{abstract}
In this paper, we propose new geometrically unfitted space-time Finite Element methods for partial differential equations posed on moving domains of higher order accuracy in space and time. As a model problem, the convection-diffusion problem on a moving domain is studied. For geometrically higher order accuracy, we apply a parametric mapping on a background space-time tensor-product mesh. Concerning discretisation in time, we consider discontinuous Galerkin, as well as related continuous (Petrov-)Galerkin and Galerkin collocation methods. For stabilisation with respect to bad cut configurations and as an extension mechanism that is required for the latter two schemes, a ghost penalty stabilisation is employed. The article puts an emphasis on the techniques that allow to achieve a robust but higher order geometry handling for smooth domains. We investigate the computational properties of the respective methods in a series of numerical experiments. These include studies in different dimensions for different polynomial degrees in space and time, validating the higher order accuracy in both variables.
\end{abstract}



\section{Introduction}
Many problems in physics, engineering, chemistry and biology can be described in terms of Partial Differential Equations (PDEs) posed on moving domains. Examples include multi-phase flows \cite{grossreusken11}, blood flow in the human heart \cite{MITTAL20161065} or evolving biological cells. Finite Element methods (FEM) provide a powerful framework to solve these problems numerically.
Traditional \emph{fitted} mesh approaches, where the computational mesh also prescribes a parametrisation of the geometry, require - when applied
to moving domain problems - mesh adaptations or remeshings in every time step. Such approaches are typically employed in Arbitrary Lagrangian-Eulerian (ALE) methods \cite{HIRT1974227}, and are computationally attractive for problems involving small deformations or without topology changes, but become difficult to realize when strong deformations or topology changes occur.

In recent years, geometrically unfitted Finite Element methods have been introduced and studied for a variety of problems under different names such as CutFEM \cite{burman15}, Finite Cell Method \cite{parvizian2007finite}, XFEM \cite{fries_xfem_review}, fictitious domain methods \cite{burman2010fictitious,BURMAN2012328} and TraceFEM \cite{olshanskii2017trace}. 
In these methods, the computational mesh is not aligned with the geometry of the problem, avoiding the necessity for a mesh generation or remeshing procedure. However, unfitted FEM have to deal with different issues such as robustness w.r.t. the position of the geometry relative to the computational mesh, the imposition of interface and boundary conditions, accurate geometry handling and time integration on moving domains. 
In this manuscript we present a class of unfitted Finite Element methods aimed at addressing the last two aspects: Robust and higher order\footnote{Here and in the following, we use the terms ``high order'' and ``higher order'' interchangeably.} space-time integration for moving domain problems.

\subsection{Approaches in the literature}
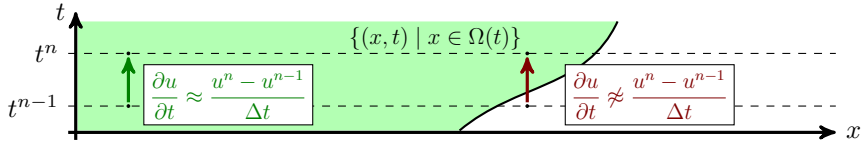
\begin{figure}
 \begin{center}  
  \begin{tikzpicture}
    [
    scale=1,
    rotate=0,
    ]
  
    \def\dw{0.025}
  
    \draw[draw=none,fill=green!30] (5.1,\dw) to[out=40,in=-115] (7.2,1.5-\dw) -- (\dw ,1.5-\dw) -- (\dw,\dw) -- cycle;

    \draw[black,thick] (5.1,\dw) to[out=40,in=-115] (7.2,1.5-\dw);

  
    \draw[->,very thick,>=stealth'] (-0.1,0)  -- (10.1,0) node(xline)[right]{$x$};
  
    \draw[->,very thick,>=stealth'] (0,-0.1)  -- (0,1.6) node(yline)[left]{$t$};

    \draw[dashed] (-0.1,0.35) node(xline)[left]{$t^{n-1}$} -- (10.1,0.35) ;

    \draw[dashed] (-0.1,1.05) node(xline)[left]{$t^{n}$} -- (10.1,1.05) ;

      \filldraw (0.7,0.35) circle (0.5pt);
      \filldraw (0.7,1.05) circle (0.5pt);
      \draw [->,green!50!black,very thick,>=stealth'] (0.7,0.4) -- (0.7,1.0);

      \filldraw (6,0.35) circle (0.5pt);
      \filldraw (6,1.05) circle (0.5pt);
      \draw [->,red!50!black,very thick,>=stealth'] (6,0.4) -- (6,1.0);

      \draw (6.5,0.5) node[red!50!black,right,fill=white,draw=black,scale=0.7] { \large $\displaystyle \frac{\partial u}{\partial t} \not\approx \frac{u^{n}-u^{n-1}}{\Delta t}$}; 
    \draw (6,1.25) node[left,scale=0.6] {\Large \color{black} $\{(x,t)\mid x\in  \Omega(t)\}$};
  \draw (0.9,0.5) node[green!50!black,right,fill=white,draw=black,scale=0.7] { \large $\displaystyle \frac{\partial u}{\partial t} \approx \frac{u^{n}-u^{n-1}}{\Delta t}$}; 
  \end{tikzpicture}
\end{center}  
\vspace*{-0.3cm}  
\caption{Sketch of the Method-of-Lines approach for moving domain problems. A finite difference stencil at a fixed point $x$ is reasonable if all involved space-time points are inside the (moving) domain (left). If this is not the case a standard finite difference stencil can not be applied (right).}   \label{fig:fd}  
\vspace*{-0.6cm}
\end{figure}

Most fitted FEM discretisations use a version of the Method-of-Lines approach to deal with the time derivative exploiting that every unknown has a well-defined history as it is associated to an entity of the mesh.
For unfitted FEM, due to the decoupling of the computational mesh and the geometry, the situation is entirely different. A naive Method-of-Lines approach may not even be well-defined when unknowns associated to entities of the computational mesh enter or leave the time-dependent domain, cf. \cref{fig:fd} for a sketch. 

In \emph{characteristic}-based methods also known as semi-Lagrangian methods such as those discussed in \cite{ma2021high,BFG22} instead of the partial derivative $\frac{ \partial}{\partial t}$ the material derivative $ \frac{d}{dt} =   \frac{\partial}{\partial t} + \vect{w}\cdot\nabla$ is approximated. The material derivative is the rate of change along the trajectories of the flow field $\vect{w}$ that also advects the geometry. Provided that these trajectories can be resolved numerically in an efficient and accurate way it can be ensured that a finite difference approximation, as in the Method-of-Lines approach, is a reasonable approximation of the material derivative. 

\emph{Extension}-based methods use a finite difference approximation of the partial time derivative, but apply an additional extension to the solution within every time step. The extension makes sure that the discrete solution is well-defined in a small neighbourhood of the domain, providing a reasonable history for each unknown involved in the finite difference approximation at the next time step. Such a method has been introduced in \cite{LO_ESAIM_2019} for a scalar convection-diffusion equation. Stokes and Navier-Stokes problems with this approach have been considered in \cite{Schott17,burman2020eulerian,IMAJNA_vWRL_2021}.

\emph{Space-time} methods do not apply finite difference-type approximations to the time derivative, but opt for a variational time discretisation which together with an FE-based space discretisation results in a variational formulation in a space-time domain. Often (but not always) the space-time domain is decomposed into \emph{time slabs} that correspond to small time intervals and a tensor-product discretisation (w.r.t. meshes and FE spaces) is chosen. This allows to keep the structure of a time stepping scheme and to build on components of related spatial problems (basis functions, quadrature rules, etc.). The former aspect implies that the PDE solution on the space-time domain can still be solved for with a series of problems with an essential complexity of a spatial PDE problem.

In this manuscript we restrict to the latter class of space-time methods in the context of geometries that are described by level set functions. 
Among this class of methods, \cite{LR_SINUM_2013,L_SISC_2015} consider Discontinuous Galerkin approaches in time of second order accuracy for two-phase flows. Moreover, in \cite{HLZ16,Z18} space-time methods for moving surface/ coupled surface-bulk problems are presented for second order accuracy as far as the level set description is concerned, and third order accuracy in the context of a geometry description involving splines.
In \cite{Lundholm2021} unfitted space-time methods based on a geometry description with overlapping meshes with piecewise linear-in-time motion are considered. Recently, in \cite{anselmann2021cutfemNSE} discontinuous-in-time methods with higher order accuracy have been investigated numerically for a Navier-Stokes problem.
%

\subsection{Main contributions}
In this manuscript, we aim to address three inherent challenges of \emph{higher order} unfitted space-time methods:

1. The derivation of discretisations that are robust, i.e.\ stable and accurate independent of the essentially arbitrary cut configurations that can occur in space-time, is a crucial component to obtain a reasonable method.

2. One inherent challenge of unfitted methods where the geometry is described by level set functions - especially for the space-time setting - is that of finding proper means to handle the \emph{implicit} space-time geometry robustly and accurately for the potentially complicated domains of numerical integration. 

3. Another challenge is the computational complexity. Although we can maintain a time stepping structure, in the higher order (in time) case we will still have considerably more unknowns to solve for in each time step than for a standard time stepping scheme. This renders space-time methods typically much more expensive compared to other time stepping schemes.

To obtain stable variational formulations, we will build on the ideas of the previous works \cite{LR_SINUM_2013,HLZ16} in the context of DG-in-time methods of lower order. One major contribution in this work is the extension of these ideas to higher order in space and time and to Petrov-Galerkin methods with discrete solutions that are (higher order) continuous in time. The transition from discontinuous solutions to continuous solutions or even solutions with higher regularity in time reduces the number of unknowns that are involved in each time step which compensates for the higher complexity of space-time methods. The second major contribution in this work is the extension to the space-time setting of a higher order accurate geometry handling based on isoparametric mappings for level set domains as introduced in \cite{L_CMAME_2016}. Thereby we can ensure to obtain arbitrary high order accuracy in space and time, including considerations of geometry approximation. Parts of the ideas in this manuscript have also been discussed in the Master's theses \cite{preuss18,heimann20}.

\subsection{Structure of this paper}
The remainder of this paper is organised as follows: In \cref{sect:modprob}, the convection-diffusion problem is introduced as a model problem. Furthermore, the different time discretisations are illustrated by semi-discrete (discrete in time) variants of the space-time methods coming up in the subsequent sections. \Cref{sect:geom_handling} is devoted to the higher order approximation of moving space-time domains described by level sets. Afterwards, in \cref{sect:weak_forms}, the fully discrete methods are defined in their discrete variational form. In \cref{subsect_num_int}, a specific strategy for higher order integration in time is discussed.
Then, \cref{sect:numexp} presents a series of numerical experiments, which demonstrate the higher order accuracy of each method, and investigate aspects of accuracy and computational complexity. Finally, \cref{sect:concl} concludes the paper with a summary and an outlook on open problems.

\section{Model problem and time discretisation} \label{sect:modprob}

\subsection{The model problem}
Before we come to the specific methods, we introduce the model problem used in this study. The model problem is given as a convection-diffusion equation of a species of concentration $u$ on a moving domain $\Omega(t)$. We denote the corresponding space-time domain as 
\hypertarget{def:Q}{$Q := \bigcup_{t} \Omega(t) \times \{ t \}$}.
As initial condition, we assume a function $u_0: \Omega(0) \to \mathbb{R}$ to be given, and for the boundary condition we assume that no transport of the species across the domain boundary occurs. We arrive at the following problem: Find $u: Q \to \mathbb{R}$, s.t.
\begin{align}
 &&\partial_t u + \vect{w} \cdot \nabla u - \Delta u &= f && \textnormal{in } \Omega(t) \textnormal{ for } t \in (0,T], \label{strongformproblem} &&&&&\\
 &&\nabla u \cdot \vect{n}_{\partial \Omega} &= 0 && \textnormal{on } \partial \Omega(t) \textnormal{ for } t \in (0,T], \nonumber & &&&&
\end{align}
with initial data $u(\cdot, 0) = u_0$ in $\Omega(0)$.
Here, $\vect{w}$ represents a divergence-free convection field with $\Vert \vect{w} \Vert_{L^\infty(Q)} \leq c < \infty$
and $f$ sources or sinks for the species. We further assume that the convection field $\vect{w}$ coincides with the domain motion on the spatial boundary, so that on $\Gamma_* :=\bigcup_{t\in (0,T]} \partial\Omega(t) \times \{ t \}$  there is $(\vect{w}, 1)\perp\vect{n}^*$ where $\vect{n}^*$ is the space-time normal to $\Gamma_*$.

In order to formulate a well-posed variational formulation in a space-time setting, we define the Sobolev space of functions with weak derivatives in spatial direction:
\begin{equation}\label{def:Hoz}
  H^{1\!,0}(\Q)\! :=\! \overline{C^1(\overline{\Q})}^{\Vert \cdot \Vert_{\!H^{1\!,0}(\Q)}} \text{ with } \Vert v \Vert_{H^{1,0}(\Q)}^2 \!:=\! \Vert \nabla v \Vert_{L^2(\Q)}^2 \!+\! \Vert v \Vert_{L^2(\Q)}^2, ~ v\!\in\! C^1(\overline{\Q}),
\end{equation}
and denote its dual as \hypertarget{def:Hmoz}{$H^{-1,0}(\Q) := (\Hoz(\Q))^*$}. The weak formulation then reads: 
Find $u \in \Hoz(\Q)$ with $\partial_t u \in \Hmoz(\Q)$ with $u(\cdot,0) = u_0 \in L^2(\Omega(0))$, so that
\begin{equation} \label{eq:weakspacetime}
  \langle \partial_t u, v \rangle + (\vect{w} \cdot \nabla u, v)_{\Q} + ( \nabla u, \nabla v)_{\Q} = \langle f, v \rangle \quad \forall~ v \in \Hoz(\Q),
\end{equation}  
where $\langle \cdot, \cdot \rangle$ denotes the duality pairing between $\Hoz(\Q)$ and $\Hmoz(\Q)$ and $(\cdot,\cdot)_G$ denotes the $L^2$ inner product on a domain $G$.
This weak formulation is well-posed, cf.\ also \cite[Section 10.3]{grossreusken11} for the more involved case of a moving interface problem.
\begin{remark}
Let us stress that although this model is comparably simple, as it is a scalar parabolic equation with very simple boundary conditions and the domain motion is assumed to be known, the techniques that are to be discussed in the remainder of this manuscript can also be applied to more complex problems, e.g.\ vectorial free boundary problems with possibly different boundary conditions. This especially holds for the geometry handling, which is a main feature of this study.
\end{remark}
Next, we introduce three different suitable time-discrete variational formulations of \cref{eq:weakspacetime}. With these, several features of moving domain problems in an Eulerian setting and the considered time discretisation approaches become visible. In this sense this section also serves as a preparation for the fully discrete formulations in \cref{sect:weak_forms}.

\subsection{Discontinuous Galerkin in time formulation}
We assume that the time interval $[0,T]$ is subdivided into time intervals, \hypertarget{def:In}{$I_n = (t_{n-1}, t_n]$}, $n=1, \dots,N$ where $0 = t_0 < t_1 < \dots < t_{n-1} < t_N = T$, and, for ease of presentation, we further assume that all slices have the same length, $\Delta t = t_{n} - t_{n-1} = T/N$.
Let \hypertarget{def:Qn}{$Q^n := \bigcup_{t\in \In{n}} \Omega(t) \times \{t\}$} be the space-time domain corresponding to one time interval $\In{n}$. Further, set
\hypertarget{def:Osn}{$\Omega_{{\mysquare}}^n =  \bigcup_{t \in \In{n}} \Omega(t)$}
and let
\hypertarget{def:Qsn}{$Q_{{\mysquare}}^n = \Osn{n} \times I_n$} be the smallest tensor-product domain containing $\Qn{n}$, see \cref{fig:slabs-TP-domains-sketch} for an illustration.
\begin{figure}[b!]
  \begin{center}
  \vspace*{-0.3cm}
  \tikzexternalexportnextfalse
    \begin{tikzpicture}[yscale=0.8]
       \begin{scope}
	 \draw (0,1) node[left] {$t_{n+1}$};
         \draw (0,0) node[left] {$t_{n}$};
         \draw (0,-1) node[left] {$t_{n-1}$};
         \draw (8.9,0.65) node[left] { \scriptsize \color{red} $\Gamma_*$};
         \draw (3.725,0.65) node[left] { \scriptsize \color{red} $\Gamma_*$};

         \node[] (a) at (2.825,0) {};
         \node[] (b) at (4.2,1) {};
         \node[] (c) at (9.4,1) {};
         \node[] (d) at (8.2,0) {};
         \fill[green!80!gray,opacity=0.25] (a.center) to[out=45,in=-120] (b.center) -- (c.center) to[out=-120,in=45] (d.center) -- cycle;
	
         \node[] (e) at (4.5,1.5) {};
         \node[] (f) at (9.7,1.5) {};
         \fill[green!80!gray,opacity=0.25] (b.center) to[out=60,in=-120] (e.center) -- (f.center) to[out=-120,in=60] (c.center) -- cycle;

         \node[] (g) at (2.425,-1) {};
         \node[] (h) at (7.8,-1) {};
	 \fill[green!80!gray,opacity=0.25] (g.center) to[out=85,in=-135] (a.center) -- (d.center) to[out=-135,in=75] (h.center) -- cycle;
        
         \node[] (i) at (2.25,-1.7) {};
         \node[] (j) at (7.6,-1.7) {};
         \draw[red,thick] (i.center) to[out=55,in=-95] (g.center); 
         \draw[red,thick] (g.center) to[out=85,in=-135] (a.center); 
         \draw[red,thick] (a.center) to[out=45,in=-120] (b.center); 
         \draw[red,thick] (b.center) to[out=60,in=-120] (e.center); 
	 \draw[red,thick] (j.center) to[out=65,in =-105]  (h.center); 
         \draw[red,thick] (h.center) to[out=75,in=-135] (d.center); 
         \draw[red,thick] (d.center) to[out=45,in=-120] (c.center); 
         \draw[red,thick] (c.center) to[out=60,in=-120] (f.center); 

         \fill[green!80!gray,opacity=0.25] (i.center) to[out=55,in=-95] (g.center) -- (h.center) to[out=-105,in=75] (j.center) -- cycle;

         \node[] (tl) at (0,-1.7) {};
         \node[] (tu) at (0,1.5) {};
	       \fill[gray!50,opacity=0.25] (tl.center) -- (i.center) to[out=55,in=-95] (g.center) to[out=85,in=-135] (a.center) to[out=45,in=-120] (b.center) to[out=60,in=-120] (e.center) -- (tu.center) -- cycle;
        
        \node[] (tlr) at (11.1,-1.7) {};
        \node[] (tur) at (11.1,1.5) {};
	\fill[gray!50,opacity=0.25] (tlr.center) -- (j.center) to[out=65,in=-105] (h.center) to[out=75,in=-135] (d.center) to[out=45,in=-120] (c.center) to[out=60,in=-120] (f.center) -- (tur.center) -- cycle;
       
	\draw[blue,thick] (g.center) -- (d |- 52, 52 |- h) -- (d.center) -- (g |- 52, 52 |- a) -- cycle;
	\draw[] (2.0,-0.5) node[] { \textcolor{blue}{ \scriptsize $Q^{n}_{\mysquare}$}};
	\draw[blue,thick] (a.center) -- (c |- 52, 52 |- d) -- (c.center) -- (a |- 52, 52 |- b)  -- cycle; 
	\draw[] (2.35,0.5) node[] { \textcolor{blue}{ \scriptsize $Q^{n+1}_{\mysquare}$}};

         \draw [thick,decoration={brace,mirror,raise=0.125cm},
           decorate] (a.center) --  (c |- 52, 52 |- d) 
	       node [pos=0.5,anchor=north,yshift=-0.15cm] { \scriptsize \textcolor{black}{$ \Osn{n} $}};
         \draw [thick,decoration={brace,mirror,raise=0.125cm},
           decorate] (g.center) --  (d |- 52, 52 |- h) 
	       node [pos=0.5,anchor=north,yshift=-0.15cm] { \scriptsize \textcolor{black}{$ \Osn{n-1} $}};
         \draw [thick,decoration={brace,raise=0.125cm},
           decorate] (g |- 52, 52 |- a) --   (c |- 52, 52 |- d)  
	       node [pos=0.5,anchor=north,yshift=0.7cm] { \scriptsize \textcolor{black}{$ \Osnp{n} $}};

        \draw[thick,dotted] (-0.1,1)  -- (a |- 52, 52 |- b) node(xline)[right]{};
	\draw[thick,dotted] (c.center)  -- (11.1,1.0) node(xline)[right]{};
        \draw[thick,dotted] (-0.1,0)  -- (g |- 52, 52 |- a) node(xline)[right]{};
	\draw[thick,dotted] (c |- 52, 52 |- d)  -- (11.1,0) node(xline)[right]{};
        \draw[thick,dotted] (-0.1,-1)  -- (g.center) node(xline)[right]{};
        \draw[thick,dotted] (d |- 52, 52 |- h) -- (11.1,-1) node(xline)[right]{};
	\draw[->,very thick,>=stealth'] (0,-1.8)  -- (0,1.6) node(yline)[above]{$t$};
      
       \end{scope} 
     \end{tikzpicture}
  \end{center}
  \vspace*{-0.3cm}
  \caption{Schematic illustration of (tensor-product) space-time domains on the time slabs.}
 \label{fig:slabs-TP-domains-sketch}
\end{figure}
We define the ansatz and trial space so that on each time interval the functions are of tensor-product form and polynomial in time with $\ktsol \in \mathbb{N}_0$
\begin{equation}\label{def:VN}
  W^n := H^1(\Osn{n}) \otimes \mathcal{P}^{\ktsol}(\In{n}).
\end{equation}
Here and in the following, we denote by $\mathcal{P}^k(E)$ the set of polynomials on the geometric entity $E$ of order smaller or equal $k$.
The space-time formulation is obtained from reducing \eqref{eq:weakspacetime} to one time slab 
and adding an \emph{upwind} stabilisation in time.
For one time interval $n\in\{1,..,N\}$ it takes the form: Find $u \in W^n$, s.t. for all $v\in W^n$ there holds
\begin{align}
   (\partial_t u \!+\! \vect{w} \!\cdot\! \nabla u, v)_{\Qn{n}} \!\!+\! ( \nabla u, \!\nabla v)_{\Qn{n}} \!\!+\! (\!u_+^{n-1}\!\!,\! v_+^{n-1})_{\Omega(t_{n\!-1}\!)} \!\!
   =\! ( f, v )_{\Qn{n}} \!\!+\! (u_-^{n-1}\!\!,\! v_+^{n-1}\!)_{\Omega(t_{n\!-1}\!)},\!\!\!\!\!\!\! \label{eq:DGweakspacetime}
\end{align}
where $u_+^{n-1}(\cdot) = \lim_{t \searrow t^{n-1}} u(\cdot,t)$, and similarly for $v$. Here, $u_-^{n-1}$ takes the role of the initial data for the current time slab $n$. For $n=1$ it is simply $u_0$, the given initial data, for $n>1$ it is $\lim_{t \nearrow t^{n-1}} w(\cdot,t)$ for $w\in W^{n-1}$ the discrete solution from the previous time step.
Subsequent solution of \eqref{eq:DGweakspacetime} for $n=1,..,N$ yields the global (discontinuous-in-time) space-time solution. 
\begin{remark}
With similar techniques as in \cite{LR_SINUM_2013}\footnote{where also a spatial discretisation is considered for a slightly more involved problem} it is fairly easy to show unique solvability and some basic error estimates in the $H^{1,0}$ norm. We are however not aware of any analysis providing superconvergence results also for the $L^2(\Omega(T))$-norm as known from the stationary domain case (see e.g.\ \cite{T97}). 
 \end{remark}

\subsection{A Continuous Galerkin in time formulation} \label{CG_cont_intro}
From the case of a stationary domain $\Omega$ it is known that a Continuous Galerkin (CG) (in time) trial space can be used while keeping a time stepping structure if a discontinuous (in time) test space (of one degree less) is used, see e.g.\ \cite{AM89}. In this case within a time slab one \dof~ (in time) corresponds to an $H^1(\Omega)$ function known from the previous time slab. Hence, only $k_t$ (instead of $k_t+1$) unknown $H^1(\Omega)$ functions remain.

Difficulties arise if one wants to apply this strategy to the case of a moving domain. To illustrate this, let us consider a time slab $\Qn{n}$, $n > 1$. To solve on $\Qn{n}$ we want to use the tensor-product space $H^1(\Osn{n}) \otimes \mathcal{P}^{\ktsol}(\In{n})$ again. However, in the previous time slab we only solved for a solution in $H^1(\Osn{n-1}) \otimes \mathcal{P}^{k_t}(\In{n-1})$ so that the initial data for the current time slab $n$ is only defined on $\Osn{n-1}\supset\Omega(t^{n-1})$, but we can have $\Osn{n} \nsubseteq \Osn{n-1}$, i.e.\ the initial data is insufficient to remove the first \dof~(in time). 
The geometrical situation is sketched in \cref{fig:slabs-TP-domains-sketch}.
To overcome this issue we combine the solution step on every time slab $n$ with a Sobolev extension \hypertarget{def:ext}{$\mathcal{E}$} of the solution $u(t^{n}) \in H^1(\Osn{n-1})$ on $H^1(\Osnp{n})$ with $\Osnp{n}$ s.t. $\Osn{n},\Osn{n-1}\subset \hypertarget{def:Osnp}{\Osnp{n}}$. Thereby we can match up the \hypertarget{def:dof}{degrees of freedom (\dofs)} (in time) for trial and test space again, which (on each time slab $n$) we choose as
\begin{align}
  U^n &= U_0^n + \ext u_-^{n-1} \cdot \phi_0(t) ,~~  U_0^n := H^1(\Osn{n}) \otimes \mathcal{P}_0^{\ktsol}(\In{n}), \label{def:VNWNCG}~~ 
   V^n = H^1(\Osn{n}) \otimes \mathcal{P}^{\ktsol-1}(\In{n}), \nonumber
\end{align}
where $u_-^{n-1}\in H^1(\Osn{n-1})$ is the solution from the previous time slab (or given initial data), $\phi_0$ is a function in $\mathcal{P}^{\ktsol}(\In{n} )$ with $\phi_0(t^{n-1}) = 1$ and $\mathcal{P}_0^{\ktsol}(\In{n} )$ is the space of polynomials up to degree $\ktsol$ with value zero at $t^{n-1}$.
The space-time formulation on time slab $n\in\{1,..,N\}$ then takes the form: Find $u \in U^n$, s.t.
\begin{align} \label{eq:CGweakspacetime}
  (\partial_t u + \vect{w} \cdot \nabla u, v)_{\Qn{n}} + ( \nabla u, \nabla v)_{\Qn{n}} = ( f, v )_{\Qn{n}} \quad \forall~ v \in V^n. 
\end{align}  
Again, subsequent solution for $n=1,..,N$ yields the global space-time solution.

\subsection{Galerkin Collocation formulations with higher regularity}
In the recent decade space-time finite element methods of Petrov-Galerkin-type with higher regularity became popular \cite{BMW19,anselmann2020gccpostprocwave,BM21,AB21}. The idea is to impose a higher regularity $\ktreg \geq 1$ in time on the trial space which allows to reduce the number of temporal unknowns per time slab. Straightforward generalizations of the previously introduced scheme are however not necessarily stable which leads to the following modification: Instead of using test functions of degree $\ktsol -\ktreg  - 1$ to match the \dofs~(in time), a smaller degree $\kttest$ is used and the remaining \dofs~are used to impose the PDE (in weak form) at $\ktcol$ fixed time instances (the collocation points) yielding a Galerkin-Collocation (GCC) method.

With $U_0^n := H^1(\Osn{n}) \otimes \mathcal{P}_{(\ktreg),0}^{\ktsol}(\In{n})$ we choose the trial and test spaces
\begin{align}
  U^n = U_0^n + \sum_{l=0}^{\ktreg} \ext u_{(l),-}^{n-1} \cdot \phi_l(t), \quad \label{def:VNWNCG2} ~ 
   V^n = H^1(\Osn{n}) \otimes \mathcal{P}^{\kttest}(\In{n}), 
\end{align}
where $u_{(l),-}^n\in H^1(\Osn{n-1}), l=0,..,\ktreg$ are given from the previous time slab (or given initial data), $\phi_l$ is a function in $\mathcal{P}^{\ktsol}(\In{n} )$ with $\phi_l^{(m)}(t^{n-1}) = \delta_{l,m}$, $l,m = 0,..,\ktreg$ and $\mathcal{P}_{(\ktreg),0}^{\ktsol}(\In{n} )$ is the space of polynomials up to degree $k_t$ which have vanishing derivatives up to degree $\ktreg$ at $t^{n-1}$.
The space-time formulation for time slab $n\in\{1,..,N\}$ takes the form: Find $u \in U^n$, s.t. for $l=1,..,\ktcol = \ktsol - \ktreg - \kttest - 1$
\begin{subequations}
\begin{align} \label{eq:CG2weakspacetime}
  (\partial_t u + \vect{w} \cdot \nabla u, v)_{\Qn{n}} + ( \nabla u, \nabla v)_{\Qn{n}} &= ( f, v )_{\Qn{n}} \quad \forall~ v \in V^n,\\ 
  (\partial_t u + \vect{w} \cdot \nabla u, v)_{\Omega(t_l^n)} + ( \nabla u, \nabla v)_{\Omega(t_l^n)} &= ( f, v )_{\Omega(t_l^n)} \quad \forall~ v \in H^1(\Osn{n}), 
\end{align}  
\end{subequations}
where $t_l^n\in\In{n}$, $l=1,..,\ktcol$ are the collocation points.
Note that for $\ktreg=\ktcol=0$ we recover the CG method from \Cref{CG_cont_intro}.

\section{Geometry handling} \label{sect:geom_handling}
Before we consider temporally and spatially discrete versions of the formulations in the previous section, we discuss the discrete geometry handling in the case of an \emph{unfitted} geometry. In this setting the geometry description is separated from the computational mesh. \hypertarget{def:tOmega}{More specifically, we assume that we are given a level set function on a background domain $\tOmega$}. The level set function $\phi : \tOmega \times [0,T] \to \mathbb{R}$ implicitly describes the space-time geometry by
\begin{equation}
 \Q = \{ (x,t) \in \tOmega \times [0,T] \mid \phi(x,t) < 0 \}.
\end{equation}
Due to the implicit nature of the geometry description it has no explicit parametrisation. It is hence not immediately clear how to realize numerical integration on $\Q$ or $\Qn{n}$ or suitable approximations. In this section we introduce an approach that allows for higher order accuracy in space-time while still guaranteeing positive quadrature weights. This approach relies on the idea of isoparametric unfitted FEM that has been established for the stationary case in \cite{L_CMAME_2016}. After a brief repetition of these techniques we extend the approach to the space-time setting exploiting tensor-product structure.


\subsection{Isoparametric mappings for stationary level set domains}
\begin{figure}
  \begin{center}
    \includegraphics[width=0.65\textwidth]{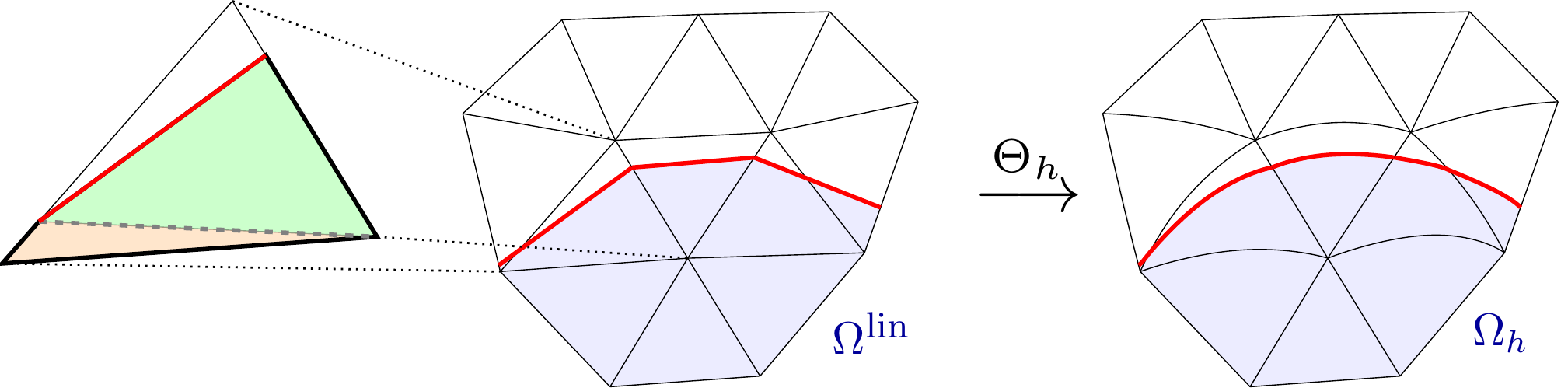}
  \end{center}
  \vspace*{-0.2cm}
  \caption{Sketch of a simple subdivision that allows for setting up quadrature on $T \cap \Omega^{\text{lin}}$ (left), the discrete geometry $\Omega^{\text{lin}}$ (middle) and the mapped geometry $\Omega_h$ obtained from the mapping $\Theta_h$ (right); stationary case.}
  \vspace*{-0.25cm} \label{fig:isossketch}
\end{figure}

In the following, we summarise the method introduced in \cite{L_CMAME_2016} which prepares the techniques and notations used later on in the space-time setting. 
We assume a domain $\Omega$ with smooth boundary to be given as well as a background domain $\tOmega \supset \Omega$ with a shape regular triangulation \hypertarget{def:Th}{$\mathcal{T}_h$}. The domain $\Omega$ is described by a (within the scope of this subsection only spatial) level set function $\phi: \Omega \to \mathbb{R}$, such that $\Omega = \{ x \in \tOmega \, | \, \phi(x) < 0 \}$.
%
%
The geometry description in terms of $\phi$ is only implicit and hence does not offer a parametrisation that can be directly used for setting up numerical integration routines. To enable the use of robust numerical integration routines, we define a piecewise linear approximation of $\phi$ in two steps:
\begin{enumerate}
 \item We introduce a polynomial approximation of higher order $\phi_h \in \Vhks{\khlset}$, where \hypertarget{def:Vhks}{$V_h^{\khlset} := \{ v \in H^1(\tOmega) \mid v|_T \in \mathcal{P}^{\khlset}(T) ~ \forall ~ T\in \Th \}$}. This will lead to an approximation $\Omega^{h} = \{ x \in \tOmega \, | \, \phi_h(x) < 0 \}$, where $\mathrm{dist}(\partial \Omega, \partial \Omega^h) \lesssim h^{\khlset+1}$. As with $\phi$, this description of the approximated domain is still implicit so that it cannot be exploited for numerical integration directly. However, in the following we will construct an approximation which is \textit{as good as} $\Omega^h$ and has an explicit parametrisation.
 \item Next, we define by \hypertarget{def:philins}{$\phi^{\text{lin}} \in \Vhks{1}$ a continuous elementwise {\color{gray}(multi-)}linear approximation of $\phi_h$} leading to the geometry approximation \hypertarget{def:Olins}{$\Omega^{\text{lin}} = \{ x \in \tOmega \, | \, \philins(x) < 0 \}$.}
\end{enumerate}
On simplex meshes $\Olins$ is polygonal and each cut element can easily be decomposed into uncut elements on which numerical integration can easily be applied, cf.\ the sketch on the left of \cref{fig:isossketch}. Also on non-simplex meshes setting up numerical integration is simplified for $\Olins$ compared to $\Omega$ (or $\Omega^h$), cf. \cite{HL_ENUMATH_2019}. 

This simplification in the geometry handling comes at the price of accuracy. By construction $\Olins$ is only a second order approximation to $\Omega$.
To improve the accuracy of the approximation we introduce a vectorial finite element function that serves as a mesh transformation $\Theta_h \colon \tOmega \to \tOmega, \Theta_h \in [\Vhks{\khlset}]^d$. The transformation $\Theta_h$ is constructed and applied on $\Th$ such that $\Omega_h := \Theta_h(\Olins)$ is a higher order approximation to $\Omega$, cf. \cref{fig:isossketch} for a sketch. Note that for technical reasons discussed in \cite{L_CMAME_2016}, $\Omega_h$ is not exactly equal to $\Omega^{h}$ stemming from $\phi_h$. However, their difference is asymptotically small, $\mathrm{dist}(\partial \Omega_h, \partial \Omega^h) \lesssim h^{\khlset+1}$, i.e.\ small enough to yield 
$ \mathrm{dist}(\partial \Omega_h, \partial \Omega) \lesssim h^{\khlset+1}$.

Problems of numerical integration on $\Omega_h$ can then be reformulated as numerical integration problems on $\Olins$ with transformation factors involving $\Theta_h$. Corresponding to the deformation of the mesh the involved finite element spaces are also mapped as usual in isoparametric finite element methods. 

For the construction, initially a preliminary mapping is set up on a subset of elements in $\Th$ which is typically chosen as the domain of cut elements $\mathcal{T}_h^{\Gamma} := \{ T \in \Th \mid \operatorname{meas}_{d} (T \cap \partial \Olins) > 0 \}$.
In a \emph{blending step}, a transition to the identity mapping is applied on all neighbours to cut elements. On all remaining elements the mapping is set to the identity.
For the construction of the mapping we introduce the map 
\begin{equation} \label{def:C}  
 \mathcal{C} : C^1(\Th) \times \Vhks{1} \times \operatorname{Pow}(\Th) \to [\Vhks{\khlset}]^d
\end{equation}  
so that $\Theta_h = \mathcal{C}(\phi_h,\philins,\mathcal{T}_h^{\Gamma})$, where we refer to the power set of $\Th$ by $\operatorname{Pow}(\Th)$.
We refer the reader to \cite{L_CMAME_2016} and \cite{LR_IMAJNA_2018} for a detailed account on the construction and analysis of the mapping for the stationary case.



\subsection{Isoparametric mappings for space-time level set domains}

We want to transfer the previous ideas to the space-time setting exploiting a tensor-product structure. 
We consider a time interval $\In{n} = (t_{n-1}, t_n]$ and define the space-time tensor-product mesh of prisms by extrusion $\{ T \times \In{n} \mid T \in \Th \}$.
In this setting, we assume that our space-time domain $\Qn{n}$ is described by a level set function $\phi: \tOmega \times \In{n} \to \mathbb{R}$ which is approximated by a space-time level set function in a tensor-product finite element space, $\phi_h \in \Vhks{\khlset} \otimes \mathcal{P}^{ \ktlset} $, where  $\ktlset$ is the order of the temporal accuracy. In the following we also use $\klset = (\khlset,\ktlset)$ to summarise the orders of spatial and temporal geometrical accuracy and write $\klset=r$ if $\khlset\!=\!\ktlset\!=\!r$.
Additionally we make the following assumption.
\begin{assumption} \label{continterprem}
  The approximated level set function 
 $\phi_h : \tOmega \times \In{n} \to \mathbb{R}$ is continuous across all time slabs and hence ensures a connected approximated space-time geometry, cf. \cref{Cuttop} for a sketch.
\end{assumption}
 We are now going to construct a space-time mapping $\Theta_h^n: \tOmega \times \In{n} \to \tOmega$ based on the tensor-product structure of $\phi_h$ and hence look for a space-time transformation $\Theta_h^n \in [\Vhks{\khlset} \otimes \mathcal{P}^{\ktlset}]^d$.
Let $\{\ell_0,\dots,\ell_{ \ktlset } \}$ be a basis of $\mathcal{P}^{\ktlset }(\In{n})$, then we can write 
\begin{equation}\label{eq:phihtp}  
  \phi_h (x,t) = {\sum}_{i=0}^{\ktlset} \ell_i(t) \cdot \phi_h^i(x), \quad \text{for functions } \phi_h^i \in \Vhks{\khlset}.
 \end{equation}
Similarly to the stationary case discussed in the previous section, we will use a level set function that is elementwise (multi-)linear in space to set up a reference configuration. Note that we keep the order in time fixed, i.e.\ higher order. 
Let $I_h^1: C(\tOmega) \to V_h^1$ be the nodal (spatial) interpolation operator. We define 
\begin{equation}
  \phi^{\text{lin}} (x,t) := (I_h^1 \phi_h(\cdot,t))(x) ~~ \Longleftrightarrow ~~
 \phi^{\text{lin}} (x,t) = {\sum}_{i=0}^{ \ktlset } \ell_i(t) \cdot I_h^1 \phi_h^i(x).
 \label{philindef}
\end{equation}
%
%
We develop the transformation $\Theta_h^n$ in the same tensor-product structure,
\begin{equation}
  \Theta_h^n(x,t) = {\sum}_{i=0}^{\ktlset} \ell_i(t) \cdot \Theta_{h,i}(x)
\label{eq:collect_sliced_mappings}
\end{equation}
and construct each spatial transformation $\Theta_{h,i}(x)$ based on $\phi_h^i(x)$ and $I_h^1 \phi_h^i (x)$ as in the stationary case using the map $\mathcal{C}$. 
It is crucial to make sure that the set of active elements is the same within one time slab, cf. \cref{thirdargumentofCremark} below. We hence define
\begin{equation} \label{eq:new_st_isoparam_regions}
 \mathcal{T}_{h,n}^{\Gamma} := \{ T \in \Th \mid \operatorname{meas}_{d} (T \cap  \Olins(t)) > 0 \textnormal{ for some } t \in \In{n} \},
\end{equation}  
and choose $\Theta_{h,i} = \mathcal{C}(\phi_h^i(x),I_h^1 \phi_h^i (x), \mathcal{T}_{h,n}^{\Gamma}),~i=0,\dots,\ktlset $.
%
%
%
\begin{figure}
  \begin{center}
  \tikzexternalexportnextfalse
 \begin{tikzpicture}[xscale=0.8]
  \node (nn) at (-2,0) {$\phi(x,t)$};
  \node (n0) at (0,0) {$\phi_h(x,t)$};
  \node (n1) at (2.5,1) {$\phi_h^0(x)$};
  \node (n2) at (2.5,0.333333) {$\phi_h^1(x)$};
  \node (n3) at (2.5,-0.333333) {\dots};
  \node (n4) at (2.5,-1) {$\phi_h^{\ktlset}(x)$};
  \draw[->] (n0) -- (n1.west); \draw[->] (n0) -- (n2.west); \draw[->] (n0) -- (n4.west);
  \node at (1.35,1.5) {restrict,\eqref{eq:phihtp}};
  \node at (-1,0.75) {approximate};
  \node (n5) at (7.3,1) {$\Theta_{h,0} =\mathcal{C} (\phi_h^0(x), I_{h}^1 \phi_h^0(x),\ThnG{n} )$};
  \node (n6) at (7.3,0.333333) {$\Theta_{h,1}  =\mathcal{C} (\phi_h^1(x), I_{h}^1 \phi_h^1(x), \ThnG{n})$};
  \node (n7) at (7.3,-0.333333) {\dots};
  \node (n8) at (7.3,-1) {$\Theta_{h,\ktlset} =  \mathcal{C} (\phi_h^{\ktlset} (x), I_{h}^1 \phi_h^{\ktlset}(x), \ThnG{n})$};
  \draw[->, dashed] (n1.east) -- (n5.west); \draw[->, dashed] (n2.east) -- (n6.west); \draw[->, dashed] (n4.east) -- (n8.west);
  \node[text width=4cm, align=center] at (4.75,1.75) {apply $\mathcal{C}$, \eqref{def:C} };
  \node (n9) at (12,0) {$\Theta_{h}^n$};
  \draw[->, dotted] (n5.east) -- (n9); \draw[->, dotted] (n6.east) -- (n9); \draw[->, dotted] (n8.east) -- (n9);
  \node at (11.5,1.75) {sum up, \cref{eq:collect_sliced_mappings}};
  \draw[->] (nn.east) -- (n0.west);
 \end{tikzpicture}
\end{center}
\vspace{-0.3cm}
 \caption{Schematic illustration of the construction of the space-time isoparametric mapping.}
 \label{space_time_mapping}
 \vspace{-0.5cm}
\end{figure}
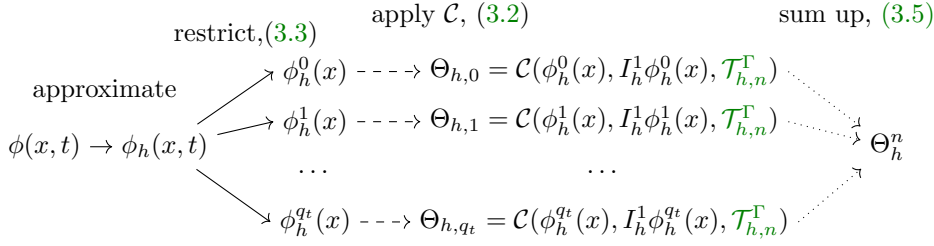
The overall procedure is applied to each time slice $I_n$, leading to a mapping $\Theta_h^n\colon \tOmega \times I_n \to \tOmega$ for each $I_n$ respectively. 
A schematic illustration is provided in \cref{space_time_mapping}. Moreover, by $\Theta_h^{n,\text{st}}$ we denote the function $(x,t) \mapsto (\Theta_h^n(x,t), t)$.

Gathering the mappings of all time slabs, we denote by $\Theta_h$ the mapping defined on $\tOmega \times [0,T]$ which takes the value of $\Theta_h^n$ for $t \in [t_{n-1},t_n)$.\footnote{Inside a time slice interval, $\Thehst$ is uniquely defined, but not on the slice boundaries $t_{n-1}$.
 That we opt for the value from $I_n$ for these time points is a matter of notational convenience.} Moreover, by \hypertarget{def:Thehst}{$\Theta_h^{\text{st}}$} we denote the function mapping each $(x,t)$ in $\tOmega \times [0,T]$ to $( \Theta_h (x,t), t)$.
In terms of $\phi^{\text{lin}}$, we furthermore introduce the discrete regions \hypertarget{def:Omlin}{$\Omega^{\text{lin}}(t) = \{ x \in \tilde \Omega \, | \, \phi^{\text{lin}}(x,t) < 0 \}$}, and
\begin{align}
  \Omega^h(t) := \Theta_h (\Omlin{t}, t), \quad
  \hypertarget{def:Qlinn}{Q^{\text{lin}, n} = {\bigcup}_{t \in I_n} \Omega^{\text{lin}}(t) \times \{ t \}},
 \quad \hypertarget{def:Qhn}{Q^{h,n} = \Theta_h^{\text{st}}(Q^{\text{lin},n})} \label{eq:introQhn}.
\end{align}
The constructed mapping will in general not be continuous in time due to different active elements on different time slabs. However, $\Omlin{t}$ and  $\Omega^h(t)$ are continuous. 

\begin{figure}
   \centering
	\subfloat[Domain evolution.] 
{\label{kiteillu}
  \begin{tikzpicture}[xscale=1.4, yscale=1.4, meshtrig/.style={fill=Set1-B!10!white, draw=Set1-B},
     backgrtrig/.style={draw=gray!50!white, thin}, reg_Fh/.style={draw=Set1-D, thick, dashed},
     Gamma_hbeg/.style={draw=Set1-A, thick, {Circle[length=3pt, width=3pt, open, color=Set1-A]}-{Circle[length=3pt, width=3pt, open, color=Set1-A]}, shorten <=-1.5pt, shorten >=-1.5pt},
     Gamma_hend/.style={draw=Set1-E, thick, {Circle[length=3pt, width=3pt, open, color=Set1-E]}-{Circle[length=3pt, width=3pt, open, color=Set1-E]}, shorten <=-1.5pt, shorten >=-1.5pt}]
\begin{scope}
	 \input{tikz_mesh}

  \draw[Set1-A, domain=-3.141:3.141,smooth,variable=\t]
plot ({ cos(\t r)+(1 - sin(\t r)^2)*0.00 },{sin(\t r)});
  \draw[Set1-E, domain=-3.141:3.141,smooth,variable=\t]
plot ({ cos(\t r)+(1 - sin(\t r)^2)*0.125 },{sin(\t r)});
  \end{scope}
 \begin{scope}[yshift = -2.25cm]
 \input{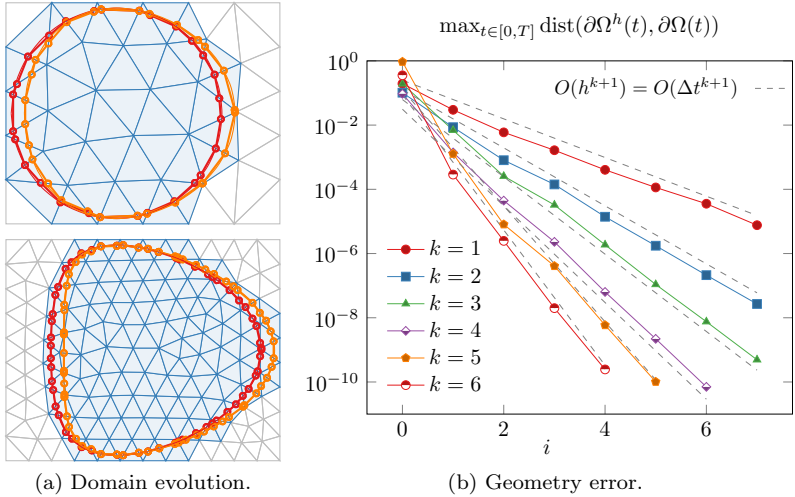}

  \draw[Set1-A,domain=-3.141:3.141,smooth,variable=\t]
plot ({ cos(\t r)+(1 - sin(\t r)^2)*(0.50 -0.125) },{sin(\t r)});
   \draw[Set1-E,domain=-3.141:3.141,smooth,variable=\t]
plot ({ cos(\t r)+(1 - sin(\t r)^2)*(0.50) },{sin(\t r)});
 \end{scope}
    \end{tikzpicture}
	  }
\subfloat[Geometry error.]
  { \label{plot_kite_geom_approx}%
  \begin{tikzpicture}[scale=0.825]
    \begin{semilogyaxis}[ xlabel=$i$, title={$\max_{t \in [0,T]} \mathrm{dist}(\partial \Omega^h(t), \partial \Omega(t))$}, legend entries ={ $k=1$, $k=2$, $k=3$, $k=4$, $k=5$, $k=6$}, legend style={anchor=north,legend columns=1, draw=none}, legend pos = south west, x label style={at={(axis description cs:0.45,0.0)},anchor=east},
     x tick label style={yshift=0cm,xshift=0.0em},
      y label style={at={(axis description cs:0.25,0.85)},anchor=east},
      ymax = 1e0,
      ymin = 1e-11
     ]
     \addplot table [x index = 0, y index =5] {num_exp/out/conv_kite_DG_ks1_kt1_both_nref8_gamma0.05_sm2.dat};
     \addplot table [x index = 0, y index =5] {num_exp/out/conv_kite_DG_ks2_kt2_both_nref8_gamma0.05_sm2.dat};
     \addplot table [x index = 0, y index =5] {num_exp/out/conv_kite_DG_ks3_kt3_both_nref8_gamma0.05_sm2.dat};
     \addplot table [x index = 0, y index =5] {num_exp/out/conv_kite_DG_ks4_kt4_both_nref7_gamma0.05_sm2.dat};
     \addplot table [x index = 0, y index =5] {num_exp/out/conv_kite_DG_ks5_kt5_both_nref6_gamma0.05_sm2.dat};
     \addplot table [x index = 0, y index =5] {num_exp/out/conv_kite_DG_ks6_kt6_both_nref5_gamma0.05_sm2.dat};
     \addplot[gray, dashed, domain=0:7] {(1/2^(x+1)))^2};
     \addplot[gray, dashed, domain=0:7] {(1/2^(x+1)))^3};
     \addplot[gray, dashed, domain=0:7] {(1/2^(x+1)))^4};
     \addplot[gray, dashed, domain=0:6] {(1/2^(x+1)))^5};
     \addplot[gray, dashed, domain=0:5] {(1/2^(x+0.5)))^6};
     \addplot[gray, dashed, domain=0:4] {(1/2^(x+0.5)))^7};
    \end{semilogyaxis}
                  \node[scale=0.75] at (4.5,5.25) {$O(h^{k+1})= O(\Delta t^{k+1})$};
     \draw[scale=0.75, gray, dash=on 2.25pt off 2.25pt phase 0pt, line width=0.4*0.75pt] (6.25/0.75,5.25/0.75) -- (6.8/0.75,5.25/0.75);
   \end{tikzpicture}
}
\vspace*{-0.2cm}
\caption{Example of a deforming kite geometry: Subfigure
(a) displays the evolution of the geometry and discrete regions on the intervals $I_1 = [0, \Delta t]$ (top) and $I_{N} = [T - \Delta t,T]$ (bottom), respectively.
The quality of the discrete geometry approximation is investigated in (b).
}
\label{fig:kite-geom-test}
\vspace*{-0.5cm}
\end{figure}

In numerical experiments (cf. \cref{fig:kite-geom-test}, \cref{sect:numexp}), we confirm that $\Qhn$ satisfies the following approximation error bound
\begin{equation}
 \max_{t \in I_n} \mathrm{dist}(\partial \Omega^h(t), \partial \Omega(t)) \lesssim h^{\khlset + 1} + \Delta t^{\ktlset +1}. \label{discretestregionapproxprop}
\end{equation}
\vspace{-0.6cm}
\begin{remark} \label{thirdargumentofCremark}
 A subtle, but very important detail is the decision on the third argument of $\mathcal{C}$, the subset of actively deformed elements.
A seemingly straightforward approach would associate to every $i \in \{0,\dots,\ktlset \}$ a set of active elements $\mathcal{T}_h^i$ according to the cut elements associated to $I_h^1 \phi_h^i$. This would easily allow to make sure that the constructed transformation is also continuous in time\footnote{With $\ell_i(t_{n-1}) = \delta_{i0}$ and $\ell_i(t_{n}) = \delta_{i\ktlset}$, $\Theta_{h,i} := \mathcal{C}(\phi_h^i,I_h^1 \phi_h^i,\mathcal{T}_h^i)$ would yield $\Theta_{h,\ktlset}^{n-1}=\Theta_{h,0}^{n}$.}. However, we observed in numerical experiments that this method would not yield the optimal approximation quality \cref{discretestregionapproxprop} for higher orders.
In the remainder we consider the previously presented construction with a discontinuous-in-time mapping. The question if 
a similar construction yielding a continuous deformation can be found which preserves the desired accuracy is left open for future research.
\end{remark}
A proper treatment of the discontinuity in the mesh deformation is discussed next.
\subsection{Discontinuous-in-time deformations and mesh transfer operations}
\label{subsect_discont_def}
\begin{figure} \centering
  \def\h{1}
  \begin{tikzpicture}[xscale=0.8, yscale=1]
  \draw[gray] (2,0*\h) grid[xstep=1, ystep=\h] (8,2*\h);
  \begin{scope}[xscale=-1, xshift=-10cm]
  \draw[thick, Set1-A] (4.5,0*\h) -- (4.7,1*\h) -- (6.1,2*\h);
  \draw[thick, dashed] (4.5,0*\h)
                .. controls (4.6,0.5*\h) .. (4.7,1*\h)
                .. controls (4.8,1.5*\h) and (5.5,1*\h) .. (6.1,2*\h);
  \fill[Set1-B, opacity=0.2] (4,0*\h) rectangle (5,1*\h);
  \fill[gray, opacity=0.4] (4,0*\h) rectangle (3,2*\h);
  \fill[gray, opacity=0.4] (5,0*\h) rectangle (6,1*\h);
  \fill[gray, opacity=0.4] (7,1*\h) rectangle (8,2*\h);
  \fill[Set1-B, opacity=0.2] (7,1*\h) rectangle (4,2*\h);
   \end{scope}
   \node at (1.4,0*\h) {$t_{n-1}$};
   \node at (1.4,1*\h) {$t_{n}$};
   \node at (1.4,2*\h) {$t_{n+1}$};
   \node at (6,1.25*\h) {$\{\phi = 0\}$};
   \node[Set1-A] at (5.5,1.6*\h) {$\{\phi_h = 0\}$};
 \end{tikzpicture} \hspace*{-0.07\textwidth}
 \begin{tikzpicture}[xscale=0.8, yscale=1]
  \draw[gray] (2,0*\h) grid[xstep=1, ystep=\h] (8,2*\h);
  \begin{scope}[xscale=-1, xshift=-10cm]
  \draw[thick, Set1-A] (4.5,0*\h) -- (4.7,1*\h);
  \draw[thick, Set1-A] (5.02,1*\h) -- (6.2,2*\h);
  \draw[thick, dashed] (4.5,0*\h)
                .. controls (4.6,0.5*\h) .. (4.7,1*\h)
                .. controls (4.8,1.5*\h) and (5.5,1*\h) .. (6.1,2*\h);
  \fill[Set1-B, opacity=0.2] (4,0*\h) rectangle (5,1*\h);
  \fill[gray, opacity=0.4] (4,0*\h) rectangle (3,1*\h);
  \fill[gray, opacity=0.4] (5,1*\h) rectangle (4,2*\h);
  \fill[gray, opacity=0.4] (5,0*\h) rectangle (6,1*\h);
  \fill[gray, opacity=0.4] (7,1*\h) rectangle (8,2*\h);
  \fill[Set1-B, opacity=0.2] (7,1*\h) rectangle (5,2*\h);
   \end{scope}
  \vphantom{
  \node at (1.4,0*\h) {$t_{n-1}$};
  \node at (1.4,1*\h) {$t_{n}$};
  \node at (1.4,2*\h) {$t_{n+1}$};
  }
  \node at (6,1.25*\h) {$\{\phi = 0\}$};
   \node[Set1-A] at (5.5,1.6*\h) {$\{\phi_h = 0\}$};
   \node[align=left, text width=1.5cm, anchor=west,scale=0.8] at (9,0.5*\h) {$\mathcal{T}_h^{\Gamma}$};
   \fill[Set1-B, opacity=0.2] (8.5, 0.1*\h) rectangle (9, 0.9*\h);
   \node[align=left, text width=1.2cm, anchor=west,scale=0.8] at (9,1.5*\h) {cont. extension};
   \fill[gray, opacity=0.4] (8.5, 1.1*\h) rectangle (9, 1.9*\h);
 \end{tikzpicture}
 \vspace*{-0.8cm}
 \label{Cuttop}
 \caption{Cut topologies for different interpolations of the same space-time level set geometry. Left: The case of \cref{continterprem}. Right: A more general interpolation is used.}
 \vspace{-0.5cm}
\end{figure}
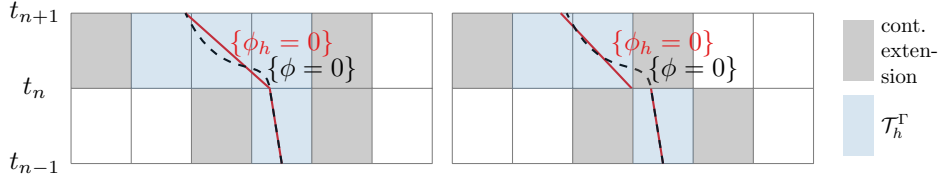

With the previously presented construction of the deformation, for successive time intervals $\In{n} = (t_{n-1}, t_{n}], \In{n+1} = (t_n, t_{n+1}]$, we will obtain deformations $\Theta^n_h(x,t)$ and $\Theta^{n+1}_h(x,t)$, respectively, where
\begin{equation}
 \Theta^{-}_h (t_n) := \Theta^{n}_h (\dots, t_n) \neq \Theta^{n+1}_h (\dots, t_n) =: \Theta^{+}_h (t_n)
\end{equation}
in general. One specific $T \in \mathcal{T}_h$ can be associated to different regions of the isoparametric mapping construction, e.g. $T \in \mathcal{T}_{h,n+1}^{\Gamma}$ but $T \not\in \mathcal{T}_{h,n}^{\Gamma}$, cf. \cref{Cuttop} (left) for an illustration of that situation.

The challenge within this setting stems from the parametrically mapped discrete function spaces. To take the deformation into account, a discrete space-time function defined on $I_n$, restricted to $t_n$, will have the structure $u^- \in \Vhks{\khlset} \circ (\Theta^{-}_h (t_n))^{-1}$, whilst we aim for a function $\Pi_u (u^-) \in \Vhks{\khlset} \circ (\Theta^+_h(t_n))^{-1}$ for a discrete representation of the initial data for $I_{n+1}$ for the CG and GCC method. A standard $L^2$ projection from one deformed mesh to the other would be a valid option.
Within the context of Eulerian time stepping methods for higher order methods in space and time, a similar transfer operation is needed. In \cite[Section 5]{LL_ARXIV_2021} a transfer operation between differently deformed meshes is presented (and thoroughly analysed). It is explicit and requires only essentially element-local operations leading to an operation that is more efficient to evaluate than a standard $L^2$ projection and has negligible costs compared to the remaining operations required to solve the discretisations discussed in the remainder.
We can fortunately exploit this operation also in our (in comparison simpler\footnote{The transfer operation \hypertarget{def:Piu}{$\Pi_u$} introduced in \cite[Section 5]{LL_ARXIV_2021} is designed to translate even between deformations of distinctive time points, whilst we apply it to a setting of two deformations relating to one time point.}) setting as well which we especially did in the numerical examples. For further details of the construction and analysis of the transfer operator we refer to  \cite{LL_ARXIV_2021}.

Concerning the DG method, there is no need for a discrete projection of $u^- \in \Vhks{\khlset} \circ (\Theta^{-}_h (t_n))^{-1}$, as the initial data is only imposed weakly. It is however computationally convenient to have a discrete approximation living on the mesh deformed by $\Theta^+_h(t_n)$. This does not need to be a continuous discrete function, so that we recommend to apply the previously discussed projection only element-wise and skip the final averaging step, called $P_h$ in \cite[Subsection 5.1.3]{LL_ARXIV_2021}. The so-constructed transferred discrete projection is denoted as $\Pi_u^\ast (u^-)$ from now on, so that $\Pi_u (u^-) = P_h \Pi_u^\ast (u^-)$.

\begin{remark}[Accumulation of projection errors from the transfer operation]
In general the application of a transfer operator in each time step of a time stepping schemes can lead to unbounded accumulated projection errors for increasing number of time steps, cf. for instance \cite{brenner2014priori}. However, one important insight in \cite[Section 3.4]{LL_ARXIV_2021} is that the special situation of only very slightly different mesh deformations considered here leads to the following important property: For a fixed element the number of times that the change of the mesh transformation is not Lipschitz-continuous in time is bounded independent of the number of time steps. Furthermore, the projection error accumulation of transfer operations corresponding to Lipschitz-continuous changes in the mesh deformation as well as the projection error accumulation due to the small number of discontinuous changes in the mesh deformation stay bounded independent of the time step size. Hence, overall the accumulation of projection errors stays bounded.
\end{remark}  

%
\section{Fully discrete weak forms} \label{sect:weak_forms}
In this section, we introduce the discrete weak forms of the three methods discussed in the semi-discrete setting in \cref{sect:modprob}. The fully discrete version of the DG method is presented in \Cref{ssection:DG}, the CG method is introduced in \Cref{ssection:CG} and the GCC method in \Cref{ssection:GCC}.
In the following when we distinguish the methods as DG, CG or GCC we refer to the discretisation in time only. In space we always only consider the same continuous finite element space. 

\subsection{Discontinuous Galerkin method} \label{ssection:DG}
\begin{figure}
 \centering
 \begin{tikzpicture}[xscale=1.2, yscale=0.66]
  \draw[gray] (0,0) grid[xstep=1, ystep=1] (10,4);
  \begin{scope}[xscale=-1, xshift=-10cm]
  \draw[thick] (3.5,0) .. controls (4.3, 0.5) and (4.3,0) .. (4.5,1)
                .. controls (4.6,1.5) .. (4.7,2)
                .. controls (4.8,2.5) and (5.5,2) .. (6.1,3)
                .. controls (6.4,3.5) .. (6.9,4);
   \end{scope}
   \draw[gray!50!white, pattern = north east lines, pattern color = gray!50!white] (5,0) rectangle (7,1) -- (5,1) rectangle (6,2) -- (3,2) rectangle (6,3) -- (3,3) rectangle (4,4);
   \draw[dash dot, Set1-B, thick] (7,0) -- (7,1) -- (6,1) -- (6,3) -- (4,3) -- (4,4) -- (10,4) -- (10,0) -- cycle;
   \draw[very thick, dashed, Set1-C] (7,0) -- (7,1) (6,0) -- (6,3) (5,2) -- (5,3) (4,2) -- (4,4);
   \draw[very thick, dashed, Set1-C!50!black] (8,0) -- (8,1) (7,2) -- (7,3) (8,2) -- (8,3);

   \draw[{Bar[]}-, thick] (5.0,1.15) --(6.8,1.15);
   \draw[-{Bar[]}, thick] (8.1,1.15) -- (10, 1.15);
   \node[fill=white]  at (7.5,1.15) {$\mathcal{E}(\Omega^{\text{lin},2})$};

   \node[Set1-A,fill=white] at (8.5,0.4) {$\bigcup_{n=1}^N \mathcal{E}(Q^{\text{lin}, n})$};
   \node[Set1-B,fill=white] at (8.5,3.6) {$\bigcup_{n=1}^N \mathcal{I}(Q^{\text{lin}, n})$};
   \node[gray,fill=white] at (3,1.5) {$\bigcup_{n=1}^N \mathcal{E}(Q^{\text{lin}, n}) \backslash \mathcal{I}(Q^{\text{lin},n})$};
   \node[Set1-C,fill=white] at (6.35,2.5) {$\mathcal{F}^{n}_R$};
   \node[Set1-C!50!black] at (8.85,2.51) {$\mathcal{F}^{n, \text{ext}}_R \backslash \mathcal{F}^{n}_R$};

   \draw[Set1-A, thick] (5,0) -- (5,2) -- (3,2) -- (3,4) -- (10,4) -- (10,0) -- cycle;
   
   \node at (10.25, 0) {$t_0$};
   \node at (10.25, 1) {$t_1$};
   \node at (10.25, 2) {$t_2$};
   \node at (10.25, 3) {$t_3$};
   \node at (10.25, 4) {$t_4$};
   
   \node at (10.1, 1.75) {$t$};
   \draw[{Bar[]}-, thick] (5.34,1.75) --(7.1,1.75);
   \draw[-{Bar[]}, thick] (7.9,1.75) -- (10, 1.75);
   \node at (7.5,1.75) {$\Omega(t)$};

 \end{tikzpicture}\vspace*{-0.43cm}
 \caption{Discrete regions in the definition of the ghost-penalty of the plain DG method.}
 \label{classicalgpsketch} \vspace*{-0.5cm}
 \end{figure}
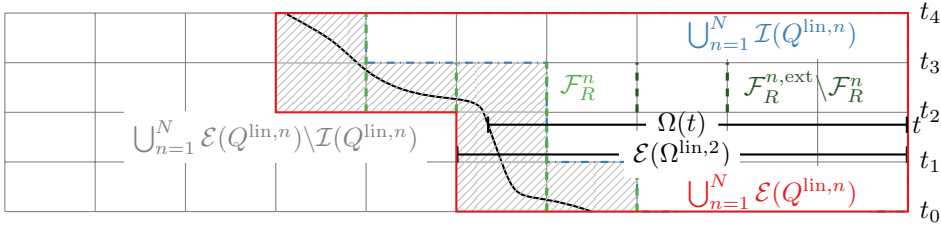

To state the discrete variational formulation of the DG method, we have to introduce several discrete regions. Some of them rely on the spatially linear interpolation of the level set function $\phi^{\text{lin}}$. 
%
First, we define an extended space-time domain, which contains all points in $\Qlinn$ and has tensor-product structure per time slice:
 \begin{align*}
  \hypertarget{def:EOmn}{\mathcal{E}(\Omega^{\text{lin}, n})} :&= \{ x \in T \textnormal{ for some } T \in \mathcal{T}_h \textnormal{ such that } (T \times I_n) \cap \Qlinn \neq \varnothing \}.
 \end{align*}
In addition, we define an interior space-time domain, which is the largest domain with time-slice tensor-product structure which is completely contained in $\Qlinn$:
 \begin{align*}
  \hypertarget{def:IOmn}{\mathcal{I}(\Omega^{\text{lin}, n})} :&= \{ x \in T \textnormal{ for some } T \in \mathcal{T}_h \textnormal{ such that } (T \times I_n) \subseteq \Qlinn \}.
 \end{align*}
 For both spatial regions, we define a space-time counterpart:
\begin{align*}
   \hypertarget{def:EQlinn}{\mathcal{E}(Q^{\text{lin}, n})} &:= \EOmn \times I_n, \quad \hypertarget{def:IQlinn}{\mathcal{I}(Q^{\text{lin},n})} = \IOmn \times I_n. 
\end{align*}
 As usual in unfitted methods, we use a ghost penalty stabilisation to handle ill-posed cut configurations. It is defined almost on the facets between interior and exterior,
  \begin{align*}
  \mathcal{F}^{n}_R = \{F \! \in\! \mathcal{F} \textnormal{ s.t. } \exists T_1 \neq T_2,& T_1 \!\in\! \mathcal{E}(\Omega^{\text{lin}, n}) \backslash \mathcal{I}(\Omega^{\text{lin}, n}), T_2 \!\in\! \mathcal{E}(\Omega^{\text{lin},n}) \textnormal{ with } F = T_1 \! \cap\! T_2 \}.
 \end{align*}~\\[-3ex]
 For the case  when $\mathcal{E}(\Omega^{\text{lin}, n})\setminus \mathcal{I}(\Omega^{\text{lin}, n})$
 has a width that is significantly larger than $h$ we include a slightly larger set of nearby interior facets, $\mathcal{F}^{n, \text{ext}}_R \supseteq \mathcal{F}^{n}_R$, so that the number of facets in the interior is proportional to that in the exterior. We refer to \cite[Assumption 3.1]{heimann20} for a precise explanation. These regions are illustrated for a setting of a one-dimensional domain in \cref{classicalgpsketch}.
 Now, we are in a position to introduce the discrete function spaces.
%
 To this end, first a finite element space for each time slab is defined as
 \begin{equation}
  \hypertarget{def:bgVhnks}{W_h^{n,\ksol}} := \Vhks{\khsol} \otimes \mathcal{P}^{\ktsol}([t_{n-1},t_n]),
 \end{equation}
 where $\ksol = (\khsol,\ktsol)$ is the tuple describing the combination of the polynomial orders in space and time and again for $\khsol\!=\!\ktsol\!=\!r$ we write $\ksol=r$.
 Now we introduce a general notation for cut discrete spaces depending on a geometrical domain $\mathcal{E} \subseteq \tilde \Omega \times [0,T]$:
 \begin{equation}
  \hypertarget{def:Whcut}{W_{h, \text{cut}}^{n,\ksol} (\mathcal{E})} = \{ v \in \bgVhnks{n}{\ksol} \, | \, v \textnormal{ vanishes outside of } \mathcal{E} \textnormal { as element of } \bgVhnks{n}{\ksol} \}.
 \end{equation}
 This means that only those \dofs~that are associated to $\mathcal{E}$ stay \emph{active} whereas all others are set to zero. 
 Note that with this definition of the space only its restriction to $\mathcal{E}$ is uniquely defined\footnote{outside of $\mathcal{E}$ the decay to zero depends on the chosen basis for $\bgVhnks{n}{\ksol}$.}, which however suffices as functions in $\Whcut{n}{\mathcal{E}}$ will not be involved outside of $\mathcal{E}$.
This allows to define the finite element space for the DG method as follows:
\begin{equation}
 W_h^n := \Whcut{\ksol}{\EQlinn} \circ (\Theta_h^{n, \text{st}})^{-1}.
\end{equation}
%
 Using the discrete regions and function spaces introduced above, we are now able to define the first components of the DG method for one time slab:
 \begin{align}
  B^n(u,v)\!&:=\!(\partial_t u \!+\! \vect{w} \cdot \nabla u, v)_{\Qhn} 
  \!+ (\nabla u, \nabla v)_{\Qhn},&
  \hspace*{-0.45cm} f^n(v)\!&:=\!(f,v)_{\Qhn}, \label{DGbilinearforms1} \\
   B_{\text{upw}}^{n}(u,v)\!&:=\!({u}^{n-1}_+\!, v^{n-1}_+)_{\Omega^{h}\!(t_{n-1}\!)}, &
  \hspace*{-0.45cm} f_{\text{upw}}^{n}(v)\!&:=\!( \Pi_u^\ast ({u}^{n-1}_-), v^{n-1}_+)_{\Omega^{h}\!(t_{n-1}\!)}, \nonumber
 \end{align}
where $u^n_{\pm} := \lim_{s \searrow t_n} u(\cdot, s)$
  and $u^{n-1}_-$ is the initial value for the current time slab. For $n=1$ this is the initial value of the global problem whereas for $n>1$ this stems from the solution of the previous time step. The transfer operation $\Pi_u^\ast$ is applied here as $\Theta_h^{\text{st}}$ may be discontinuous across $t^{n-1}$.

 In addition, we define a ghost penalty stabilisation. We use a variant which exploits facet patches, which we call direct version.\footnote{Other versions appear in the literature, such as the original version \cite{burman2010ghost} or the normal-derivative jump version \cite{BURMAN2012328}. See \cite[Remark 6]{preuss18} for a comparison in terms of analysis.} Let $F = \overline{T_1} \cap \overline{T_2}$ be a facet of the undeformed triangulation, $T_1, T_2 \in \mathcal{T}_h$. Then, the facet patch is defined as $\omega_F := T_1 \cup T_2$. Taking into account the mesh deformation, we define also $\omega_F^{h}(t) = \Theta_h(\omega_F,t)$. In addition, we need a \emph{volumetric} jump operation on these curved elements. Let us assume $u = \hat u \circ (\Theta^{\text{st}}_h)^{-1}$, and we are interested in the jump for some point $x$ in the mapped element $T_1$. Then, \vspace{-0.2cm}
 \begin{equation*}
  \jump{u}_{\omega_F^h(t)}|_{\Theta_h(T_1,t)} (x) := u(x) - \mathcal{E}^p(\hat u|_{T_2})( (\mathcal{E}^p \Theta_h(t)|_{T_2}) ^{-1}(x)),
 \end{equation*}
 where the latter term extends the (mapped) polynomial on $T_2$ to $T_1$ based on
 the canonical extension of polynomials from their element $T \in \mathcal{T}_h$ to $\mathbb{R}^d$(denoted by $\mathcal{E}^p$) combined with the involved mesh transformations. Then,
  \begin{align*}
 j_h(\mathcal{F}, t; u,v) &\!:=\! \sum_{F \in \mathcal{F}} \int_{\omega_F^h(t)} \frac{1}{h^2} \jump{u}_{\omega_F^h(t)} \jump{v}_{\omega_F^h(t)}dx,~ J^n(u,v) \!=\! \int_{t_{n-1}}^{t_n} \!\!\tilde \gamma_J j_h(\mathcal{F}^{n, \text{ext}}_R, t; u,v) dt,
\end{align*}
defines the ghost penalty bilinear form which - roughly speaking - allows to carry over control from one element to its neighbour.
	Here, $\tilde \gamma_J = \left( 1 + \frac{\Delta t}{h} \right) \gamma_J$, where we refer to \cite[Section 4.4]{LO_ESAIM_2019} for an explanation of the scaling with $\Delta{t}$ and $h$ that is relevant in the anisotropic case $h\not\sim \Delta t$. Note that the analysis guarantees that the chosen scaling is sufficient for stability and bounded condition-numbers. The question if this scaling is also necessary is left open for future research. $\gamma_J$ is a stabilisation constant to be chosen sufficiently large, c.f. \cite{preuss18, heimann20} for an analysis and \Cref{sssect:stab_choice} for a computational study.
In total, the discrete time slab problem reads: Find $u \in W_h^n$ such that
\begin{equation}
 B^n(u,v) + B_{\text{upw}}^n(u,v) + J^n(u,v) = f^n(v) + f^n_{\text{upw}}(v) \quad \forall v \in W_h^n. \label{discreteproblem}
\end{equation}
The time slab solutions can be computed one after another yielding a global space-time solution $u_h$.
The same applies to all following discrete problem statements.
\subsection{Continuous Galerkin method}\label{ssection:CG}
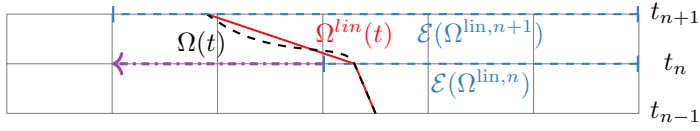
\begin{figure} \centering
  \begin{tikzpicture}[xscale=1.4, yscale=0.66]
  \draw[gray] (2,0) grid[xstep=1, ystep=1] (8,2);
  \begin{scope}[xscale=-1, xshift=-10cm]
  \draw[thick, Set1-A] (4.5,0) -- (4.7,1) -- (6.1,2);
  \draw[thick, dashed] (4.5,0)
                .. controls (4.6,0.5) .. (4.7,1)
                .. controls (4.8,1.5) and (5.5,1) .. (6.1,2);
   \end{scope}
   \draw[{Bar[]}-{Bar[]}, thick, dashed, Set1-B] (8,1) -- (5,1);
   \node[Set1-B] at (6.5,0.6) {$\mathcal{E}(\Omega^{\text{lin}, n})$};
   \draw[{Bar[]}-{Bar[]}, thick, dashed, Set1-B] (8,2) -- (3,2);
   \node[Set1-B] at (6.5,1.6) {$\mathcal{E}(\Omega^{\text{lin}, n+1})$};
   
   \draw[->, very thick, dashdotted, Set1-D] (5,1) -- (3,1);
   
   \node at (8.35,0) {$t_{n-1}$};
   \node at (8.35,1) {$t_{n}$};
   \node at (8.35,2) {$t_{n+1}$};
   \node at (3.85,1.4) {$\Omega(t)$};
   \node[Set1-A] at (5.3,1.6) {$\Omega^{lin}(t)$};
 \end{tikzpicture}\vspace*{-0.45cm}
\caption{Illustration of the domains of definition of the initial data at $t_n$ as needed and given from the previous time step if the DG-regions were used. The region where an extension of $u$ from $I_n$ is needed is depicted in purple.}
\label{EPG_intro}\vspace*{-0.5cm}
 \end{figure}

In general, the structure of the discrete regions for the CG method resembles that of the DG method. For the \dof~associated to the end point of a time slice, we  introduce $\EplOmn$, an extension of $\EOmn$, $\EOmn \subseteq \EplOmn$, which additionally should satisfy the constraint
\begin{equation}
 \mathcal{E}(\Omega^{\text{lin}, n+1}) \subseteq \EplOmn \quad n=1,\dots,N-1. \label{eq:CCEdomconstr}
\end{equation}
This is the discrete version of the domain $\Osnp{n}$ introduced in \Cref{CG_cont_intro}, with the only difference that the discrete extended domain aligns with the mesh, as illustrated in \cref{EPG_intro}. 
We start with a variant of $\Omlin{t}$,
 $\Omega^{\text{lin}}_\epsilon(t) := \{ x \in \tilde \Omega \, | \, \phi^{\text{lin}}(x,t) < \epsilon \}$,
which is used to define a strip around the discrete boundary $\mathcal{S}^n_\delta$, and $\EplOmn$:
\begin{subequations}
\begin{align}
 \mathcal{S}^n_\delta :&= \Omega^{\text{lin}}_{\delta}(t_n) \backslash \Omega^{\text{lin}}_{-\delta}(t_n),\\
 \mathcal{S}(\Omega^{\text{lin},n}) :&= \{ x \in T \textnormal{ for some } T \in \mathcal{T}_h \textnormal{ s.t. } T \cap \mathcal{S}^{n}_\delta \neq \varnothing \}, \\
 \hypertarget{def:EplOmn}{\mathcal{E}^{+}(\Omega^{\text{lin},n}) :&= \{ x \in T \textnormal{ for some } T \in \mathcal{T}_h \textnormal{ s.t. } T \cap \Omega^{\text{lin}}_{\delta}(t_{n}) \neq \varnothing \} }.
\end{align}
\end{subequations}
Note that for a time slice $I_n$, the regions $\EplOmn$ and $\mathcal{S}(\Omega^{\text{lin},n})$ depend only on the discrete function $\phi^{\text{lin}}$ at $t_n$, the end time point of the interval.
Moreover, $\delta$ is assumed to be given as follows:
\begin{assumption}
 We assume $\delta$ to be some small number such that $\EOmn \subseteq \EplOmn$ and the constraint \cref{eq:CCEdomconstr} is satisfied.
\end{assumption}
 If $\phi^{\text{lin}}$ would be a signed distance function, we could set
$
  \delta \geq \Delta t \| w\|_{L^\infty(0,T,L^\infty(\tOmega))}
$. 
For practical computations, we recommend to multiply this expression by a factor $\epsilon_f >1$ and to check \cref{eq:CCEdomconstr} in each time step to optimize $\epsilon_f$ when necessary. For a sufficiently large $\epsilon_f$, our framework reduces to a global extension on the whole mesh, such as suggested in \cite{anselmann2021cutfemNSE}. We comment on this option in \Cref{subsectnze}.

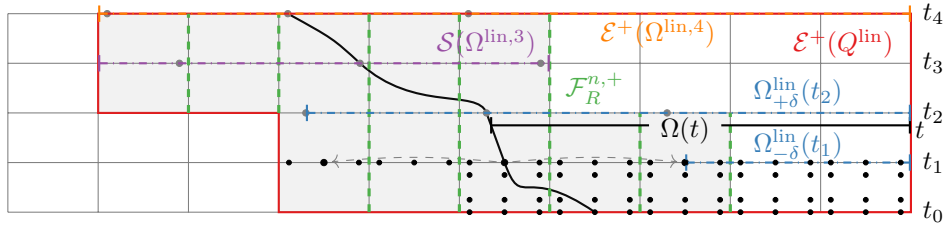
\begin{figure}[h!] 
 \centering
 \begin{tikzpicture}[xscale=1.2, yscale=0.66]
  \draw[gray] (0,0) grid[xstep=1, ystep=1] (10,4);
  \begin{scope}[xscale=-1, xshift=-10cm]
  \draw[thick] (3.5,0) .. controls (4.3, 1) and (4.3,0) .. (4.5,1)
                .. controls (4.6,1.5) .. (4.7,2)
                .. controls (4.8,2.5) and (5.5,2) .. (6.1,3)
                .. controls (6.4,3.5) .. (6.9,4);
   \end{scope}
   
   \filldraw (5.5,1) ellipse [y radius=1.8pt, x radius=1pt];
   \filldraw (3.5,1) ellipse [y radius=1.8pt, x radius=1pt];
   \filldraw (7.5,1) ellipse [y radius=1.8pt, x radius=1pt];
   \draw[->, gray, shorten >= 3pt, shorten <= 3pt, dashed] (5.5,1) to[bend left=10] (7.5,1);
   \draw[->, gray, shorten >= 3pt, shorten <= 3pt, dashed] (5.5,1) to[bend right=10] (3.5,1);
   
   \filldraw[gray] (5.3,2) ellipse [y radius=1.8pt, x radius=1pt];
   \filldraw[gray] (3.3,2) ellipse [y radius=1.8pt, x radius=1pt];
   \filldraw[gray] (7.3,2) ellipse [y radius=1.8pt, x radius=1pt];
   
   \filldraw[gray] (3.9,3) ellipse [y radius=1.8pt, x radius=1pt];
   \filldraw[gray] (1.9,3) ellipse [y radius=1.8pt, x radius=1pt];
   \filldraw[gray] (5.9,3) ellipse [y radius=1.8pt, x radius=1pt];
   
   \filldraw[gray] (3.1,4) ellipse [y radius=1.8pt, x radius=1pt];
   \filldraw[gray] (1.1,4) ellipse [y radius=1.8pt, x radius=1pt];
   \filldraw[gray] (5.1,4) ellipse [y radius=1.8pt, x radius=1pt];
   
   
   \draw[Set1-A, thick] (3,0) -- (3,2) -- (1,2) -- (1,4) -- (10,4) -- (10,0) -- cycle; 
   
   
   \fill[gray, fill opacity=0.1] (3,0) rectangle (8,2) -- (1,2) rectangle (6,4); 

   \draw[very thick, dashed, Set1-C] (2,2) -- (2,4) (3,2) -- (3,4) (4,0) -- (4,4) (5,0) -- (5,4) (6,0) -- (6,4) (7,0) -- (7,2) (8,0) -- (8,2);
   
   \node[Set1-A] at (9.25,3.4) {$\mathcal{E}^{+}(Q^{\text{lin}})$};
    \node[Set1-C] at (6.5,2.5) {$\mathcal{F}^{n,+}_R$};
    
       \node at (10.25, 0) {$t_0$};
   \node at (10.25, 1) {$t_1$};
   \node at (10.25, 2) {$t_2$};
   \node at (10.25, 3) {$t_3$};
   \node at (10.25, 4) {$t_4$};
   
   \node at (10.1, 1.65) {$t$};
   \draw[{Bar[]}-, thick] (5.34,1.75) --(7.1,1.75);
   \draw[-{Bar[]}, thick] (7.9,1.75) -- (10, 1.75);
   \node at (7.5,1.65) {$\Omega(t)$};
   
   \draw[Set1-B, dashdotted, thick, {Bar[]}-{Bar[]}] (7.5,1) -- (10,1);
   \node[Set1-B] at (8.75,1.35) {$\Omega^{\text{lin}}_{-\delta}(t_1)$};

   \draw[Set1-B, dashdotted, thick, {Bar[]}-{Bar[]}] (3.3,2) -- (10,2);
   \node[Set1-B] at (8.75,2.4) {$\Omega^{\text{lin}}_{+\delta}(t_2)$};
   
   \draw[Set1-D, dashdotted, thick, {Bar[]}-{Bar[]}] (1,3) -- (6,3);
   \node[Set1-D] at (5.3,3.4) {$\mathcal{S}(\Omega^{\text{lin,3}})$};
   
   \draw[Set1-E, dashdotted, thick, {Bar[]}-{Bar[]}] (1,4) -- (10,4);
   \node[Set1-E] at (7.2,3.6) {$\mathcal{E}^{+}(\Omega^{\text{lin,4}})$};
   
   \foreach \c in {5,6,7,8,9} {
   \foreach \y in {0,0.25,0.75,1} {
   \filldraw (\c + 0.112702, \y) ellipse [y radius=0.75*1.8pt, x radius=0.75pt];
   \filldraw (\c + 0.5,\y) ellipse [y radius=0.75*1.8pt, x radius=0.75pt];
   \filldraw (\c + 0.887298, \y) ellipse [y radius=0.75*1.8pt, x radius=0.75pt]; } }
   
   
   \filldraw (3.112702, 1) ellipse [y radius=0.75*1.8pt, x radius=0.75pt];
   \filldraw (3.5, 1) ellipse [y radius=0.75*1.8pt, x radius=0.75pt];
   \filldraw (3.887298, 1) ellipse [y radius=0.75*1.8pt, x radius=0.75pt];
   
   \filldraw (4.112702, 1) ellipse [y radius=0.75*1.8pt, x radius=0.75pt];
   \filldraw (4.5, 1) ellipse [y radius=0.75*1.8pt, x radius=0.75pt];
   \filldraw (4.887298, 1) ellipse [y radius=0.75*1.8pt, x radius=0.75pt];
   
 \end{tikzpicture}\vspace*{-0.4cm}
 \caption{Discrete regions in the definition of the ghost-penalty of the CG method. 
 The black circles represent the \pdofs~of the discrete function space for the first time step.}
 \label{newgpsketch}\vspace*{-0.2cm}
 \end{figure}

 The discrete extension should be realised by a merely spatial ghost penalty stabilisation. The according facet set is hence given as
 \begin{align*}
  \mathcal{F}^{n, +}_R \!=\! \{F \!\in\! \mathcal{F} \textnormal{ s.t. } \exists T_1 \neq T_2, T_1 \subseteq \EplOmn , T_2 \subseteq \mathcal{S}(\Omega^{\text{lin},n}) \textnormal{ with } F = T_1 \cap T_2 \}.
 \end{align*}
 Again, we give a sketch of this construction in \cref{newgpsketch}.
 
 In order to introduce the discrete function spaces next, we pose an assumption on the basis of the function space $\mathcal{P}^{\ktsol}([t_{n-1},t_n])$ contained in $\bgVhnks{n}{k}$:
 \begin{assumption} \label{assumption_on_Pkt}
  The basis $\{ p_{0},\dots,p_{\ktsol} \}$ of the space $\mathcal{P}^{\ktsol}([t_{n-1},t_n])$ contains two functions $p_{i_\text{lo}}(t)$ and $p_{i_\text{up}}(t)$ which can be associated to the lower and upper domain boundary $t_{n-1}$ and $t_n$, s.t. $p_{i}(t_{n-1}) = \delta_{i_\text{lo}i}$ and $p_{i}(t_{n}) = \delta_{ i_\text{up}i}$
  ~$\forall i=0, \dots, \ktsol$.
  \end{assumption}
 The discrete trial function space is then defined as follows in terms of the nomenclature known from the DG method: $U_h^n = U_{h,0}^n + u_{\text{init}}$ with     \vspace*{-0.1cm}   
 \begin{align*}
  U_{h,0}^{n} := \{ v \in &\Whcut{\ksol}{ \EQlinn \cup (\EplOmn \times \{ t_n \})} \circ (\Theta_h^{n, \text{st}})^{-1} \mid v^{n-1}_+ = 0 \} \vspace*{-0.1cm}
 \end{align*}
 and $u_{\text{init}}(x,t) = \Piu (u^{n-1}_-(x)) \cdot p_{i_{\text{lo}}}(t) $
 where $\Piu (u^{n-1}_-)$ serves as initial data for the solution as for the DG method, but with a strong imposition.
 The relevant discrete region for each time slice is given as $\EQlinn \cup (\EplOmn \times \{ t_n \})$, where the expression $\EplOmn \times \{ t_n \}$ will allow for the discrete extension motivated above at the end of each time slice. The \dofs~from $\bgVhks{\ksol}$ related to the region $\EplOmn \times \{ t_n \}$ will stem from $p_{i_{\text{up}}}$, so the construction is well-posed by \cref{assumption_on_Pkt}. The structure of the domain is illustrated in \cref{newgpsketch} by black circles.
 %
 For the space of test functions, continuity beyond time slice boundaries is not required, so we define:\vspace*{-0.1cm}   
 \begin{align*}
  V_h^{n} := \Whcut{(\khsol,\ktsol-1)}{  \EQlinn \cup (\EplOmn \times \{ t_n \})} \circ (\Theta_h^{n, \text{st}})^{-1}.\vspace*{-0.1cm}
 \end{align*} 
 Note that the polynomial order in time is one less as usual with CG methods.
 
 Now we have all tools in place to state the variational formulation of the CG method within one time slab. It reads: Find $u = u_0 + u_{\text{init}}$ with $u_0 \in U_{h,0}^{n}$ s.t. \vspace*{-0.1cm}   
\begin{equation*}
 B^{n}(u_0,v) + J^n(u_0,v) + j_h(\mathcal{F}_R^{n,+},t_n; (u_0)_-^n, v_-^n)  = f^n(v) - B^{n}(u_{\text{init}},v) \quad \forall v \in V_h^{n}. \vspace*{-0.1cm}
\end{equation*}
Here, both ghost penalty parts only act on the \dofs~ of $u_0$ that needs to be solved for, but not on $u_{\text{init}}$. We further note that the decomposition $u = u_0 + u_{\text{init}}$ and hence the solution $u$ depends on definition of the basis function $p_{i_{\text{lo}}}$.
 \begin{remark}\label{rem:CGtg}
  There exists a variant of the suggested CG method which we denote as CG$\square$, where only the ghost penalty stabilisation $J^{n\square}(u,v) = \int_{t_{n-1}}^{t_n} $ $j_h(\mathcal{F}_R^{n,+},t; u, v)$ $\mathrm{d}t$ is applied, which is defined on the smallest tensor-product domain which contains both domains used before.
  This can be regarded as an implementational simplification, although more \dofs~are needed as elaborated on in \Cref{subsectnze}. Further details about this variant are given in the supplementary material.
 \end{remark}
 \subsection{Galerkin-Collocation method}\label{ssection:GCC}
 The GCC method builds on the same idea as the CG method, but extends the strategy by one step: In addition to the continuity of the discrete solution function, we ask for continuity also in higher order time derivatives up to order $\ktreg$.
 Correspondingly, we strengthen the assumption of \cref{assumption_on_Pkt} w.r.t. to the basis of $\mathcal{P}^{\ktsol}(\In{n})$:
 \begin{assumption} \label{strong_assumption_on_Pkt}
The basis $\{ p_{1},\dots,p_{\ktsol+1} \}$ of the space $\mathcal{P}^{\ktsol}([t_{n-1},t_n])$ contains $\ktreg +1$ functions $p_{i_\ell}(t)$ for $\ell=0,\dots,\ktreg$ which can be associated to the $\ell$th derivative at the boundary $t_{n-1}$: $\forall \ell=0,\dots,\ktreg, i=1,\dots,\ktsol+1, ~ \partial^{\ell} p_i (t_{n-1}) = \delta_{i_\ell i}$.
 \end{assumption}
We start by introducing a corresponding trial function space: $U_h^n = U_{h,0}^n + u_{\text{init}}$ with
 \begin{align*}
  U_{h,0}^{n} \!\!:= \{ v \in &\Whcut{\ksol}{ \EplOmn \!\times\! I_n} \circ\! (\Theta_h^{n, \text{st}})^{-1} \!\mid
   ( \partial_t^{\ell}  v )^{n-1}_+ = 0, \ell = 0,\dots,\ktreg\},
 \end{align*}
and $u_{\text{init}}(x,t) =\sum_{\ell=0}^{\ktreg}     \Piu (\partial_t^{\ell} u_-^{n-1}(x)) \cdot p_{i_\ell}(t)$
where 
$\Piu (\partial_t^\ell u^{n-1}_-)$ serve as initial data for the solution as for the CG and the DG method. We note however that we have to prescribe $\ktreg+1$ functions as initial data here.
 For the space of test functions, we now pick only order $\kttest\leq k_t -\ktreg -1$:
 \begin{align*}
   V_h^{n} := \Whcut{(\khsol,\kttest)}{\EplOmn \times I_n} \circ (\Theta_h^{n, \text{st}})^{-1}.
 \end{align*}
 For the ghost penalty we use the tensor-product variant CG$\square$ from \cref{rem:CGtg}\footnote{One could also construct GCC methods which generalise the original CG method and not CG$\square$ regarding the structure of function space and stabilisation. However, this would necessitate further extensions/ stabilisations in time derivatives of $u$, which is why we opt for the presented variant.}.  The formulation of the GCC method then reads: Find $u = u_0 + u_{\text{init}}$ with $u_0 \in U_{h,0}^{n}$ s.t. \vspace{-0.1cm}
\begin{align*}
 B^{n}(u,v) + J^{n\square}(u_0,v) &= f^n(v) \quad \forall v \in V_h^{n} \text{  and  } \\
 B^{n}_{\text{col}}(t_l^n;u,v) + j_h(\mathcal{F}_R^{n,+},t_l^n;u_0,v) &= f^n_{\text{col}}(t_l^n;v) \quad \forall v \in \Vhks{\khsol}, l=1,..,\ktcol,\\
  \text{with } B^n_{\text{col}}(t;u,v) := (\partial_t u + \vect{w}\cdot \nabla u, &v)_{\Oht{t}} + (\nabla u, \nabla v)_{\Oht{t}},~f^n_{\text{col}}(t;v) := (f,v)_{\Oht{t}},  
\end{align*}
with collocation points $\{t^n_1,\dots,t^n_{\ktcol}\}$. Here, similar to the CG method, we made sure that the ghost penalty only acts on the \dofs~of $u_0$. Again, the solution $u$ will depend on the concrete choice of the basis functions $p_{i_\ell},~\ell=0,..,\ktreg$.  

A possible choice for the triple $(\ktsol,\ktreg,\ktcol$) is $(3,1,1)$ which we also consider in the numerical experiments below. The one collocation point in this setting is chosen as the end point of each time interval, $t^n_1 = t^n_{\ktcol} = t^n$.
%
\section{Implementation of numerical integration in time} \label{subsect_num_int}
The discrete (bi)-linear forms presented in the previous section involve summands defined on the space-time domains $\Qhn$, c.f. e.g. \cref{DGbilinearforms1}. Hence, a numerical integration procedure is needed for those domains. More specifically, as $\Qhn$ was defined as the image of $\Qlinn$, cf. \cref{eq:introQhn}, and finite element computations are done element-wise, we apply the following calculation:
\begin{align}
 (f, v)_{\Qhn} &= \int_{\Qhn} f v \, \mathrm{d}x = \int_{\Qlinn} |\det{D\Thehst}| (f \circ \Thehst) \hat v \, \mathrm{d}\hat x, \quad v = \hat v \circ \Thehst \\
 &= \sum_{T \in \mathcal{T}_h} \int_{\Qlinn \cap (T \times I_n) } |\det{D\Thehst}| (f \circ \Thehst) \hat v \, \mathrm{d} \hat x
\end{align}
For the other summands, similar equations hold. Hence, we only need to implement an algorithm to approximate
\begin{equation}
 \int_{\Qlinn \cap (T \times I_n) }  \tilde f(x,t) \mathrm{d}(\hat x,t) = \int_{t_{n-1}}^{t_n}  \int_{\Omlin{t} \cap T} \, \tilde f(\hat x,t) \, \mathrm{d}\hat x \,\mathrm{d}t \,\label{taskoftint}
\end{equation}
for each $T \in \Th$ for some integrand $\tilde f$. For a given time $\hat t$, an integration rule on $\Omlin{\,\hat t\,} \cap T$ is assumed to be available, e.g. on simplices by tessellation. Hence, the inner integral in \cref{taskoftint} is well-suited for numerical approximation, and it remains to replace the integration $\int_{t_{n-1}}^{t_n}$ by a numerical quadrature. In the following, we will define two variants of doing so, where the first approach is mentioned for comparison with existing literature \cite{HLZ16,Z18} and the second is our suggested approach.

\begin{definition}[Topology-insensitive time integration]
 In the \emph{topology-insen\-sitive time integration}, $\int_{t_{n-1}}^{t_n}$ is approximated by Gaussian quadrature. Let $(t_i, \omega_i)_i$ be the Gaussian quadrature rule on $[t_{n-1},t_n]$ of order $k_{t,int}$. Then, \cref{taskoftint} will be approximated by
 \begin{equation}
  \sum_i \omega_i \int_{\Omlin{t_i} \cap T} \, f(\hat x,t_i) \, \mathrm{d}\hat x.
 \end{equation}
\end{definition}
Note however, that the topology of $\Omlin{t} \cap T$ might change over $t \in I_n$; for instance the domain boundary can leave or enter an element. This will cause the integrand $\int_{\Omlin{t} \cap T} \, f(\hat x,t)\, \mathrm{d}\hat x$ to show a discontinuity in the first derivative even if  $f$ is smooth, which impedes the accuracy of the topology-insensitive time integration. To circumvent this, we define a time integration that takes the cut topology changes of $\Qlinn$ into account.
\begin{definition}[Topology-preserving time integration]
  To approximate \cref{taskoftint}, we start by calculating the set $\mathcal{R}^\ast$ of subintervals of $I_n$ represented by a tuple $(t_a,t_b)$ in which no topology changes take place with \cref{alg1}.
\begin{algorithm}
\caption{Subdivision of the time interval}
\label{alg:time_I_subdiv}
\begin{algorithmic}[1]
\STATE{Set $\mathcal{R} = \{ t_{n-1}, t_n\} $ (Set of time points with topology changes).}
\STATE{Let $V$ be the set of vertices of $T$.}
\FOR{ $v \in V$}
  \STATE{Define the scalar function $\phi_v: I_n \to \mathbb{R}, t \mapsto \phi^{\text{lin}} (v, t)$. }
  \STATE{Search for the roots $\mathcal{R}_v$ of $\phi_v$.}
  \STATE{Set $\mathcal{R} \leftarrow \mathcal{R} \cup \mathcal{R}_v$.}
\ENDFOR
\STATE{Define $\mathcal{R}^\ast$ as the set of intervals with endpoints according to $\mathcal{R}$.} 
\RETURN{$\mathcal{R}^\ast$.}
\end{algorithmic}
\label{alg1}
\end{algorithm}\\[-0.3cm]
Denote by $I_n^s = (t_l, t_u) \in \mathcal{R}^\ast$ a subinterval in $\mathcal{R}^\ast$. 
Then, \cref{taskoftint} is approximated by the \emph{topology-preserving time integration} as
\begin{equation}
 \sum_{I_n^s \in \mathcal{R}^\ast} \sum_i \omega_i^s \int_{\Omlin{t_i} \cap T} \, f(\hat x,t_l^s) \, \mathrm{d}\hat x,
\end{equation}
where $(t_i^s,\omega_i^s)$ is the Gaussian integration rule transformed to the interval $I_n^s$.
\end{definition}
This time integration procedure exploits the fact that a (multi-)linear-in-space level set function is fully described by its function values on vertices of the finite element $T$. The functions $\phi_v$ within \cref{alg:time_I_subdiv} will in general be of higher polynomial order.

\begin{figure} \centering
\begin{tikzpicture}[scale=0.5]
 \begin{axis}[xmin=0,xmax=1,ymin=0,ymax=1, ylabel={$t$}, xlabel={$x$}, legend entries ={$t_i$}, legend style={at={(0.5,-0.2)},anchor=north,legend columns=1, draw=none}]
  \addplot[dashed] {0.330009}; \addplot[dashed] {0.0694318};
  \addplot[dashed] {1 - 0.330009}; \addplot[dashed] {1 - 0.0694318};
  \addplot[Set1-A,ultra thick] table {num_exp/visu/interf_out_172.dat};
  \addplot[Set1-B, only marks, mark options={scale=1.5}] table {num_exp/visu/rule_out_172_naive.dat};
 \end{axis}

\end{tikzpicture}
\begin{tikzpicture}[scale=0.5]
 \begin{axis}[xmin=0,xmax=1,ymin=0,ymax=1, yticklabels={,,}, xlabel={$x$}, legend entries ={$\mathcal{R}^\ast$}, legend style={at={(0.5,-0.2)},anchor=north,legend columns=1, draw=none}]
 \addplot[dashdotted] {0.27}; \addplot[dashdotted] {0.52};
  \addplot[Set1-A,ultra thick] table {num_exp/visu/interf_out_172.dat};
  \addplot[Set1-B, only marks, mark options={scale=1.5}] table {num_exp/visu/rule_out_172_tsint.dat};
 \end{axis}
\end{tikzpicture}\hspace*{0.2cm}
\begin{tikzpicture}[scale=2.5]

  \coordinate (A1) at (-0.50, 0.50, 0.00);
  \coordinate (B1) at ( 0.50, 0.50, 0.00);
  \coordinate (C1) at ( 0.00, 0.50,-0.80);

  \coordinate (A2) at (-0.50, 1.50, 0.00);
  \coordinate (B2) at ( 0.50, 1.50, 0.00);
  \coordinate (C2) at ( 0.00, 1.50,-0.80);

  \coordinate (A3) at (-0.50, 1.05, 0.00);
  \coordinate (B3) at ( 0.50, 1.05, 0.00);
  \coordinate (C3) at ( 0.00, 1.05,-0.80);


  \coordinate (R1) at ($(A2)!0.5!(C2)$);
  \coordinate (R2) at ($(B2)!0.25!(C2)$);
  \coordinate (R3) at ($(A1)!0.4!(A2)$);
  \coordinate (R4) at ($(B1)!1.0!(B3)$);

  \coordinate (S3) at ($(A3)!0.31!(C3)$);

  \node[left,scale=0.7] at (A3) {\small $t^2_\ast$};
  \node[left,scale=0.7] at (R3) {\small $t^1_\ast$};
  \node[left,scale=0.7] at (A1) {\small $t^{n\!-\!1}$};
  \node[left,scale=0.7] at (A2) {\small $t^{n}$};

  \draw[thin,densely dotted, opacity=0.5] (B1) -- (C1);
  \draw[thin,densely dotted, opacity=0.5] (C1) -- (A1);

  \draw[thin,densely dotted, opacity=0.5] (C1) -- (C2);

  \draw[thin,densely dotted] (B3) -- (C3);
  \draw[thin,densely dotted] (C3) -- (A3);


  \draw[thin,fill=Set1-C,opacity=0.25] (S3) -- (B3) -- (C3) -- cycle;
  \draw[thin,fill=Set1-D,opacity=0.25] (S3) -- (A3) -- (B3) -- cycle;

  \draw[thin,fill=Set1-C,opacity=0.25] (R1) -- (R2) -- (C2) -- cycle;
  \draw[thin,fill=Set1-D,opacity=0.25] (R1) -- (R2) -- (B2) -- (A2) -- cycle;

  \draw[thin,fill=Set1-C,opacity=0.25] (A1) -- (B1) -- (C1) -- cycle;

  \draw[fill=Set1-A,opacity=0.1] (R1) -- (R2) to[in=105,out=-85] (R4) to[in=5,out=-165] (R3) to[in=-100,out=33] (R1);
  \draw[draw=Set1-A!80!black,fill=none] (R1) -- (R2) to[in=105,out=-85] (R4) to[in=5,out=-165] (R3) to[in=-100,out=33] (R1);

  \draw[thin] (A3) -- (B3);

  \draw[Set1-A!60!Set1-C,thick] (B3) -- (S3);

  \draw[thick] (A1) -- (B1);

  \draw[thick] (A2) -- (B2);
  \draw[thick] (B2) -- (C2);
  \draw[thick] (C2) -- (A2);

  \draw[thick] (A1) -- (A2);
  \draw[thick] (B1) -- (B2);

  \coordinate (A3B3) at ($(A3)!0.5!(B3)$);
  \coordinate (A2B2) at ($(A2)!0.5!(B2)$);
  \coordinate (A1B1) at ($(A1)!0.5!(B1)$);

  \coordinate (q) at ($(R1)!0.5!(R2)!0.5!(C2)$);
  \node[right=0.5cm,scale=0.7] (ql) at (q) {$\Omlin{t_i} \cap T$};
  \draw[->] (ql) to[in=50,out=150] (q);

  \coordinate (p) at ($(R1)!0.5!(R2)!0.1!(B2)!0.5!(A2)$);

  \node[below=0.15cm,scale=0.7] at (A1B1) {$T \times I_n$};

\end{tikzpicture} \vspace*{-0.35cm}
\caption{Comparison of the topology-insensitive time integration in 1+1D (left) with the topology-preserving time integration in 1+1D (center) and 2+1D (right). In the left and middle picture, blue dots represent integration points, and the red line represents the space-time interface. In the left picture, vertical lines indicate the Gaussian quadrature points (in time), which are independent of the cut configuration. In the middle, vertical lines represent the subdivision of the time interval into subintervals without topology change, on which the Gaussian rules (in time) are then applied. On each quadrature node in time a quadrature rule in space on the cut domain is applied resulting in the space-time integration points indicated by blue dots.
}\vspace*{-0.5cm}
\label{num_int_comp}
\end{figure}
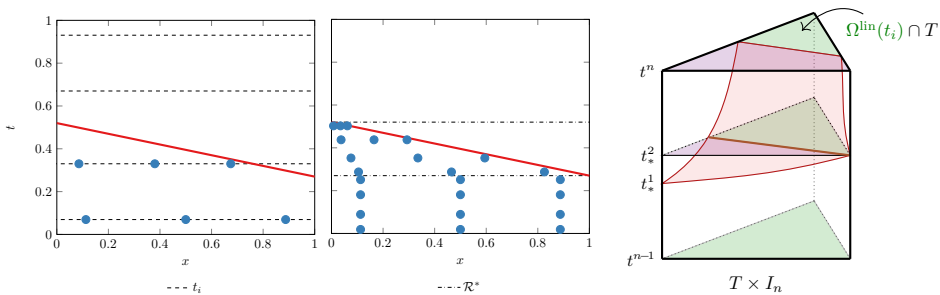

The topology-preserving time integration is our method of choice and for the numerical studies in this work we choose integration rules in time of exactness degree $2 (k_t + 1)$. In \cref{num_int_comp}, we show the quadrature points generated by both approaches for simple examples in one and two space dimensions. Later in this paper, we investigate numerically differences to the topology-insensitive time integration in \Cref{tint_numexp}.


\section{Numerical investigations} \label{sect:numexp}
In this section, we want to demonstrate the numerical properties of the proposed methods. We consider three problem geometries to investigate specific aspects, each in one of the following subsections.
Recall that the notation $\ksol=r$ means $\khsol\!=\!\ktsol\!=\!r$ and analogously $\klset=r$ means $\khlset\!=\!\ktlset\!=\!r$.

All experiments are performed with \texttt{ngsxfem} \cite{xfem_joss}, an unfitted finite element extension of \texttt{ngsolve}. Reproduction data are available at \url{https://gitlab.gwdg.de/fabian.heimann/repro-ho-unf-space-time-fem}.

The linear systems are solved with the direct solvers \texttt{umfpack} and \texttt{pardiso} of IntelMKL. As a basis for $\mathcal{P}^{k_t}([t_{n-1},t_n])$, we implement a Lagrange basis with the Gauß-Lobatto points between $t_{n-1}$ and $t_n$ of order $k_t$ as Lagrange points for DG and CG. These satisfy \cref{assumption_on_Pkt}. For GCC of order 3, we implement cubic Hermite spline polynomials, which satisfy \cref{strong_assumption_on_Pkt}.
\subsection{Kite geometry} \label{sect:kite}
The first geometry is a circle, which deforms with time into a kite-shaped object in two spatial dimensions. It is described in terms of the following level set function $\phi$
\begin{equation*}
 \rho(t, y) = (1 - y^2) \cdot t, \quad r = \sqrt{(x- \rho)^2 + y^2}, \quad \phi = r - r_0, \quad r_0 = 1,
\end{equation*}
where $\tilde \Omega = [-3.5, 3.5] \times [-1.5,1.5]$, $\vect{w} = (\frac{\partial \rho}{\partial t}, 0)$, $T=0.5$. The geometry as well as some discrete regions for the DG method are shown in \cref{kiteillu}. The right-hand-side is calculated such that $u = \cos (\pi \frac{r}{r_0}) \cdot \sin (\pi t)$, and the stabilisation constant is chosen as $\gamma_J = 0.05$. We use unstructured simplicial meshes with mesh size $h=0.9 \cdot 0.5^{i_s}$, $\Delta t = 0.5 \cdot 2^{-i_t-1}$. In the CG and GCC methods, the extension factor is $\epsilon_f = 1.1$.
\subsubsection{Geometry approximation quality}
In a first experiment, we investigate the quality of the discrete domain approximation by measuring the distance to the exact interface, $\max_{t \in [0,T]} \mathrm{dist}(\partial \Omega^h(t), \partial \Omega(t))$, for the kite geometry with $\ksol=\klset$ and refinements in space and time, $i=i_t=i_s$. As can be seen in \cref{plot_kite_geom_approx}, the discrete regions satisfy the asymptotic error bound given in \cref{discretestregionapproxprop}.

\subsubsection{Convergence of methods of matching space and time order}
First, we evaluate convergence of the numerical error for $i_s=i_t=i$, $\ksol=\klset$. The range of parameters $k \in \{1,2,\dots,6\}$ is applied to the DG and CG methods, whilst the GCC method is investigated for $k= 3$. In \cref{plot_dg_kite_conv} (left) and \cref{plot_cg_kite_conv}, we display the numerical errors for simultaneous space-time refinements $i=1,2,\dots$.
    \begin{figure}
     \begin{tikzpicture}[scale=0.75]
    \begin{semilogyaxis}[ xlabel=$i$, ylabel=$\| u_h(T) - u (T)\|_{L^2(\Omega(T))}$, legend entries ={ $k=1$, $k=2$, $k=3$, $k=4$, $k=5$,$k=6$
    }, legend style={anchor=north,legend columns=1, draw=none, legend cell align=left, fill=none}, legend pos =south west,
    x label style={at={(axis description cs:0.5,0.15)},anchor=center},
    ymax=3, ymin=3e-11
     ]
     \addplot+[thick,mark options={scale=1.5}] table {num_exp/out/conv_kite_DG_ks1_kt1_both_nref8_gamma0.05_sm2.dat};
     \addplot+[thick,mark options={scale=1.5}] table {num_exp/out/conv_kite_DG_ks2_kt2_both_nref8_gamma0.05_sm2.dat};
     \addplot+[thick,mark options={scale=1.5}] table {num_exp/out/conv_kite_DG_ks3_kt3_both_nref8_gamma0.05_sm2.dat};
     \addplot+[thick,mark options={scale=1.5}] table[skip coords between index={6}{7}] {num_exp/out/conv_kite_DG_ks4_kt4_both_nref7_gamma0.05_sm2.dat};
     \addplot+[thick,mark options={scale=1.5}] table[skip coords between index={5}{6}] {num_exp/out/conv_kite_DG_ks5_kt5_both_nref6_gamma0.05_sm2.dat};
     \addplot+[thick,mark options={scale=1.5}] table[skip coords between index={4}{5}] {num_exp/out/conv_kite_DG_ks6_kt6_both_nref5_gamma0.05_sm2.dat};
     \addplot[gray, dashed, domain=0:7] {(1/2^(x-0.2)))^2};
     \addplot[gray, dashed, domain=0:7] {(1/2^(x+0.5)))^3};
     \addplot[gray, dashed, domain=0:7] {(1/2^(x+1)))^4};
     \addplot[gray, dashed, domain=0:5] {(1/2^(x+1)))^5};
     \addplot[gray, dashed, domain=0:4] {(1/2^(x+1)))^6};
     \addplot[gray, dashed, domain=0:3] {(1/2^(x+1)))^7};
    \end{semilogyaxis}
     \node[scale=0.75] at (4.5,5.25) {$O(h^{k+1})= O(\Delta t^{k+1})$};
     \draw[scale=0.75, gray, dash=on 2.25pt off 2.25pt phase 0pt, line width=0.4*0.75pt] (6.25/0.75,5.25/0.75) -- (6.8/0.75,5.25/0.75);
   \end{tikzpicture} \hspace{0.2cm} \tikzexternalexportnextfalse
     \begin{tikzpicture}[scale=0.75]
    \begin{semilogyaxis}[ xlabel=$i$, ylabel=$\| u_h(T) - u (T)\|_{L^2(\Omega(T))}$, legend entries ={$\gamma_J = 5 \cdot 10^4$,$\gamma_J = 5$, $\gamma_J = 5 \cdot 10^{-2}$,$\gamma_J = 5 \cdot 10^{-4}$}, legend style={anchor=north,legend columns=1, draw=none, fill=none},legend pos =south west, x label style={at={(axis description cs:0.6,0.15)},anchor=center}]
     \addplot+[thick,mark options={scale=1.5}] table {num_exp/out/conv_kite_DG_ks4_kt4_both_nref5_gamma50000.0_sm2.dat};
     \addplot+[thick,mark options={scale=1.5}] table {num_exp/out/conv_kite_DG_ks4_kt4_both_nref5_gamma5.0_sm2.dat};
     \addplot+[thick,mark options={scale=1.5}] table {num_exp/out/conv_kite_DG_ks4_kt4_both_nref5_gamma0.05_sm2.dat};
     \addplot+[thick,mark options={scale=1.5}] table {num_exp/out/conv_kite_DG_ks4_kt4_both_nref5_gamma0.0005_sm2.dat};
     \addplot+[thick,mark options={scale=1.5}] table {num_exp/out/conv_kite_DG_ks4_kt4_both_nref5_gamma5e-10_sm2.dat}; \label{gamma5e-10}
     \addplot+[thick,mark options={scale=1.5}] table {num_exp/out/conv_kite_DG_ks4_kt4_both_nref3_gamma0.0_sm2.dat}; \label{gamma0}
     \addplot[gray, dashed, domain=0:4] {(1/2^(x+1.3)))^5};

          \node at (rel axis cs: 0.75,0.8) {\shortstack[r]{
$\gamma_J = 5 \cdot 10^{-10}$ \ref{gamma5e-10} \\ $\gamma_J = 0$ \ref{gamma0}}};
    \end{semilogyaxis}
              \node[scale=0.75] at (4.5,5.25) {$O(h^5)=O(\Delta t^5)$};
     \draw[scale=0.75, gray, dash=on 2.25pt off 2.25pt phase 0pt, line width=0.4*0.75pt] (6./0.75,5.25/0.75) -- (6.55/0.75,5.25/0.75);
   \end{tikzpicture}
    \vspace*{-0.8cm}
   \caption{Left: Convergence of the DG method for the kite case. Right: Influence of the stabilisation parameter $\gamma_J$ on the convergence behaviour of the DG method for $\ksol=\klset = 4$.}
   \label{plot_dg_kite_conv}
   \vspace*{-0.25cm}
  \end{figure}
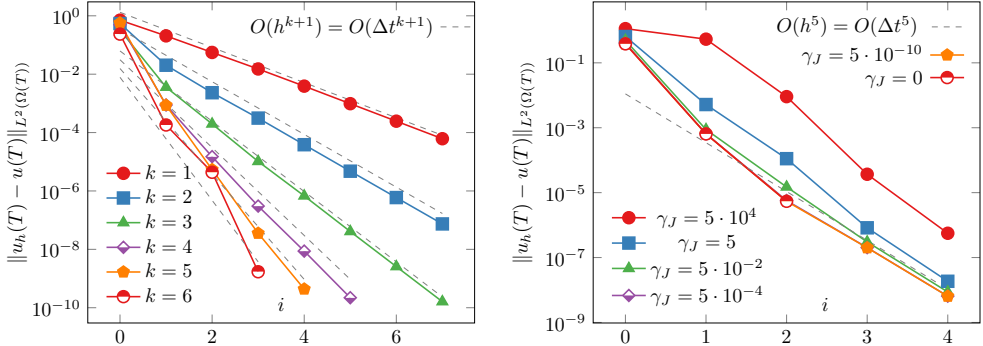
        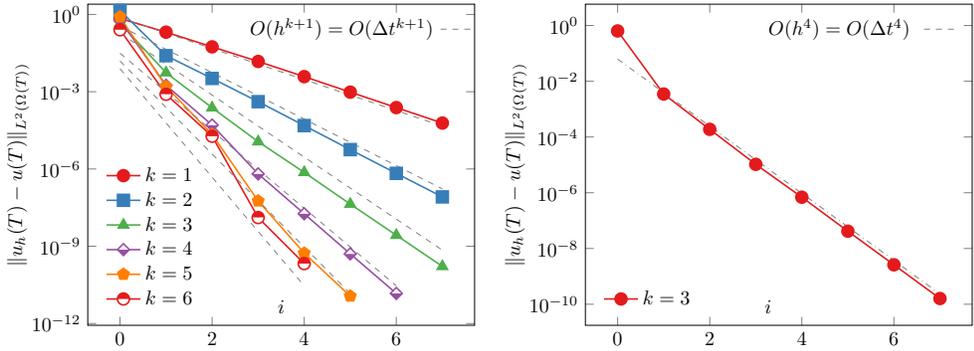
\begin{figure}
     \begin{tikzpicture}[scale=0.75]
    \begin{semilogyaxis}[ xlabel=$i$, ylabel=$\| u_h(T) - u (T)\|_{L^2(\Omega(T))}$, legend entries ={ $k=1$, $k=2$, $k=3$, $k=4$, $k=5$, $k=6$}, legend style={anchor=north,legend columns=1, draw=none, legend cell align=left, fill=none},legend pos =south west,
    x label style={at={(axis description cs:0.5,0.15)},anchor=center},
     ymax =2.5
    ]
     \addplot+[thick,mark options={scale=1.5}] table {num_exp/out/conv_kite_CG_ks1_kt1_both_nref8_gamma0.05_sm2.dat};
     \addplot+[thick,mark options={scale=1.5}] table {num_exp/out/conv_kite_CG_ks2_kt2_both_nref8_gamma0.05_sm2.dat};
     \addplot+[thick,mark options={scale=1.5}] table {num_exp/out/conv_kite_CG_ks3_kt3_both_nref8_gamma0.05_sm2.dat};
     \addplot+[thick,mark options={scale=1.5}] table {num_exp/out/conv_kite_CG_ks4_kt4_both_nref7_gamma0.05_sm2.dat};
     \addplot+[thick,mark options={scale=1.5}] table {num_exp/out/conv_kite_CG_ks5_kt5_both_nref6_gamma0.05_sm2.dat};
     \addplot+[thick,mark options={scale=1.5}] table {num_exp/out/conv_kite_CG_ks6_kt6_both_nref5_gamma0.05_sm2.dat};
     \addplot[gray, dashed, domain=0:7] {(1/2^(x+0.2)))^2};
     \addplot[gray, dashed, domain=0:7] {(1/2^(x+0.5)))^3};
     \addplot[gray, dashed, domain=0:7] {(1/2^(x+0.6)))^4};
     \addplot[gray, dashed, domain=0:6] {(1/2^(x+1)))^5};
     \addplot[gray, dashed, domain=0:5] {(1/2^(x+1)))^6};
     \addplot[gray, dashed, domain=0:4] {(1/2^(x+1)))^7};
    \end{semilogyaxis}
              \node[scale=0.75] at (4.5,5.25) {$O(h^{k+1})= O(\Delta t^{k+1})$};
     \draw[scale=0.75, gray, dash=on 2.25pt off 2.25pt phase 0pt, line width=0.4*0.75pt] (6.25/0.75,5.25/0.75) -- (6.8/0.75,5.25/0.75);
   \end{tikzpicture} \hspace{0.2cm}
 \hspace*{-0.2cm}  
     \begin{tikzpicture}[scale=0.75]
    \begin{semilogyaxis}[ xlabel=$i$, ylabel=$\| u_h(T) - u (T)\|_{L^2(\Omega(T))}$, legend entries ={ $k=3$}, legend style={anchor=north,legend columns=1, draw=none}, legend pos =south west, x label style={at={(axis description cs:0.5,0.15)},anchor=east}]
     \addplot+[thick,mark options={scale=1.5}] table {num_exp/out/conv_kite_GCC_ks3_kt3_both_nref8_gamma0.05_sm2.dat};
     \addplot[gray, dashdotted, domain=0:7] {(1/2^(x+1.)))^4};
    \end{semilogyaxis}
         \node[scale=0.75] at (4.5,5.25) {$O(h^4)=O(\Delta t^4)$};
     \draw[scale=0.75, gray, dash=on 2.25pt off 2.25pt phase 0pt, line width=0.4*0.75pt] (6./0.75,5.25/0.75) -- (6.55/0.75,5.25/0.75);
   \end{tikzpicture} 
     \vspace*{-0.4cm}
   \caption{Convergence of the higher order CG (left) and GCC (right) methods for the kite case.}
      \label{plot_cg_kite_conv}
    \end{figure}
    We measure the numerical error in terms of the following norms:
\begin{align*}
 \| u \|_{L^2(\Omega(T))}^2 := \int_{\Omega(T)} u(\mathbf{x}, T)^2 \ \mathrm{d}\mathbf{x}, \quad
 \| u \|_{L^2( L^2(\Omega(t)), 0, T)}^2 := \int_0^T \int_{\Omega(t)} u(\mathbf{x}, t)^2 \ \mathrm{d}\mathbf{x} \ \mathrm{d}t.
\end{align*}
We find both errors to decrease as follows, where $k = \khsol = \ktsol = \khlset = \ktlset$, $h \sim \Delta t$,
\begin{equation}
 \| u - u_h \|_{L^2(\Omega(T))} + \| u - u_h \|_{L^2( L^2(\Omega(t)), 0, T)} = \mathcal{O}(h^{k+1}) = \mathcal{O}(\Delta t^{k+1}). \label{eq:generalerrorbound} 
\end{equation}
This confirms the expected higher order convergence property of the suggested methods.
\subsubsection{Choice of the stabilisation parameter} \label{sssect:stab_choice}
The stabilisation parameter of the ghost penalty has to be chosen large enough to fulfill sufficient conditions for stability of the method and to obtain reasonable bounds on the condition numbers of arising linear systems \cite{burman15, Schott17, preuss18}.
In \cref{plot_dg_kite_conv} (right) we illustrate the dependence of the numerical error on $\gamma_J$ for the problem considered in this subsection using $i=i_s=i_t$ and $\ksol = \klset = 4$ for the DG method.
Overall, $\gamma_J$ influences the absolute values of the errors, whilst the asymptotic behaviour remains intact even if the parameter is chosen relatively large. A particular limit is that of stabilisation constant $\gamma_J \to 0$. There, we observe that the method is stable in general. However, condition numbers may blow up, which impedes the practical solution of the linear system at some refinement level/ polynomial degree. (In this case, $i=2$ is the last level solved properly.) We explain this behavior as follows: For mere stability, i.e. assuming exact arithmetics, ghost penalty stabilization is not even necessary. Note that in \cite{LR_SINUM_2013} even robust, i.e. cut-position independent, a priori error bounds have been derived without ghost penalty stabilization in a comparable setting (involving the assumption of exact integration). However, taking computer arithmetics, i.e. round-off errors and their impact, into account the possible ill-conditioning of linear systems is practically hardly acceptable. A more detailed investigation of these observations (potentially also in combination with different time integration strategies) in terms of a full numerical analysis would be an interesting task for future research.

In conclusion, we may pick $\gamma_J$ small as long as a sufficiently accurate (and efficient) solution of linear systems is possible. We choose $\gamma_J=0.05$ for the remaining studies.
\subsubsection{Superconvergence}
\begin{figure}\vspace*{-0.25cm}
   \tikzexternalexportnextfalse
     \begin{tikzpicture}[scale=0.75]
    \begin{semilogyaxis}[ xlabel=$i_t$, ylabel=$\| u_h(T) - u (T)\|_{L^2(\Omega(T))}$, legend entries ={DG $\ktsol=1$, DG $\ktsol=2$, CG $\ktsol=1$, CG $\ktsol=2$}, legend style={,anchor=north,legend columns=1, draw=none},legend pos =south west, x label style={at={(axis description cs:0.5,0.15)},anchor=center}]
     \addplot+[thick,mark options={scale=1.5}] table {num_exp/out/conv_kite_DG_ks3_kt1_time_nref6_gamma0.05_sm2_ktls3_ro-1.dat};
     \addplot+[thick,mark options={scale=1.5}] table {num_exp/out/conv_kite_DG_ks5_kt2_time_nref5_gamma0.05_sm2_ktls5_ro-1.dat};
     \addplot+[thick,mark options={scale=1.5}] table {num_exp/out/conv_kite_CG_ks3_kt1_time_nref6_gamma0.05_sm2_ktls3_ro-1.dat};
     \addplot+[thick,mark options={scale=1.5}] table {num_exp/out/conv_kite_CG_ks5_kt2_time_nref5_gamma0.05_sm2_ktls5_ro-1.dat};
     \addplot[very thick, dashed, domain=2:4] {(1/2^(x+3.2)))^3.8}; \label{order38}
     \addplot[very thick, dashdotted, domain=3:5] {(1/2^(x+2.5)))^2.8}; \label{order28}
     \addplot[gray, dashed, domain=0:4.25] {(1/2^(x+3.7)))^3};
     \addplot[gray, dashed, domain=0:4] {(1/2^(x+3.3)))^4};

     \addplot[gray, dashed, domain=0:5] {(1/2^(x+2.5)))^2};

      \node[gray] at (rel axis cs: 0.85,0.67) {$O(\Delta t^2)$};
      \node[gray] at (rel axis cs: 0.9,0.2) {$O(\Delta t^3)$};
      \node[gray] at (rel axis cs: 0.65,0.07) {$O(\Delta t^4)$};

      \node[black] at (rel axis cs: 0.75,0.85) {\shortstack[r]{ \hypersetup{linkcolor=black} \ref{order28} $O(\Delta t^{2.8})$ \\ \hypersetup{linkcolor=black} \ref{order38} $O(\Delta t^{3.8})$}};

    \end{semilogyaxis}
   \end{tikzpicture} \hspace{0.2cm}
     \begin{tikzpicture}[scale=0.75]
     \begin{semilogyaxis}[ xlabel=$i_t$, ylabel=$\| u_h - u\|_{L^2( L^2(\Omega(t)), 0, T) }$, legend entries ={DG $\ktsol=1$, DG $\ktsol=2$, CG $\ktsol=1$, CG $\ktsol=2$}, legend style={anchor=north, legend columns=1, draw=none},legend pos =south west, x label style={at={(axis description cs:0.5,0.15)},anchor=center}]
     \addplot+[thick,mark options={scale=1.5}] table [x index = 0, y index =2] {num_exp/out/conv_kite_DG_ks3_kt1_time_nref6_gamma0.05_sm2_ktls3_ro-1.dat};
     \addplot+[thick,mark options={scale=1.5}] table [x index = 0, y index =2] {num_exp/out/conv_kite_DG_ks5_kt2_time_nref5_gamma0.05_sm2_ktls5_ro-1.dat};
     \addplot+[thick,mark options={scale=1.5}] table [x index = 0, y index =2] {num_exp/out/conv_kite_CG_ks3_kt1_time_nref6_gamma0.05_sm2_ktls3_ro-1.dat};
     \addplot+[thick,mark options={scale=1.5}] table [x index = 0, y index =2] {num_exp/out/conv_kite_CG_ks5_kt2_time_nref5_gamma0.05_sm2_ktls5_ro-1.dat};
     \addplot[gray, dashed, domain=0:5] {(1/2^(x+3.0)))^2};
     \addplot[gray, dashed, domain=0:4] {(1/2^(x+3.2)))^3};

      \node[gray] at (rel axis cs: 0.85,0.55) {$O(\Delta t^2)$};
      \node[gray] at (rel axis cs: 0.88,0.22) {$O(\Delta t^3)$};

     \end{semilogyaxis}
    \end{tikzpicture} \vspace*{-0.8cm}
   \caption{Superconvergence study of the DG and CG methods for $\ktsol=1,2$ where $\khsol = 2 \ktsol +1$. Refinement in time for a fine mesh. On the left-hand side, the error is measured in the $L^2(\Omega(T))$-norm, whilst on the right-hand side the same results are investigated in the $L^2( L^2(\Omega(t)), 0, T)$-norm.}
      \label{plot_ks3kt1_conv_time}
    \end{figure}
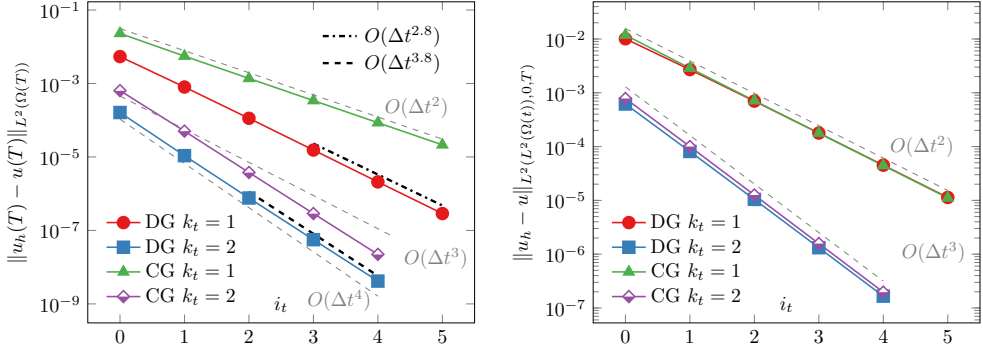

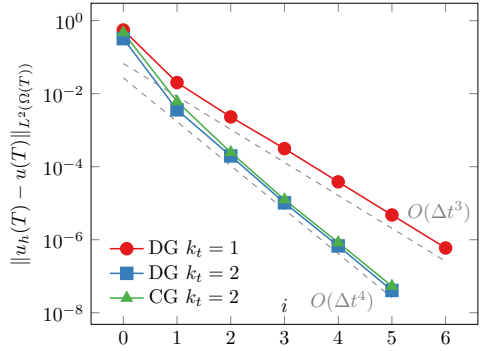
\begin{figure}\vspace*{-0.15cm} \centering
   \tikzexternalexportnextfalse
     \begin{tikzpicture}[scale=0.75]
    \begin{semilogyaxis}[ xlabel=$i$, ylabel=$\| u_h(T) - u (T)\|_{L^2(\Omega(T))}$, legend entries ={DG $\ktsol=1$, DG $\ktsol=2$, CG $\ktsol=2$}, legend style={,anchor=north,legend columns=1, draw=none},legend pos =south west, x label style={at={(axis description cs:0.5,0.15)},anchor=center}]
     \addplot+[thick,mark options={scale=1.5}] table {num_exp/out/conv_kite_DG_ks2_kt1_both_nref7_gamma0.05_sm2_ktls2.dat};
     \addplot+[thick,mark options={scale=1.5}] table {num_exp/out/conv_kite_DG_ks3_kt2_both_nref6_gamma0.05_sm2_ktls3.dat};
     \addplot+[thick,mark options={scale=1.5}] table {num_exp/out/conv_kite_CG_ks3_kt2_both_nref6_gamma0.05_sm2_ktls3.dat};
     \addplot[gray, dashed, domain=0:6] {(1/2^(x+1.3)))^3};
     \addplot[gray, dashed, domain=0:5] {(1/2^(x+1.3)))^4};


      \node[gray] at (rel axis cs: 0.9,0.35) {$O(\Delta t^3)$};
      \node[gray] at (rel axis cs: 0.65,0.07) {$O(\Delta t^4)$};


    \end{semilogyaxis}
   \end{tikzpicture} \vspace*{-0.2cm}
   \caption{Superconvergence study of the DG and CG methods for $\ktsol=1,2$ where $\khsol = \ktsol +1$. Refinement in space and time. The error is measured in the $L^2(\Omega(T))$-norm.}
      \label{plot_ks2kt1_conv_both}
      \vspace*{-0.2cm}
    \end{figure}

The investigations in the first part of this subsection confirmed the numerical error to scale at least with $\mathcal{O}(\Delta t^{\ktsol + 1})$ in the mentioned norms. A well-known property of DG type discretisations is that of superconvergence, i.e. the convergence in time with a higher order (at specific time instances, e.g.\ nodes; cf. also \cite{AM89} for a superconvergence result for a CG method of the heat equation.). We investigate numerically whether our method shows such a behaviour in a two-stage procedure: In a first experiment, we perform only time-refinements for a fine mesh and a high order discretisation in space. This experiment allows to estimate the order of (super)-convergence in time, although asymptotically this order will be hard to reach exactly because a saturation from the spatial error will occur. Hence, we validate the first estimate for the superconvergence in a second step by a space-time-refinement study involving the suiting spatial discretisation order. The results for the first study, involving the choices of order $(\ktsol,\khsol)=(1,3), (2,5)$, $\klset = \khsol$ and a fixed mesh in space ($i_s = \mathrm{const}$) are shown in \cref{plot_ks3kt1_conv_time}.
Whilst in the space-time norm $L^2( L^2(\Omega(t)), 0, T)$ the estimate of \cref{eq:generalerrorbound} is confirmed, for $L^2(\Omega(T))$ we observe for the DG method and the CG method with $\ktsol = 2$, $\| u - u_h \|_{L^2(\Omega(T))} \approx \mathcal{O}(\Delta t^{\ktsol + 1.8})$, which in light of the aforementioned saturation effect suggests
\begin{equation}
 \| u - u_h \|_{L^2(\Omega(T))} = \mathcal{O}(\Delta t^{\ktsol + 2}), \textnormal{for DG, CG with } \ktsol = 2. \label{eq:superconverrorbound}
\end{equation}
For CG with $\ktsol = 1$, the estimate \cref{eq:generalerrorbound} is also observed to be sharp, there is no superconvergence. In a second step, we confirm \cref{eq:superconverrorbound} to hold with a space-time refinement study involving $(\ktsol,\khsol)=(1,2), (2,3)$, $\klset = \khsol$, see \cref{plot_ks2kt1_conv_both}.
We mention that theoretical results of findings similar to \cref{eq:superconverrorbound} exist for static domains, c.f. e.g. \cite{eriksson1985time,ern2016discontinuous}, and leave a further theoretical study of these phenomena in this particular setting for future research.

\subsection{Moving $n$-sphere}
The second example geometry is a moving circle in 2D/ sphere in 3D/ interval in 1D. These cases involve no change in deformation, but exemplify a moving domain setting. In dimension $d$, they are described by
\begin{subequations}
\begin{align}
 \rho(t) &= \frac{1}{\pi} \sin(2 \pi t), \quad \Omega = [-1,1] \times [-0.6,0.6]^{d-1}, \quad \mathbf{w} = (\frac{\partial \rho}{\partial t}, 0, \dots, 0), \\
 r(\mathbf{x}, t) &= [ (\mathbf{x}_1 - \rho(t))^2 + {\sum}_{i=2}^d \mathbf{x}_i^2]^{1/2} \quad \phi = r - 0.5.
\end{align}
\end{subequations}
The manufactured solution $u$ remains $u = \cos (\pi \frac{r}{r_0}) \cdot \sin (\pi t)$, leading to slightly different expressions for $f$ in each dimension. Moreover, $h= 0.5^{(i_s + 1)}, \Delta t = 2^{-i_t - 2}, \gamma_J = 0.05, \epsilon_f=1.1$.
\subsubsection{On 1D and 3D}
\begin{figure}
     \begin{tikzpicture}[scale=0.75]
    \begin{semilogyaxis}[ xlabel=$i$, ylabel=$\| u_h(T) - u (T)\|_{L^2(\Omega(T))}$, legend entries ={ $k=1$, $k=2$, $k=3$, $k=4$}, legend style={anchor=north,legend columns=1, draw=none, fill=none},legend pos =south west, x label style={at={(axis description cs:0.5,0.15)},anchor=center}, ymax=1]
     \addplot+[thick,mark options={scale=1.5}] table {num_exp/out/conv_moving_sphere_DG_ks1_kt1_both_nref6_gamma0.05_sm2.dat};
     \addplot+[thick,mark options={scale=1.5}] table {num_exp/out/conv_moving_sphere_DG_ks2_kt2_both_nref5_gamma0.05_sm2.dat};
     \addplot+[thick,mark options={scale=1.5}] table {num_exp/out/conv_moving_sphere_DG_ks3_kt3_both_nref4_gamma0.05_sm2.dat};
     \addplot+[thick,mark options={scale=1.5}] table {num_exp/out/conv_moving_sphere_DG_ks4_kt4_both_nref3_gamma0.05_sm2.dat};
     \addplot[gray, dashed, domain=0:5] {(1/2^(x-0.2)))^2};
     \addplot[gray, dashed, domain=0:4] {(1/2^(x+0.5)))^3};
     \addplot[gray, dashed, domain=0:3] {(1/2^(x+0.75)))^4};
     \addplot[gray, dashed, domain=0:2] {(1/2^(x+1.5)))^5};
    \end{semilogyaxis}
    \node[scale=0.75] at (4.5,5.25) {$O(h^{k+1})= O(\Delta t^{k+1})$};
     \draw[scale=0.75, gray, dash=on 2.25pt off 2.25pt phase 0pt, line width=0.4*0.75pt] (6.25/0.75,5.25/0.75) -- (6.8/0.75,5.25/0.75);
   \end{tikzpicture} \hspace{0.1cm}
  \tikzexternalexportnextfalse
     \begin{tikzpicture}[scale=0.75]
    \begin{semilogyaxis}[ xlabel=$i$, ylabel=$\| u_h(T) - u (T)\|_{L^2(\Omega(T))}$, legend entries ={ $k=1$, $k=2$, $k=3$, $k=4$}, legend style={legend columns=1, draw=none}, legend pos =south west, x label style={at={(axis description cs:0.8,0.15)},anchor=center}, ymin=5e-15]
     \addplot+[thick,mark options={scale=1.5}] table {num_exp/out/conv_moving_int_CG_ks1_kt1_both_nref12_gamma0.05_sm1.dat}; \label{p1}
     \addplot+[thick,mark options={scale=1.5}] table {num_exp/out/conv_moving_int_CG_ks2_kt2_both_nref11_gamma0.05_sm1.dat}; \label{p2}
     \addplot+[thick,mark options={scale=1.5}] table {num_exp/out/conv_moving_int_CG_ks3_kt3_both_nref10_gamma0.05_sm1.dat}; \label{p3}
     \addplot+[thick,mark options={scale=1.5}] table {num_exp/out/conv_moving_int_CG_ks4_kt4_both_nref9_gamma0.05_sm1.dat}; \label{p4}
     \addplot+[thick,mark options={scale=1.5}] table {num_exp/out/conv_moving_int_CG_ks5_kt5_both_nref8_gamma0.05_sm1.dat}; \label{p5}
     \addplot+[thick,mark options={scale=1.5}] table[skip coords between index={6}{7}] {num_exp/out/conv_moving_int_CG_ks6_kt6_both_nref7_gamma0.05_sm1.dat}; \label{p6}
     \addplot+[thick,mark options={scale=1.5}] table {num_exp/out/conv_moving_int_CG_ks7_kt7_both_nref6_gamma0.05_sm1.dat}; \label{p7}
     \addplot+[thick,mark options={scale=1.5}] table {num_exp/out/conv_moving_int_CG_ks8_kt8_both_nref5_gamma0.05_sm1.dat}; \label{p8}
     \addplot[gray, dashed, domain=0:11] {(1/2^(x-0.2)))^2}; \label{order}
     \addplot[gray, dashed, domain=0:10] {(1/2^(x+0.5)))^3};
     \addplot[gray, dashed, domain=0:9] {(1/2^(x+0.6)))^4};
     \addplot[gray, dashed, domain=0:8] {(1/2^(x+0.7)))^5};
     \addplot[gray, dashed, domain=0:7] {(1/2^(x+0.25)))^6};
     \addplot[gray, dashed, domain=0:5] {(1/2^(x+0.25)))^7};
     \addplot[gray, dashed, domain=0:5] {(1/2^(x+0.25)))^8};
     \addplot[gray, dashed, domain=0:4] {(1/2^(x+0.25)))^9};

     \node at (rel axis cs: 0.65,0.78) {\shortstack[r]{
 \ref{order} $O(h^{k+1})= O(\Delta t^{k+1})$ \\ \ref{p5} $k=5$ \\
\ref{p6} $k=6$ \\
\ref{p7} $k=7$ \\ \ref{p8} $k=8$}};
    \end{semilogyaxis}
   \end{tikzpicture} 
\vspace*{-0.7cm}
   \caption{Convergence of DG method for 3+1D (left) and the CG method for 1+1D (right).}
   \label{plot_dg_mov_sph_cg_mov_int_conv}
   \vspace*{-0.5cm}
    \end{figure}
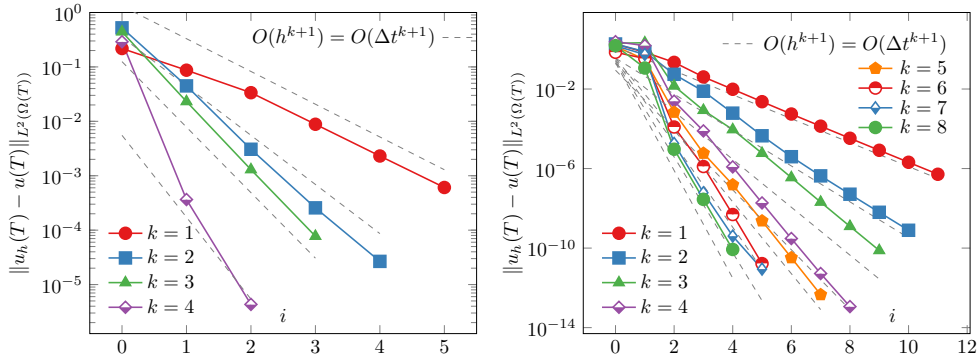

In this part, we want to show that the convergence results \cref{eq:generalerrorbound} also hold for spatially one- and three-dimensional problems. To this end, we apply the methods of consideration to the mentioned problems and study the decay of the numerical error for $\ksol=\klset$, $i_s=i_t$. In \cref{plot_dg_mov_sph_cg_mov_int_conv}, the DG method is applied in three spatial dimensions and the CG method is considered in 1D.
The respective choice of refinement levels and tested discretisation orders reflects the different computational demands. The results in the $L^2( L^2(\Omega(t)), 0, T)$-norm are equivalent and shown in 
the supplement, as well as results for the DG method in 1D and the CG method in 3D. In all cases, the convergence orders of \cref{eq:generalerrorbound} are confirmed.
\subsubsection{Non-zero entries of the matrices} \label{subsectnze}
For an evaluation of the computational costs, we compare the number of non-zero entries (\texttt{nze}s) of the discrete problem matrix of the respective methods. We fix the setting of the moving circle, $\ksol = \klset $, $i_s=2$. On each time-slice, after the system matrix is assembled, \texttt{nze}s are counted and the minimal and maximal value of this number over all slices is considered.

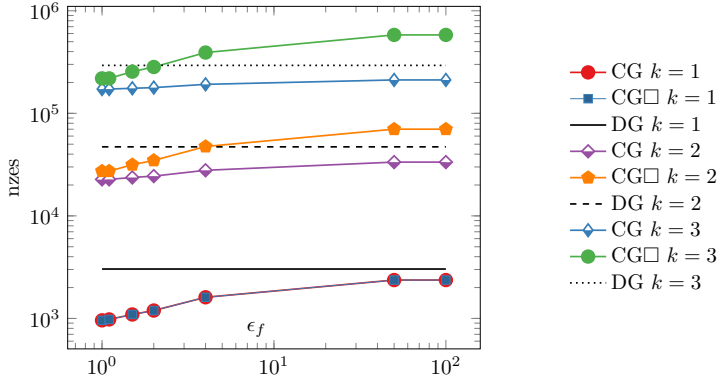
\begin{figure} \centering
     \begin{tikzpicture}[scale=0.8]
     \begin{loglogaxis}[ xlabel=$\epsilon_f$, ylabel={nzes}, legend entries ={CG $k = 1$, CG$\square$ $k = 1$, DG $k=1$, CG $k=2$, CG$\square$ $k=2$, DG $k=2$, CG $k=3$, CG$\square$ $k=3$, DG $k=3$, CG $k=4$, CGtg $k=4$, DG $k=4$}, legend style={at={(1.4,0.5)},anchor=center,legend columns=1, draw=none, legend cell align=left},x label style={at={(axis description cs:0.5,0.15)},anchor=east},y label style={at={(axis description cs:0.05,0.5)}}]
      \addplot+[mark options={scale=1.5},thick] table [x index = 0, y index =3] {num_exp/out/nze_moving_circle_CG_ks1_kt1_ktls1tp_gp0.dat};
      \addplot table [x index = 0, y index =3] {num_exp/out/nze_moving_circle_CG_ks1_kt1_ktls1tp_gp1.dat};
      \addplot[no markers, thick] table [x index = 0, y index =3] {num_exp/out/nze_moving_circle_DG_ks1_kt1_ktls1.dat};
      \addplot+[mark options={scale=1.5},thick] table [x index = 0, y index =3] {num_exp/out/nze_moving_circle_CG_ks2_kt2_ktls2tp_gp0.dat};
      \addplot+[mark options={scale=1.5},thick] table [x index = 0, y index =3] {num_exp/out/nze_moving_circle_CG_ks2_kt2_ktls2tp_gp1.dat};
      \addplot[no markers, dashed, thick] table [x index = 0, y index =3] {num_exp/out/nze_moving_circle_DG_ks2_kt2_ktls2.dat};
      \addplot+[mark options={scale=1.5},thick] table [x index = 0, y index =3] {num_exp/out/nze_moving_circle_CG_ks3_kt3_ktls3tp_gp0.dat}; 
      \addplot+[mark options={scale=1.5},thick] table [x index = 0, y index =3] {num_exp/out/nze_moving_circle_CG_ks3_kt3_ktls3tp_gp1.dat}; 
      \addplot[no markers, dotted, thick] table [x index = 0, y index =3] {num_exp/out/nze_moving_circle_DG_ks3_kt3_ktls3.dat};
     \end{loglogaxis}
   \end{tikzpicture} \vspace*{-0.3cm}
   \caption{Non-zero entries for DG and CG in comparison for different parameter choices for $\epsilon_f$. Each value represents the maximal number of non-zero entries over all time steps.}
      \label{nze_comp}\vspace*{-0.25cm}
    \end{figure}

For the CG method, in comparison to the DG method, \texttt{nze}s will depend on the factor $\epsilon_f$, as we evaluate in a first study. To this end, we choose $i_t=2$ and evaluate \texttt{nze}s for the CG method with $\epsilon_f \in \{ 1, 1.1, 1.5, \dots, 100\}$, and for comparison for the DG method, for the polynomial degrees $k=1,2,3$. The results are displayed in \cref{nze_comp}.
First, we observe that over all polynomial degrees the CG method leads to less non-zero entries as the DG method, whilst the relative difference is specifically strong for small orders $k$. This aligns with the elimination of one \dof~of the time finite element space. Second, this benefit of the CG method is reduced by a high choice of the factor $\epsilon_f$, as a larger domain of extension will induce more entries in the matrix. Third, the original CG method of \Cref{ssection:CG} performs better than the variant mentioned in \cref{rem:CGtg}. Hence, to exploit the computational benefits of the CG method fully, we suggest to use the non-tensor-product variant of the ghost penalty on a region as small as possible.

In a second study, we fix $\epsilon_f=1.1$ and compare the non-zero entries for all three methods for different time step sizes $\Delta t$ and different polynomial orders. The numbers are given in \cref{nze_study2}.
%
    \begin{table} \centering
     \small  
     \begin{tabular}{@{}r@{~}r@{ --}r@{~(}r@{)~}r@{ - }r@{~(}r@{)~}r@{ --}r@{~(}r@{)}} \toprule
      &\multicolumn{3}{c}{$i_t = 1$} 
      &\multicolumn{3}{c}{$i_t = 3$} 
      &\multicolumn{3}{c}{$i_t = 5$}\\
      &\multicolumn{2}{c}{\texttt{nze}s} &\multicolumn{1}{c}{\texttt{err}}     
      &\multicolumn{2}{c}{\texttt{nze}s} &\multicolumn{1}{c}{\texttt{err}}     
      &\multicolumn{2}{c}{\texttt{nze}s} &\multicolumn{1}{c}{\texttt{err}} 
      \\ \midrule
      \input num_exp/nze2_table_k1.tex \midrule
      \input num_exp/nze2_table_k3.tex \midrule
      \input num_exp/nze2_table_k5.tex
      \bottomrule
     \end{tabular}
     \caption{Range for non-zero entries (\texttt{nze}s) of system matrix (depends on time step) and absolute $L^2(T)$-error (\texttt{err}) in comparison for DG, CG and GCC. In this whole table, $i_s=2$.}
     \label{nze_study2}
    \vspace*{-0.9cm}
    \end{table}
Summarizing the results within each cell, the number of non-zero entries decreases from DG to CG and in turn to GCC. In particular for discretisation parameter $i_t=3,5$, this holds within a regime of similar absolute numerical errors. Within each row of the table, we note that smaller time steps decrease the numbers of non-zero entries for all methods, as the interface movement per time step decreases also. This effect is stronger for CG and GCC because of the extended regions of the ghost penalty stabilisation. For increasing polynomial order, we again observe a relatively smaller benefit of the CG and GCC method as compared to DG. We conclude that those methods are specifically interesting for orders such as $k=1,2,3$ from a computational point of view.
\subsubsection{Different time integration strategies} \label{tint_numexp}
In this part we investigate numerical differences between the topology-insensitive and the topology-preserving method of numerical quadrature in time suggested in \Cref{subsect_num_int}. For this purpose, we implement both approaches and compare the respective numerical errors. By and large, we observe competitive numerical results over most of the problem settings introduced so far. This is in agreement with the good results obtained in the literature, cf. \cite{Z18, HLZ16}. However, we observed differences in the following setup with almost no ghost-penalty stabilisation: Consider the DG method applied to the moving interval in 1D with a radius $R=0.505$, constant velocity $\rho = 0.5 t$, manufactured right-hand side for $u= (x-\rho +R)^2 \cdot (x - \rho -R)^2$, and polynomial order $\ksol = \klset = 4$. For this example, we fix $i_t = 0$ and refine in space. From $i_s=2$ on, the discrete interface is captured exactly so that the method can approximate the solution exactly as a polynomial of order 4.
 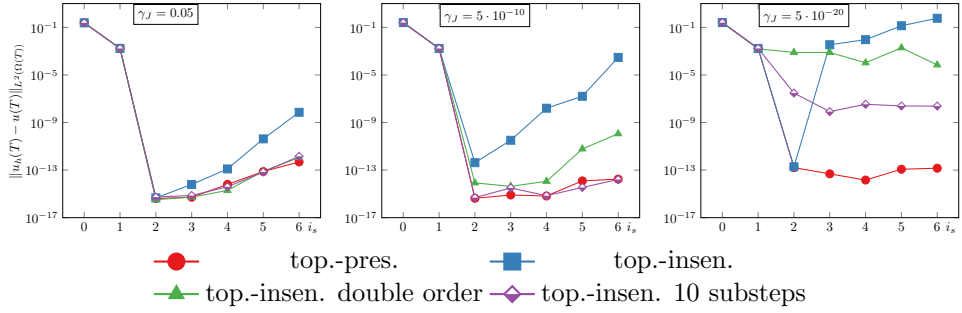
\begin{figure}
  \vspace*{-0.25cm}
  \tikzexternalexportnextfalse
   \begin{tikzpicture}[scale=0.5]
    \begin{semilogyaxis}[ xlabel=$i_s$, ylabel=$\| u_h(T) - u (T)\|_{L^2(\Omega(T))}$, 
      title style={at={(0.4,0.96)},draw=black, fill=white,anchor=north},
      x label style={at={(0.97,0.18)},anchor=west, below=5mm},      
      title={$\gamma_J = 0.05$}, legend entries ={top.-pres., top.-insen., top.-insen. double order, top.-insen. 10 substeps}, legend style={legend columns=2, draw=none}, legend to name=named, ymax=10, ymin=1e-17]
     \addplot+[thick, mark options={scale=1.5}] table {num_exp/out/conv_moving_int_sp_DG_ks4_kt4_space_nref7_gamma0.05_sm1_tint_proper_ro-7.dat};
     \addplot+[thick, mark options={scale=1.5}] table {num_exp/out/conv_moving_int_sp_DG_ks4_kt4_space_nref7_gamma0.05_sm1_tint_naive_ro-7.dat};
     \addplot+[thick, mark options={scale=1.5}] table {num_exp/out/conv_moving_int_sp_DG_ks4_kt4_space_nref7_gamma0.05_sm1_tint_naivedo_ro-7.dat};
     \addplot+[thick, mark options={scale=1.5}] table {num_exp/out/conv_moving_int_sp_DG_ks4_kt4_space_nref7_gamma0.05_sm1_tint_naiveS10_ro-7.dat};
    \end{semilogyaxis}
   \end{tikzpicture} \hspace{0.cm} \tikzexternalexportnextfalse
   \begin{tikzpicture}[scale=0.5]
    \begin{semilogyaxis}[ xlabel=$i_s$, 
      title style={at={(0.4,0.96)},draw=black, fill=white,anchor=north},
      x label style={at={(0.97,0.18)},anchor=west, below=5mm},      
      title={$\gamma_J = 5 \cdot 10^{-10}$}, ymax=10, ymin=1e-17]
     \addplot+[thick, mark options={scale=1.5}] table {num_exp/out/conv_moving_int_sp_DG_ks4_kt4_space_nref7_gamma5e-10_sm1_tint_proper_ro-7.dat};
     \addplot+[thick, mark options={scale=1.5}] table {num_exp/out/conv_moving_int_sp_DG_ks4_kt4_space_nref7_gamma5e-10_sm1_tint_naive_ro-7.dat};
     \addplot+[thick, mark options={scale=1.5}] table {num_exp/out/conv_moving_int_sp_DG_ks4_kt4_space_nref7_gamma5e-10_sm1_tint_naivedo_ro-7.dat};
     \addplot+[thick, mark options={scale=1.5}] table {num_exp/out/conv_moving_int_sp_DG_ks4_kt4_space_nref7_gamma5e-10_sm1_tint_naiveS10_ro-7.dat};
     \end{semilogyaxis}
   \end{tikzpicture} \hspace{0.cm}
   \begin{tikzpicture}[scale=0.5]
    \begin{semilogyaxis}[ xlabel=$i_s$, 
      title style={at={(0.4,0.96)},draw=black, fill=white,anchor=north},
      x label style={at={(0.97,0.18)},anchor=west, below=5mm},      
      title={$\gamma_J = 5 \cdot 10^{-20}$}, ymax=10, ymin=1e-17]
     \addplot+[thick, mark options={scale=1.5}] table {num_exp/out/conv_moving_int_sp_DG_ks4_kt4_space_nref7_gamma5e-20_sm1_tint_proper_ro-7.dat};
     \addplot+[thick, mark options={scale=1.5}] table {num_exp/out/conv_moving_int_sp_DG_ks4_kt4_space_nref7_gamma5e-20_sm1_tint_naive_ro-7.dat};
     \addplot+[thick, mark options={scale=1.5}] table {num_exp/out/conv_moving_int_sp_DG_ks4_kt4_space_nref7_gamma5e-20_sm1_tint_naivedo_ro-7.dat};
     \addplot+[thick, mark options={scale=1.5}] table {num_exp/out/conv_moving_int_sp_DG_ks4_kt4_space_nref7_gamma5e-20_sm1_tint_naiveS10_ro-7.dat};
     \end{semilogyaxis}
   \end{tikzpicture}\\  \centering
    \tikzexternalexportnextfalse
    \pgfplotslegendfromname{named}    \vspace*{-0.35cm}
\caption{Topology-insensitive and topology-preserving time integration in comparison.}
\label{naive_tint_comp}
\vspace*{-0.5cm}
\end{figure}
 The results of this study are shown in \cref{naive_tint_comp}. For the method with the topology-preserving time integration, we note that there appear only errors below a reasonable round-off error tolerance such as $10^{-12}$. Meanwhile, the method of topology-insensitive time integration leads to significant errors when very small stabilisation constants and mesh sizes are used. This is caused by cut configurations where an element is included in the active mesh because of a small cut, but the corresponding matrix entries vanish because of the numerical integration in time. In the presence of a moderate ghost penalty stabilisation, the \dofs~are still controlled by this stabilisation, otherwise stability issues appear. For an interpretation of this behaviour in terms of numerical analysis, we refer the reader to the comment made about $\gamma_J \to 0$ in \Cref{sssect:stab_choice}. Moreover, from a numerical analysis point of view we expect that the topology-preserving time integration is important for the geometrical consistency. In the numerical experiments however, we did not observe corresponding issues with the topology-insensitive time integration, at least not in combination with sufficient ghost penalty stabilization.
We are not sure if there is an underlying mechanism that can explain this behaviour or if the test case is still too mild and leave further investigations of this for future research.

In addition we compared the computational effort for both time integration strategies for the two-dimensional configuration of \cref{sect:kite}. As differences only occur on a small subset of elements in the computation domain, the overall computation time is only marginally affected by the choice of the integration strategy. 

We conclude that the topology-preserving time integration represents the discrete geometry properly, and hence is our method of choice, although the topology-insensitive time integration seems to perform equally well with sufficient stabilisation. Finally, we note that by doubling the order or introducing 10 evenly spaced substeps, the errors for the topology-insensitive time integration approach with degenerate stabilisation can be improved.\footnote{Here, doubling the order amounts to taking a plain Gaussian integration rule of order $4(k_t +1) =20$ instead of the plain rule of order $2 (k_t + 1) =10$ (c.f. Section 5). Moreover, introducing 10 evenly spaced substeps amounts to subdividing the interval $I_n$ into 10 subintervals of the same length and applying the Gaussian rule of order $2 (k_t +1) = 10$ to each of those.}

\subsection{Colliding circles}
Thirdly, we want to investigate the stability of the method in a topologically challenging example taken from \cite{LO_ESAIM_2019}: Two circles collide and separate again afterwards. This setting can be described by the following level set function:\vspace{-0.3cm}
\begin{equation}
 \phi(\mathbf{x},t) = \mathrm{min}(\| \mathbf{x} - s_1(t) \|, \| \mathbf{x} - s_2(t) \|) - R, \nonumber
 s_1(t) = (0, t - 3/4), s_2(t) = (0, 3/4 - t).
\end{equation}
The end time is set to $T = 3/2$, such that $\phi(\dots, 0) = \phi(\dots, T)$, the diffusion constant $\alpha = 0.1$, $R=0.5$, $\Omega=[-0.6,0.6] \times [-1.35,1.35]$, and the associated convection velocity is \vspace{-0.3cm}
\begin{equation}
 \mathbf{w}(x,y,t) = \begin{cases}
	              (0,-1)^T \quad \textnormal{if } (y > 0 \textnormal{ and } t \leq T/2) \textnormal{ or } (y < 0 \textnormal{ and } t > T/2) \\
		      (0,1)^T \quad \textnormal{if } (y \leq 0 \textnormal{ and } t \leq T/2) \textnormal{ or } (y > 0 \textnormal{ and } t > T/2).
                     \end{cases}
\end{equation}
Here, $u_0$ is defined to be $u_0(x,y) := \mathrm{sign}(y)$, i.e.\  1 in the upper circle and -1 in the lower one.

\begin{figure}
 \includegraphics[width=\textwidth]{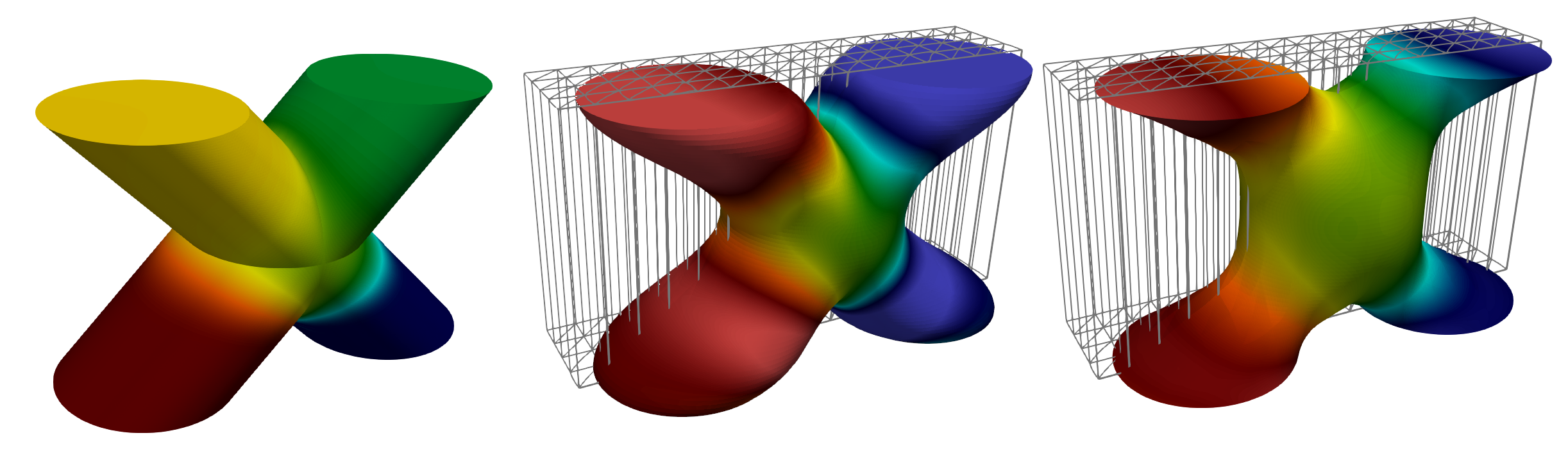}
 \vspace*{-0.8cm}
 \caption{Discrete solutions for the colliding circles test case. Left: Reference configuration. Center: Single-time-step run with $\ksol=4$. Right: Single-time-step run with $\ksol=2$.}
 \label{collid_circle_fig}
 \vspace*{-0.3cm}
\end{figure}

First, we investigate the numerical method with $\ksol = \klset = 4$, $i_s = 2$, $\Delta t= T/16$ to generate a reference solution, which is shown on the left-hand side of \cref{collid_circle_fig}. As one would expect physically, first the circles are transported with the convection field. After they merged, the diffusion successively evens out the concentration difference. By and large, the simulation does not show any numerical instability or large error. As a small detail, we notice that some numerical diffusion will appear shortly before the contact of the circles, when the numerically extended regions merge. This was also observed in \cite{LO_ESAIM_2019} for the time stepping method presented there. But for fine meshes or small stabilisation constants, this effect can be kept small. We conclude that the method performs well, despite a challenging topology change, without the need of any specific adaptation.

Second, we consider a more extreme version by setting $\Delta t = T$, so that the whole time interval will be solved in one time step. We show the numerical solutions for the parameter choices $\ksol = \klset = 2, 4$ in \cref{collid_circle_fig}. As the overall topology of the problem is still resolved properly, the method appears capable of handling topologically demanding configurations properly.
\section{Conclusion} \label{sect:concl}
In this final section, we want to summarise the findings of this work and highlight some questions that remain for future research. 

We introduced unfitted space-time methods for a 
problem posed on a moving domain. These are based on discontinuous Galerkin, continuous Galerkin, and higher order continuous Galerkin-Collocation discretisations in time. We further introduced techniques to achieve a high order geometry handling of the space-time domains based on isoparametric mappings of the background mesh. 
In extensive numerical experiments we discussed and compared several important properties of the proposed methods. Especially, the expectations of robustness and higher order convergence in space and time have been confirmed. Even in topologically challenging configurations, our methods perform reliably.
%
\noindent
These results motivate the following directions for future investigations: 

First, the methods have been investigated from a computational point of view in this paper. This suggests to investigate whether the numerically established convergence estimates can be proven from the point of view of numerical analysis. As far as the DG method is concerned, the master's theses \cite{preuss18} and \cite{heimann20} present such error estimates, although in the former an exact handling of the geometries is assumed, which is weakened in the latter in regard to the spatial approximation. 
This leaves a numerical analysis relative to the spatially and temporally discrete geometry an open task. Moreover, for the CG and GCC method no such analysis has been performed to the author's best knowledge.

Next, the numerical experiments w.r.t. the different numerical integration approaches suggest that the topology-insensitive time integration is also a reasonable and accurate approach
and the topology preserving time integration seems to be superfluous in the presence of ghost penalty stabilisation. This behavior is not fully understood and leaves room for future research.

Finally, let us mention that although we constructed a space-time mesh deformation that can be discontinuous between time slabs, it may be possible to construct a continuous mesh deformation that preserves the accuracy of the scheme. This would allow to remove the transfer operation $\Pi_u$ and hence simplify the method and potentially its analysis.  

\vspace*{-0.4cm} 

\bibliographystyle{siamplain}
\bibliography{st_article}
\end{document}


\maketitle

\appendix

\section{Specification of the CG$\square$ variant of the CG method}
The variant of the CG method discussed in \cref{rem:CGtg} is defined through the following modified definitions: $U_h^n = U_{h,0}^n + u_{\text{init}}$, with $u_{\text{init}}$ as in the original CG method and
\begin{align*}
 U_{h,0}^n :&= \{ v \in \Whcut{\ksol}{\EplOmn \times I_n} \circ (\Theta_h^{n, \mathrm{st}})^{-1} \, | \, v_+^{n-1} = 0 \} \\
V_h^{n} :&= \Whcut{(\khsol,\ktsol-1)}{\EplOmn \times I_n} \circ (\Theta_h^{n, \mathrm{st}})^{-1}.
\end{align*}
The discrete problem is then given as: Find $u = u_0 + u_{\text{init}}$ with $u_0 \in U_{h,0}^{n}$ s.t.
\begin{equation*}
B^{n}(u_0,v) + J^{n\square}(u_0,v)  = f^n(v) - B^{n}(u_{\text{init}},v) \quad \forall v \in V_h^{n}.
\end{equation*}

\section{Supplementary Plots}
In the following, we present results of numerical experiments which we use to draw our final conclusions, which are however very similar to related graphs in the original paper.

First, \cref{plot_dg_kite_conv_l2l2} and \cref{plot_cg_kite_conv_l2l2} contain evaluations of the numerical experiments of \cref{plot_dg_kite_conv} and \cref{plot_cg_kite_conv} in the $L^2( L^2(\Omega(t)), 0, T)$-norm.

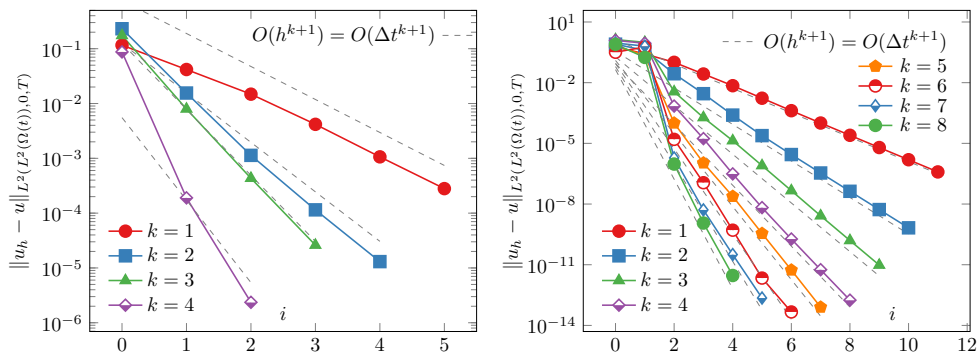
\begin{figure}
     \begin{tikzpicture}[scale=0.75]
    \begin{semilogyaxis}[ xlabel=$i$, ylabel=$\| u_h - u\|_{L^2( L^2(\Omega(t)), 0, T) }$, legend entries ={ $k=1$, $k=2$, $k=3$, $k=4$}, legend style={anchor=north,legend columns=1, draw=none, fill=none},legend pos =south west, x label style={at={(axis description cs:0.5,0.15)},anchor=center}, ymax=0.5]
     \addplot+[thick,mark options={scale=1.5}] table [x index = 0, y index =2] {num_exp/out/conv_moving_sphere_DG_ks1_kt1_both_nref6_gamma0.05_sm2.dat};
     \addplot+[thick,mark options={scale=1.5}] table [x index = 0, y index =2] {num_exp/out/conv_moving_sphere_DG_ks2_kt2_both_nref5_gamma0.05_sm2.dat};
     \addplot+[thick,mark options={scale=1.5}] table [x index = 0, y index =2] {num_exp/out/conv_moving_sphere_DG_ks3_kt3_both_nref4_gamma0.05_sm2.dat};
     \addplot+[thick,mark options={scale=1.5}] table [x index = 0, y index =2] {num_exp/out/conv_moving_sphere_DG_ks4_kt4_both_nref3_gamma0.05_sm2.dat};
     \addplot[gray, dashed, domain=0:5] {(1/2^(x+0.2)))^2};
     \addplot[gray, dashed, domain=0:4] {(1/2^(x+1)))^3};
     \addplot[gray, dashed, domain=0:3] {(1/2^(x+0.75)))^4};
     \addplot[gray, dashed, domain=0:2] {(1/2^(x+1.5)))^5};
    \end{semilogyaxis}
    \node[scale=0.75] at (4.5,5.25) {$O(h^{k+1})= O(\Delta t^{k+1})$};
     \draw[scale=0.75, gray, dash=on 2.25pt off 2.25pt phase 0pt, line width=0.4*0.75pt] (6.25/0.75,5.25/0.75) -- (6.8/0.75,5.25/0.75);
   \end{tikzpicture} \hspace{0.1cm}
  \tikzexternalexportnextfalse
     \begin{tikzpicture}[scale=0.75]
    \begin{semilogyaxis}[ xlabel=$i$, ylabel=$\| u_h - u\|_{L^2( L^2(\Omega(t)), 0, T) }$, legend entries ={ $k=1$, $k=2$, $k=3$, $k=4$}, legend style={legend columns=1, draw=none, fill=none}, legend pos =south west, x label style={at={(axis description cs:0.8,0.15)},anchor=center}, ymin=5e-15]
     \addplot+[thick,mark options={scale=1.5}] table [x index = 0, y index =2] {num_exp/out/conv_moving_int_CG_ks1_kt1_both_nref12_gamma0.05_sm1.dat};
     \addplot+[thick,mark options={scale=1.5}] table [x index = 0, y index =2] {num_exp/out/conv_moving_int_CG_ks2_kt2_both_nref11_gamma0.05_sm1.dat};
     \addplot+[thick,mark options={scale=1.5}] table [x index = 0, y index =2] {num_exp/out/conv_moving_int_CG_ks3_kt3_both_nref10_gamma0.05_sm1.dat};
     \addplot+[thick,mark options={scale=1.5}] table [x index = 0, y index =2] {num_exp/out/conv_moving_int_CG_ks4_kt4_both_nref9_gamma0.05_sm1.dat};
     \addplot+[thick,mark options={scale=1.5}] table [x index = 0, y index =2] {num_exp/out/conv_moving_int_CG_ks5_kt5_both_nref8_gamma0.05_sm1.dat}; \label{p5_l2l2}
     \addplot+[thick,mark options={scale=1.5}] table [x index = 0, y index =2] {num_exp/out/conv_moving_int_CG_ks6_kt6_both_nref7_gamma0.05_sm1.dat}; \label{p6_l2l2}
     \addplot+[thick,mark options={scale=1.5}] table [x index = 0, y index =2] {num_exp/out/conv_moving_int_CG_ks7_kt7_both_nref6_gamma0.05_sm1.dat}; \label{p7_l2l2}
     \addplot+[thick,mark options={scale=1.5}] table [x index = 0, y index =2] {num_exp/out/conv_moving_int_CG_ks8_kt8_both_nref5_gamma0.05_sm1.dat}; \label{p8_l2l2}
     \addplot[gray, dashed, domain=0:11] {(1/2^(x-0.2)))^2}; \label{order_l2l2}
     \addplot[gray, dashed, domain=0:10] {(1/2^(x+0.5)))^3};
     \addplot[gray, dashed, domain=0:9] {(1/2^(x+0.6)))^4};
     \addplot[gray, dashed, domain=0:8] {(1/2^(x+0.7)))^5};
     \addplot[gray, dashed, domain=0:7] {(1/2^(x+0.5)))^6};
     \addplot[gray, dashed, domain=0:6] {(1/2^(x+0.5)))^7};
     \addplot[gray, dashed, domain=0:5] {(1/2^(x+0.5)))^8};
     \addplot[gray, dashed, domain=0:4] {(1/2^(x+0.5)))^9};

     \node at (rel axis cs: 0.65,0.78) {\shortstack[r]{
 \ref{order_l2l2} $O(h^{k+1})= O(\Delta t^{k+1})$ \\ \ref{p5_l2l2} $k=5$ \\
\ref{p6_l2l2} $k=6$ \\
\ref{p7_l2l2} $k=7$ \\ \ref{p8_l2l2} $k=8$}};
    \end{semilogyaxis}
   \end{tikzpicture}
\vspace*{-0.2cm}
   \caption{Convergence of DG method for 3+1D (left) and the CG method for 1+1D (right), measured in the $L^2( L^2(\Omega(t)), 0, T)$-norm. This figure is the counterpart of \cref{plot_dg_mov_sph_cg_mov_int_conv}.}
   \label{plot_dg_mov_sph_cg_mov_int_conv_l2l2}
   \vspace*{-0.2cm}   
    \end{figure}

\begin{figure}
    \begin{tikzpicture}[scale=0.75]
    \begin{semilogyaxis}[ xlabel=$i$, ylabel=$\| u_h(T) - u (T)\|_{L^2(\Omega(T))}$, legend entries ={ $k=1$, $k=2$, $k=3$, $k=4$}, legend style={anchor=north,legend columns=1, draw=none},legend pos =south west, x label style={at={(axis description cs:0.5,0.15)},anchor=center}, ymax=1]
     \addplot+[thick,mark options={scale=1.5}] table [x index = 0, y index =1] {num_exp/out/conv_moving_sphere_CG_ks1_kt1_both_nref6_gamma0.05_sm2.dat};
     \addplot+[thick,mark options={scale=1.5}] table [x index = 0, y index =1] {num_exp/out/conv_moving_sphere_CG_ks2_kt2_both_nref5_gamma0.05_sm2.dat};
     \addplot+[thick,mark options={scale=1.5}] table [x index = 0, y index =1] {num_exp/out/conv_moving_sphere_CG_ks3_kt3_both_nref4_gamma0.05_sm2.dat};
     \addplot+[thick,mark options={scale=1.5}] table [x index = 0, y index =1] {num_exp/out/conv_moving_sphere_CG_ks4_kt4_both_nref3_gamma0.05_sm2.dat};
     \addplot[gray, dashed, domain=0:5] {(1/2^(x-0.2)))^2};
     \addplot[gray, dashed, domain=0:4] {(1/2^(x+0.75)))^3};
     \addplot[gray, dashed, domain=0:3] {(1/2^(x+0.75)))^4};
     \addplot[gray, dashed, domain=0:2] {(1/2^(x+1.25)))^5};
    \end{semilogyaxis}
    \node[scale=0.75] at (4.5,5.25) {$O(h^{k+1})= O(\Delta t^{k+1})$};
     \draw[scale=0.75, gray, dash=on 2.25pt off 2.25pt phase 0pt, line width=0.4*0.75pt] (6.25/0.75,5.25/0.75) -- (6.8/0.75,5.25/0.75);
   \end{tikzpicture} \hspace{0.1cm}
  \tikzexternalexportnextfalse
     \begin{tikzpicture}[scale=0.75]
    \begin{semilogyaxis}[ xlabel=$i$, ylabel=$\| u_h(T) - u(T)\|_{L^2(\Omega(T))}$, legend entries ={ $k=1$, $k=2$, $k=3$, $k=4$}, legend style={legend columns=1, draw=none, fill=none}, legend pos =south west, x label style={at={(axis description cs:0.8,0.15)},anchor=center}, ymin=5e-15]
     \addplot+[thick,mark options={scale=1.5}] table [x index = 0, y index =1] {num_exp/out/conv_moving_int_DG_ks1_kt1_both_nref12_gamma0.05_sm1.dat};
     \addplot+[thick,mark options={scale=1.5}] table [x index = 0, y index =1] {num_exp/out/conv_moving_int_DG_ks2_kt2_both_nref11_gamma0.05_sm1.dat};
     \addplot+[thick,mark options={scale=1.5}] table [x index = 0, y index =1] {num_exp/out/conv_moving_int_DG_ks3_kt3_both_nref10_gamma0.05_sm1.dat};
     \addplot+[thick,mark options={scale=1.5}] table [x index = 0, y index =1] {num_exp/out/conv_moving_int_DG_ks4_kt4_both_nref9_gamma0.05_sm1.dat};
     \addplot+[thick,mark options={scale=1.5}] table [x index = 0, y index =1] {num_exp/out/conv_moving_int_DG_ks5_kt5_both_nref8_gamma0.05_sm1.dat}; \label{p5_dg}
     \addplot+[thick,mark options={scale=1.5}] table [x index = 0, y index =1] {num_exp/out/conv_moving_int_DG_ks6_kt6_both_nref7_gamma0.05_sm1.dat}; \label{p6_dg}
     \addplot+[thick,mark options={scale=1.5}] table [x index = 0, y index =1] {num_exp/out/conv_moving_int_DG_ks7_kt7_both_nref6_gamma0.05_sm1.dat}; \label{p7_dg}
     \addplot+[thick,mark options={scale=1.5}] table [x index = 0, y index =1] {num_exp/out/conv_moving_int_DG_ks8_kt8_both_nref5_gamma0.05_sm1.dat}; \label{p8_dg}
     \addplot[gray, dashed, domain=0:11] {(1/2^(x-0.2)))^2}; \label{order_dg}
     \addplot[gray, dashed, domain=0:10] {(1/2^(x+0.5)))^3};
     \addplot[gray, dashed, domain=0:9] {(1/2^(x+0.6)))^4};
     \addplot[gray, dashed, domain=0:8] {(1/2^(x+0.7)))^5};
     \addplot[gray, dashed, domain=0:7] {(1/2^(x+0.25)))^6};
     \addplot[gray, dashed, domain=0:6] {(1/2^(x+0.25)))^7};
     \addplot[gray, dashed, domain=0:5] {(1/2^(x+0.25)))^8};
     \addplot[gray, dashed, domain=0:4] {(1/2^(x+0.25)))^9};

     \node at (rel axis cs: 0.65,0.78) {\shortstack[r]{
 \ref{order_dg} $O(h^{k+1})= O(\Delta t^{k+1})$ \\ \ref{p5_dg} $k=5$ \\
\ref{p6_dg} $k=6$ \\
\ref{p7_dg} $k=7$ \\ \ref{p8_dg} $k=8$}};
    \end{semilogyaxis}
   \end{tikzpicture}
   \vspace{0.0cm}
    
     \begin{tikzpicture}[scale=0.75]
    \begin{semilogyaxis}[ xlabel=$i$, ylabel=$\| u_h - u \|_{L^2( L^2(\Omega(t)), 0, T)}$, legend entries ={ $k=1$, $k=2$, $k=3$, $k=4$}, legend style={anchor=north,legend columns=1, draw=none, fill=none},legend pos =south west, x label style={at={(axis description cs:0.5,0.15)},anchor=center}, ymax=1]
     \addplot+[thick,mark options={scale=1.5}] table [x index = 0, y index =2] {num_exp/out/conv_moving_sphere_CG_ks1_kt1_both_nref6_gamma0.05_sm2.dat};
     \addplot+[thick,mark options={scale=1.5}] table [x index = 0, y index =2] {num_exp/out/conv_moving_sphere_CG_ks2_kt2_both_nref5_gamma0.05_sm2.dat};
     \addplot+[thick,mark options={scale=1.5}] table [x index = 0, y index =2] {num_exp/out/conv_moving_sphere_CG_ks3_kt3_both_nref4_gamma0.05_sm2.dat};
     \addplot+[thick,mark options={scale=1.5}] table [x index = 0, y index =2] {num_exp/out/conv_moving_sphere_CG_ks4_kt4_both_nref3_gamma0.05_sm2.dat};
     \addplot[gray, dashed, domain=0:5] {(1/2^(x+0.2)))^2};
     \addplot[gray, dashed, domain=0:4] {(1/2^(x+0.9)))^3};
     \addplot[gray, dashed, domain=0:3] {(1/2^(x+1.2)))^4};
     \addplot[gray, dashed, domain=0:2] {(1/2^(x+1.5)))^5};
    \end{semilogyaxis}
    \node[scale=0.75] at (4.5,5.25) {$O(h^{k+1})= O(\Delta t^{k+1})$};
     \draw[scale=0.75, gray, dash=on 2.25pt off 2.25pt phase 0pt, line width=0.4*0.75pt] (6.25/0.75,5.25/0.75) -- (6.8/0.75,5.25/0.75);
   \end{tikzpicture} \hspace{0.1cm}
  \tikzexternalexportnextfalse
     \begin{tikzpicture}[scale=0.75]
    \begin{semilogyaxis}[ xlabel=$i$, ylabel=$\| u_h - u\|_{L^2( L^2(\Omega(t)), 0, T) }$, legend entries ={ $k=1$, $k=2$, $k=3$, $k=4$}, legend style={legend columns=1, draw=none, fill=none}, legend pos =south west, x label style={at={(axis description cs:0.8,0.15)},anchor=center}, ymin=5e-15]
     \addplot+[thick,mark options={scale=1.5}] table [x index = 0, y index =2] {num_exp/out/conv_moving_int_DG_ks1_kt1_both_nref12_gamma0.05_sm1.dat};
     \addplot+[thick,mark options={scale=1.5}] table [x index = 0, y index =2] {num_exp/out/conv_moving_int_DG_ks2_kt2_both_nref11_gamma0.05_sm1.dat};
     \addplot+[thick,mark options={scale=1.5}] table [x index = 0, y index =2] {num_exp/out/conv_moving_int_DG_ks3_kt3_both_nref10_gamma0.05_sm1.dat};
     \addplot+[thick,mark options={scale=1.5}] table [x index = 0, y index =2] {num_exp/out/conv_moving_int_DG_ks4_kt4_both_nref9_gamma0.05_sm1.dat};
     \addplot+[thick,mark options={scale=1.5}] table [x index = 0, y index =2] {num_exp/out/conv_moving_int_DG_ks5_kt5_both_nref8_gamma0.05_sm1.dat}; \label{p5_dg_l2l2}
     \addplot+[thick,mark options={scale=1.5}] table [x index = 0, y index =2] {num_exp/out/conv_moving_int_DG_ks6_kt6_both_nref7_gamma0.05_sm1.dat}; \label{p6_dg_l2l2}
     \addplot+[thick,mark options={scale=1.5}] table [x index = 0, y index =2] {num_exp/out/conv_moving_int_DG_ks7_kt7_both_nref6_gamma0.05_sm1.dat}; \label{p7_dg_l2l2}
     \addplot+[thick,mark options={scale=1.5}] table [x index = 0, y index =2] {num_exp/out/conv_moving_int_DG_ks8_kt8_both_nref5_gamma0.05_sm1.dat}; \label{p8_dg_l2l2}
     \addplot[gray, dashed, domain=0:11] {(1/2^(x-0.2)))^2}; \label{order_dg_l2l2}
     \addplot[gray, dashed, domain=0:10] {(1/2^(x+0.5)))^3};
     \addplot[gray, dashed, domain=0:9] {(1/2^(x+0.6)))^4};
     \addplot[gray, dashed, domain=0:8] {(1/2^(x+0.7)))^5};
     \addplot[gray, dashed, domain=0:7] {(1/2^(x+0.5)))^6};
     \addplot[gray, dashed, domain=0:6] {(1/2^(x+0.5)))^7};
     \addplot[gray, dashed, domain=0:5] {(1/2^(x+0.5)))^8};
     \addplot[gray, dashed, domain=0:4] {(1/2^(x+0.5)))^9};

     \node at (rel axis cs: 0.65,0.78) {\shortstack[r]{
 \ref{order_dg_l2l2} $O(h^{k+1})= O(\Delta t^{k+1})$ \\ \ref{p5_dg_l2l2} $k=5$ \\
\ref{p6_dg_l2l2} $k=6$ \\
\ref{p7_dg_l2l2} $k=7$ \\ \ref{p8_dg_l2l2} $k=8$}};
    \end{semilogyaxis}
   \end{tikzpicture}
\vspace*{-0.2cm}
   \caption{Convergence of CG method for 3+1D (left) and the DG method for 1+1D (right), measured in the $L^2(\Omega(T))$-norm (top), and in the $L^2( L^2(\Omega(t)), 0, T)$-norm (bottom).}
   \label{plot_cg_mov_sph_dg_mov_int_conv}
   \vspace*{-0.2cm}   
    \end{figure}
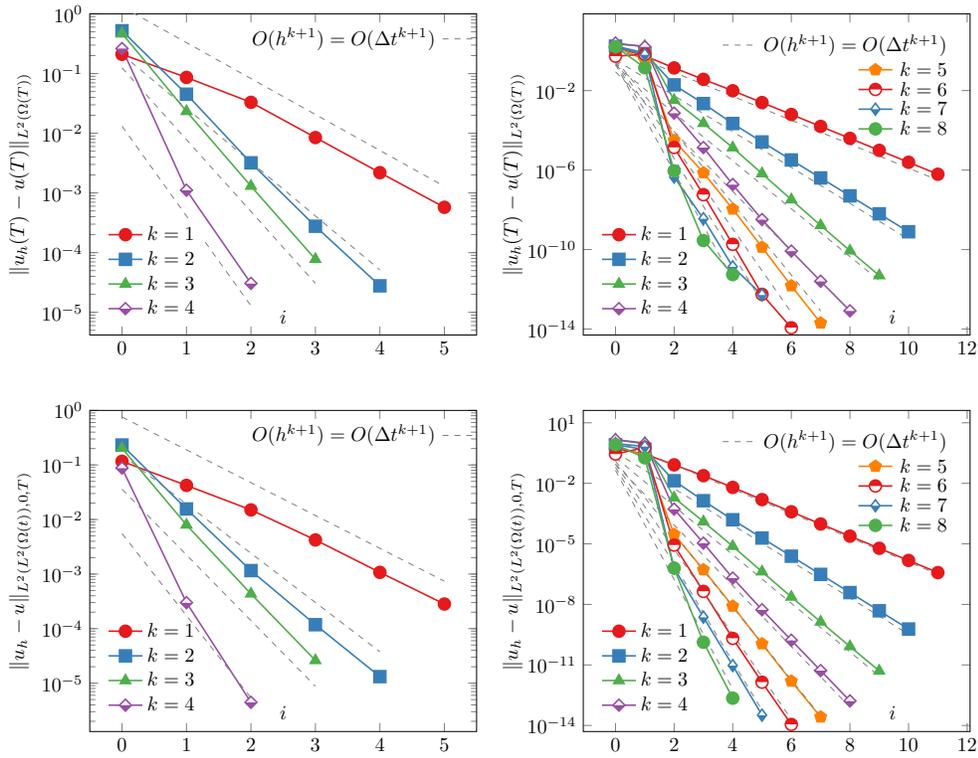

Second, \cref{plot_dg_mov_sph_cg_mov_int_conv_l2l2} and \cref{plot_cg_mov_sph_dg_mov_int_conv} supplement \cref{plot_dg_mov_sph_cg_mov_int_conv} from the paper by evaluating the $L^2( L^2(\Omega(t)), 0, T)$-norm and by investigating the other pairings of dimensionality (between 1+1D and 3+1D) and method (between DG and CG).
\begin{figure}
     \begin{tikzpicture}[scale=0.75]
    \begin{semilogyaxis}[ xlabel=$i$, ylabel=$\| u_h - u\|_{L^2( L^2(\Omega(t)), 0, T) }$, legend entries ={ $k=1$, $k=2$, $k=3$, $k=4$, $k=5$, $k=6$
    }, legend style={anchor=north,legend columns=1, draw=none, legend cell align=left, fill=none}, legend pos =south west,
    x label style={at={(axis description cs:0.5,0.15)},anchor=center},
    ymax=3, ymin=1e-11
     ]
     \addplot+[thick,mark options={scale=1.5}] table [x index = 0, y index =2] {num_exp/out/conv_kite_DG_ks1_kt1_both_nref8_gamma0.05_sm2.dat};
     \addplot+[thick,mark options={scale=1.5}] table [x index = 0, y index =2] {num_exp/out/conv_kite_DG_ks2_kt2_both_nref8_gamma0.05_sm2.dat};
     \addplot+[thick,mark options={scale=1.5}] table [x index = 0, y index =2] {num_exp/out/conv_kite_DG_ks3_kt3_both_nref8_gamma0.05_sm2.dat};
     \addplot+[thick,mark options={scale=1.5}] table [x index = 0, y index =2, skip coords between index={6}{7}] {num_exp/out/conv_kite_DG_ks4_kt4_both_nref7_gamma0.05_sm2.dat};
     \addplot+[thick,mark options={scale=1.5}] table [x index = 0, y index =2, skip coords between index={5}{6}] {num_exp/out/conv_kite_DG_ks5_kt5_both_nref6_gamma0.05_sm2.dat};
     \addplot+[thick,mark options={scale=1.5}] table [x index = 0, y index =2, skip coords between index={4}{5}] {num_exp/out/conv_kite_DG_ks6_kt6_both_nref5_gamma0.05_sm2.dat};
     \addplot[gray, dashed, domain=0:7] {(1/2^(x-0.2)))^2};
     \addplot[gray, dashed, domain=0:7] {(1/2^(x+0.5)))^3};
     \addplot[gray, dashed, domain=0:7] {(1/2^(x+1)))^4};
     \addplot[gray, dashed, domain=0:5] {(1/2^(x+1.25)))^5};
     \addplot[gray, dashed, domain=0:4] {(1/2^(x+1.5)))^6};
     \addplot[gray, dashed, domain=0:3] {(1/2^(x+1.5)))^7};
    \end{semilogyaxis}
     \node[scale=0.75] at (4.5,5.25) {$O(h^{k+1})= O(\Delta t^{k+1})$};
     \draw[scale=0.75, gray, dash=on 2.25pt off 2.25pt phase 0pt, line width=0.4*0.75pt] (6.25/0.75,5.25/0.75) -- (6.8/0.75,5.25/0.75);
   \end{tikzpicture} \hspace{0.2cm}
     \begin{tikzpicture}[scale=0.75]
    \begin{semilogyaxis}[ xlabel=$i$, ylabel=$\| u_h - u\|_{L^2( L^2(\Omega(t)), 0, T) }$, legend entries ={$\gamma_J = 5 \cdot 10^4$,$\gamma_J = 5$, $\gamma_J = 5 \cdot 10^{-2}$,$\gamma_J = 5 \cdot 10^{-4}$}, legend style={anchor=north,legend columns=1, draw=none, fill=none},legend pos =south west, x label style={at={(axis description cs:0.6,0.15)},anchor=center}]
     \addplot+[thick,mark options={scale=1.5}] table [x index = 0, y index =2] {num_exp/out/conv_kite_DG_ks4_kt4_both_nref5_gamma50000.0_sm2.dat};
     \addplot+[thick,mark options={scale=1.5}] table [x index = 0, y index =2] {num_exp/out/conv_kite_DG_ks4_kt4_both_nref5_gamma5.0_sm2.dat};
     \addplot+[thick,mark options={scale=1.5}] table [x index = 0, y index =2] {num_exp/out/conv_kite_DG_ks4_kt4_both_nref5_gamma0.05_sm2.dat};
     \addplot+[thick,mark options={scale=1.5}] table [x index = 0, y index =2] {num_exp/out/conv_kite_DG_ks4_kt4_both_nref5_gamma0.0005_sm2.dat};
     \addplot+[thick,mark options={scale=1.5}] table [x index = 0, y index =2] {num_exp/out/conv_kite_DG_ks4_kt4_both_nref5_gamma5e-10_sm2.dat}; \label{gamma5e-10B}
     \addplot+[thick,mark options={scale=1.5}] table [x index = 0, y index =2] {num_exp/out/conv_kite_DG_ks4_kt4_both_nref3_gamma0.0_sm2.dat}; \label{gamma0B}
     \addplot[gray, dashed, domain=0:4] {(1/2^(x+1.3)))^5};

               \node at (rel axis cs: 0.75,0.8) {\shortstack[r]{
$\gamma_J = 5 \cdot 10^{-10}$ \ref{gamma5e-10B} \\ $\gamma_J = 0$ \ref{gamma0B}}};
    \end{semilogyaxis}
              \node[scale=0.75] at (4.5,5.25) {$O(h^5)=O(\Delta t^5)$};
     \draw[scale=0.75, gray, dash=on 2.25pt off 2.25pt phase 0pt, line width=0.4*0.75pt] (6./0.75,5.25/0.75) -- (6.55/0.75,5.25/0.75);
   \end{tikzpicture}
    \vspace*{-0.6cm}
   \caption{Left: Convergence of the higher order DG method for the kite case in the $L^2( L^2(\Omega(t)), 0, T)$-norm. Right: Influence of the choice of the stabilisation parameter $\gamma$ on the convergence behaviour of the DG method for $\ksol = \klset = 4$, measured in the $L^2( L^2(\Omega(t)), 0, T)$. This figure is the counterpart of \cref{plot_dg_kite_conv}.}
   \label{plot_dg_kite_conv_l2l2}
    \end{figure}
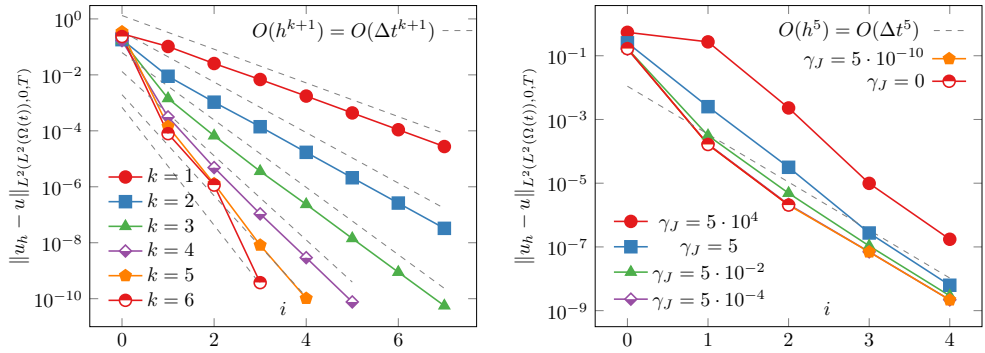
        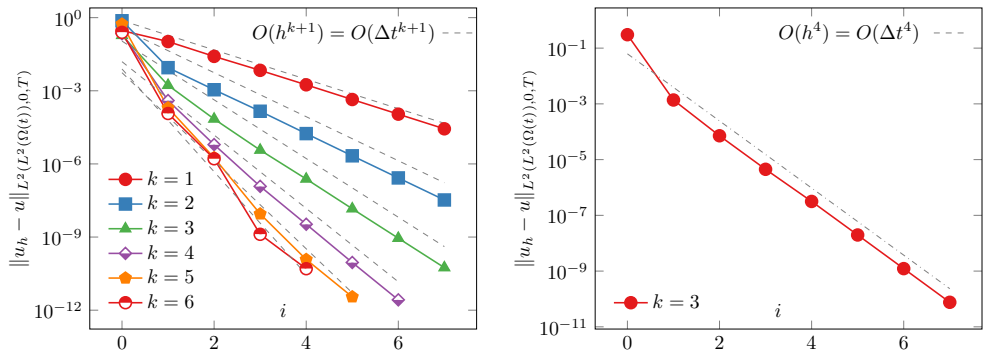
\begin{figure}
     \begin{tikzpicture}[scale=0.75]
    \begin{semilogyaxis}[ xlabel=$i$, ylabel=$\| u_h - u\|_{L^2( L^2(\Omega(t)), 0, T) }$, legend entries ={ $k=1$, $k=2$, $k=3$, $k=4$, $k=5$, $k=6$}, legend style={anchor=north,legend columns=1, draw=none, legend cell align=left, fill=none},legend pos =south west,
    x label style={at={(axis description cs:0.5,0.15)},anchor=center},
     ymax =2.5
    ]
     \addplot+[thick,mark options={scale=1.5}] table [x index = 0, y index =2] {num_exp/out/conv_kite_CG_ks1_kt1_both_nref8_gamma0.05_sm2.dat};
     \addplot+[thick,mark options={scale=1.5}] table [x index = 0, y index =2] {num_exp/out/conv_kite_CG_ks2_kt2_both_nref8_gamma0.05_sm2.dat};
     \addplot+[thick,mark options={scale=1.5}] table [x index = 0, y index =2] {num_exp/out/conv_kite_CG_ks3_kt3_both_nref8_gamma0.05_sm2.dat};
     \addplot+[thick,mark options={scale=1.5}] table [x index = 0, y index =2] {num_exp/out/conv_kite_CG_ks4_kt4_both_nref7_gamma0.05_sm2.dat};
     \addplot+[thick,mark options={scale=1.5}] table [x index = 0, y index =2] {num_exp/out/conv_kite_CG_ks5_kt5_both_nref6_gamma0.05_sm2.dat};
     \addplot+[thick,mark options={scale=1.5}] table [x index = 0, y index =2] {num_exp/out/conv_kite_CG_ks6_kt6_both_nref5_gamma0.05_sm2.dat};
     \addplot[gray, dashed, domain=0:7] {(1/2^(x+0.2)))^2};
     \addplot[gray, dashed, domain=0:7] {(1/2^(x+0.5)))^3};
     \addplot[gray, dashed, domain=0:7] {(1/2^(x+0.8)))^4};
     \addplot[gray, dashed, domain=0:6] {(1/2^(x+1.2)))^5};
     \addplot[gray, dashed, domain=0:5] {(1/2^(x+1.25)))^6};
     \addplot[gray, dashed, domain=0:4] {(1/2^(x+1)))^7};
    \end{semilogyaxis}
              \node[scale=0.75] at (4.5,5.25) {$O(h^{k+1})= O(\Delta t^{k+1})$};
     \draw[scale=0.75, gray, dash=on 2.25pt off 2.25pt phase 0pt, line width=0.4*0.75pt] (6.25/0.75,5.25/0.75) -- (6.8/0.75,5.25/0.75);
   \end{tikzpicture} \hspace{0.2cm}
     \begin{tikzpicture}[scale=0.75]
    \begin{semilogyaxis}[ xlabel=$i$, ylabel=$\| u_h - u\|_{L^2( L^2(\Omega(t)), 0, T) }$, legend entries ={ $k=3$}, legend style={anchor=north,legend columns=1, draw=none}, legend pos =south west, x label style={at={(axis description cs:0.5,0.15)},anchor=east}]
     \addplot+[thick,mark options={scale=1.5}] table[x index = 0, y index =2] {num_exp/out/conv_kite_GCC_ks3_kt3_both_nref8_gamma0.05_sm2.dat};
     \addplot[gray, dashdotted, domain=0:7] {(1/2^(x+1.)))^4};
    \end{semilogyaxis}
         \node[scale=0.75] at (4.5,5.25) {$O(h^4)=O(\Delta t^4)$};
     \draw[scale=0.75, gray, dash=on 2.25pt off 2.25pt phase 0pt, line width=0.4*0.75pt] (6./0.75,5.25/0.75) -- (6.55/0.75,5.25/0.75);
   \end{tikzpicture} 
     \vspace*{-0.6cm}
   \caption{Convergence of the higher order CG (left) and GCC (right) methods for the kite case, measured in the $L^2( L^2(\Omega(t)), 0, T)$-norm. This figure is the counterpart of \cref{plot_cg_kite_conv}.}
      \label{plot_cg_kite_conv_l2l2}
    \end{figure}

%% file: st_shared.tex
\usepackage{lipsum}
\usepackage{amsfonts}
\usepackage{amsmath}
\usepackage{amssymb}
\usepackage{graphicx}
\usepackage[caption=false]{subfig}
\usepackage{epstopdf}
\usepackage{algorithmic}

\usepackage{changes}
\ifpdf
  \DeclareGraphicsExtensions{.eps,.pdf,.png,.jpg}
\else
  \DeclareGraphicsExtensions{.eps}
\fi

\newsiamremark{remark}{Remark}
\newsiamremark{assumption}{Assumption}
\newsiamremark{hypothesis}{Hypothesis}
\crefname{hypothesis}{Hypothesis}{Hypotheses}
\newsiamthm{claim}{Claim}

\headers{Geometrically Unfitted Space-Time Methods}{F. Heimann, C. Lehrenfeld, J. Preu\ss{}}

\title{Geometrically Higher Order Unfitted Space-Time Methods for PDEs on Moving Domains\thanks{Revision 2
}}

\author{Fabian Heimann
\and Christoph Lehrenfeld
\and Janosch Preu\ss{}\thanks{Institut f\"ur Numerische und Angewandte Mathematik, Universit\"at G\"ottingen, 
(\texttt{\{f.heimann, lehrenfeld\}@math.uni-goettingen.de}, \texttt{j.preuss@ucl.ac.uk}).}}

\usepackage{amsopn}

\usepackage[]{hyperref}
\usepackage[]{cleveref}

\newcommand{\vertiii}[1]{{\left\vert\kern-0.25ex\left\vert\kern-0.25ex\left\vert #1 
\right\vert\kern-0.25ex\right\vert\kern-0.25ex\right\vert}}

\usepackage{tikz}
\usepackage{tikz-3dplot}
\usetikzlibrary{positioning,shapes,arrows,calc, arrows.meta, patterns}

\usepackage{pgfplots}
\usepgfplotslibrary{colorbrewer,groupplots}

\pgfplotsset{
cycle list/Set1-5,
cycle multiindex* list={
mark list*\nextlist
Set1-5\nextlist
},
}

\usepackage{shellesc}
\usetikzlibrary{external}

\usepackage{booktabs}

\usepackage{graphics}

    \usepackage{mathtools}

\newcommand{\dof}{\hyperlink{def:dof}{\texttt{dof}}}
\newcommand{\dofs}{\hyperlink{def:dof}{\texttt{dof}}s}
\newcommand{\pdofs}{\texttt{dof}s}  

\newcommand{\jump}[1]{[ #1 ]}

\newcommand{\vect}[1]{\mathbf{#1}}
\newcommand{\Q}{\hyperlink{def:Q}{Q}}
\newcommand{\Hoz}{\hyperref[def:Hoz]{H^{1,0}}}
\newcommand{\Hmoz}{\hyperlink{def:Hmoz}{H^{-1,0}}}
\newcommand{\In}[1]{\hyperlink{def:In}{I_{#1}}}
\newcommand{\Qn}[1]{\hyperlink{def:Qn}{Q^{#1}}}
\newcommand{\mysquare}{\square}

\newcommand{\Osn}[1]{\hyperlink{def:Osn}{\Omega_{{\mysquare}}^{#1}}}

\newcommand{\philins}{\hyperlink{def:philins}{\phi^{\text{lin}}}}
\newcommand{\Olins}{\hyperlink{def:Olins}{\Omega^{\text{lin}}}}
\newcommand{\Piu}{\hyperlink{def:Piu}{\Pi_u}}

\newcommand{\tOmega}{\hyperlink{def:tOmega}{\tilde{\Omega}}}
\newcommand{\Th}{\hyperlink{def:Th}{\mathcal{T}_h}}
\newcommand{\Vhks}[1]{\hyperlink{def:Vhks}{V_h^{#1}}}
\newcommand{\ThnG}[1]{\hyperref[eq:new_st_isoparam_regions]{\mathcal{T}_{h,#1}^{\Gamma}}}

\newcommand{\Oht}[1]{\hyperref[eq:introQhn]{{\Omega}^{h}\!(#1)}}

\newcommand{\Osnp}[1]{\hyperlink{def:Osnp}{\Omega_{{\mysquare+}}^{#1}}}
\newcommand{\EOmn}{\hyperlink{def:EOmn}{\mathcal{E}(\Omega^{\text{lin}, n})}}
\newcommand{\EQlinn}{\hyperlink{def:EQlinn}{\mathcal{E}(Q^{\text{lin}, n})}}
\newcommand{\IOmn}{\hyperlink{def:IOmn}{\mathcal{I}(\Omega^{\text{lin}, n})}}

\newcommand{\EplOmn}{\hyperlink{def:EplOmn}{\mathcal{E}^{+}(\Omega^{\text{lin}, n})}}

\newcommand{\Thehst}{\hyperlink{def:Thehst}{\Theta_h^{\text{st}}}}
\newcommand{\Qhn}{\hyperlink{def:Qhn}{Q^{h,n}}}
\newcommand{\Qlinn}{\hyperlink{def:Qlinn}{Q^{\text{lin}, n}}}
\newcommand{\Omlin}[1]{\hyperlink{def:Omlin}{\Omega^{\text{lin}}(#1)}}
\newcommand{\bgVhks}[1]{{W_h^{#1}}}
\newcommand{\bgVhnks}[2]{\hyperlink{def:bgVhnks}{W_h^{#1,#2}}}
\newcommand{\Whcut}[2]{\hyperlink{def:Whcut}{W_{h, \text{cut}}^{n,#1}(#2)}} 

\newcommand{\ext}{\hyperlink{def:ext}{\mathcal{E}}}

\newcommand{\khsol}{k_{s}}
\newcommand{\ktsol}{k_{t}}
\newcommand{\ksol}{k}
\newcommand{\khlset}{q_{s}}
\newcommand{\ktlset}{q_{t}}
\newcommand{\klset}{q}
\newcommand{\ktreg}{m_r}
\newcommand{\ktcol}{m_c}
\newcommand{\kttest}{m_t}


%% file: tikz_mesh.tex
\draw[backgrtrig] ( -1.05 , -1.05 ) -- ( -0.6166666666666668 , -1.05 ) -- ( -1.05 , -0.6300000000000003 ) -- cycle;
\draw[backgrtrig] ( 0.6833333333333327 , -1.05 ) -- ( 1.116666666666666 , -1.05 ) -- ( 0.8212531032817161 , -0.668206736929906 ) -- cycle;
\draw[backgrtrig] ( 1.55 , -1.05 ) -- ( 1.55 , -0.63 ) -- ( 1.116666666666666 , -1.05 ) -- cycle;
\draw[backgrtrig] ( 1.55 , -0.63 ) -- ( 1.55 , -0.21000000000000008 ) -- ( 1.1578746161954439 , -0.44698643381097225 ) -- cycle;
\draw[backgrtrig] ( 1.55 , -0.21000000000000008 ) -- ( 1.55 , 0.20999999999999996 ) -- ( 1.1415183492371044 , 0.016373332214704023 ) -- cycle;
\draw[backgrtrig] ( 1.55 , -0.21000000000000008 ) -- ( 1.1415183492371044 , 0.016373332214704023 ) -- ( 1.1578746161954439 , -0.44698643381097225 ) -- cycle;
\draw[backgrtrig] ( 1.55 , 0.20999999999999996 ) -- ( 1.55 , 0.6300000000000003 ) -- ( 1.1792544350003562 , 0.4754657340609835 ) -- cycle;
\draw[backgrtrig] ( 1.55 , 0.20999999999999996 ) -- ( 1.1792544350003562 , 0.4754657340609835 ) -- ( 1.1415183492371044 , 0.016373332214704023 ) -- cycle;
\draw[backgrtrig] ( 1.55 , 1.05 ) -- ( 1.1166666666666667 , 1.05 ) -- ( 1.55 , 0.6300000000000003 ) -- cycle;
\draw[backgrtrig] ( 1.1166666666666667 , 1.05 ) -- ( 0.6833333333333337 , 1.05 ) -- ( 0.8537612867784096 , 0.7056698463200347 ) -- cycle;
\draw[backgrtrig] ( -1.05 , 1.05 ) -- ( -1.05 , 0.63 ) -- ( -0.616666666666666 , 1.05 ) -- cycle;
\draw[backgrtrig] ( 1.116666666666666 , -1.05 ) -- ( 1.55 , -0.63 ) -- ( 1.1578746161954439 , -0.44698643381097225 ) -- cycle;
\draw[backgrtrig] ( 1.55 , 0.6300000000000003 ) -- ( 1.1166666666666667 , 1.05 ) -- ( 1.1792544350003562 , 0.4754657340609835 ) -- cycle;
\draw[backgrtrig] ( 1.116666666666666 , -1.05 ) -- ( 1.1578746161954439 , -0.44698643381097225 ) -- ( 0.8212531032817161 , -0.668206736929906 ) -- cycle;
\draw[backgrtrig] ( 1.1166666666666667 , 1.05 ) -- ( 0.8537612867784096 , 0.7056698463200347 ) -- ( 1.1792544350003562 , 0.4754657340609835 ) -- cycle;
\draw[meshtrig] ( -0.6166666666666668 , -1.05 ) -- ( -0.18333333333333368 , -1.05 ) -- ( -0.46573529073056247 , -0.6835123756350594 ) -- cycle;
\draw[meshtrig] ( -0.18333333333333368 , -1.05 ) -- ( 0.24999999999999956 , -1.05 ) -- ( -0.016375937719542086 , -0.646911115601791 ) -- cycle;
\draw[meshtrig] ( -0.18333333333333368 , -1.05 ) -- ( -0.016375937719542086 , -0.646911115601791 ) -- ( -0.46573529073056247 , -0.6835123756350594 ) -- cycle;
\draw[meshtrig] ( 0.24999999999999956 , -1.05 ) -- ( 0.6833333333333327 , -1.05 ) -- ( 0.42125922650351183 , -0.6311516272278445 ) -- cycle;
\draw[meshtrig] ( 0.24999999999999956 , -1.05 ) -- ( 0.42125922650351183 , -0.6311516272278445 ) -- ( -0.016375937719542086 , -0.646911115601791 ) -- cycle;
\draw[meshtrig] ( 0.6833333333333327 , -1.05 ) -- ( 0.8212531032817161 , -0.668206736929906 ) -- ( 0.42125922650351183 , -0.6311516272278445 ) -- cycle;
\draw[meshtrig] ( 0.6833333333333337 , 1.05 ) -- ( 0.25000000000000044 , 1.05 ) -- ( 0.4772704753189263 , 0.6599589103990077 ) -- cycle;
\draw[meshtrig] ( 0.6833333333333337 , 1.05 ) -- ( 0.4772704753189263 , 0.6599589103990077 ) -- ( 0.8537612867784096 , 0.7056698463200347 ) -- cycle;
\draw[meshtrig] ( 0.25000000000000044 , 1.05 ) -- ( -0.18333333333333268 , 1.05 ) -- ( 0.06406477184154608 , 0.6423728544761721 ) -- cycle;
\draw[meshtrig] ( 0.25000000000000044 , 1.05 ) -- ( 0.06406477184154608 , 0.6423728544761721 ) -- ( 0.4772704753189263 , 0.6599589103990077 ) -- cycle;
\draw[meshtrig] ( -0.18333333333333268 , 1.05 ) -- ( -0.616666666666666 , 1.05 ) -- ( -0.32685022560532856 , 0.6762126992658455 ) -- cycle;
\draw[meshtrig] ( -0.18333333333333268 , 1.05 ) -- ( -0.32685022560532856 , 0.6762126992658455 ) -- ( 0.06406477184154608 , 0.6423728544761721 ) -- cycle;
\draw[meshtrig] ( -1.05 , 0.63 ) -- ( -1.05 , 0.21000000000000008 ) -- ( -0.662003809555919 , 0.4617132239118765 ) -- cycle;
\draw[meshtrig] ( -1.05 , 0.21000000000000008 ) -- ( -1.05 , -0.20999999999999996 ) -- ( -0.650560958288211 , 0.02488343392998438 ) -- cycle;
\draw[meshtrig] ( -1.05 , 0.21000000000000008 ) -- ( -0.650560958288211 , 0.02488343392998438 ) -- ( -0.662003809555919 , 0.4617132239118765 ) -- cycle;
\draw[meshtrig] ( -0.6166666666666668 , -1.05 ) -- ( -0.46573529073056247 , -0.6835123756350594 ) -- ( -1.05 , -0.6300000000000003 ) -- cycle;
\draw[meshtrig] ( -0.616666666666666 , 1.05 ) -- ( -1.05 , 0.63 ) -- ( -0.662003809555919 , 0.4617132239118765 ) -- cycle;
\draw[meshtrig] ( -1.05 , -0.20999999999999996 ) -- ( -1.05 , -0.6300000000000003 ) -- ( -0.6971923744228338 , -0.36057467183152875 ) -- cycle;
\draw[meshtrig] ( -1.05 , -0.20999999999999996 ) -- ( -0.6971923744228338 , -0.36057467183152875 ) -- ( -0.650560958288211 , 0.02488343392998438 ) -- cycle;
\draw[meshtrig] ( -0.016375937719542086 , -0.646911115601791 ) -- ( 0.42125922650351183 , -0.6311516272278445 ) -- ( 0.16139908952035206 , -0.2112520920264157 ) -- cycle;
\draw[meshtrig] ( -0.46573529073056247 , -0.6835123756350594 ) -- ( -0.016375937719542086 , -0.646911115601791 ) -- ( -0.295081547266559 , -0.26766310993246095 ) -- cycle;
\draw[meshtrig] ( -0.016375937719542086 , -0.646911115601791 ) -- ( 0.16139908952035206 , -0.2112520920264157 ) -- ( -0.295081547266559 , -0.26766310993246095 ) -- cycle;
\draw[meshtrig] ( 1.1415183492371044 , 0.016373332214704023 ) -- ( 1.1792544350003562 , 0.4754657340609835 ) -- ( 0.7698780756308289 , 0.2979251021457717 ) -- cycle;
\draw[meshtrig] ( 1.1578746161954439 , -0.44698643381097225 ) -- ( 1.1415183492371044 , 0.016373332214704023 ) -- ( 0.6829318808011956 , -0.1985703945312043 ) -- cycle;
\draw[meshtrig] ( 1.1415183492371044 , 0.016373332214704023 ) -- ( 0.7698780756308289 , 0.2979251021457717 ) -- ( 0.6829318808011956 , -0.1985703945312043 ) -- cycle;
\draw[meshtrig] ( -0.616666666666666 , 1.05 ) -- ( -0.662003809555919 , 0.4617132239118765 ) -- ( -0.32685022560532856 , 0.6762126992658455 ) -- cycle;
\draw[meshtrig] ( 0.06406477184154608 , 0.6423728544761721 ) -- ( -0.32685022560532856 , 0.6762126992658455 ) -- ( -0.19510386774107463 , 0.2177653724852468 ) -- cycle;
\draw[meshtrig] ( 0.4772704753189263 , 0.6599589103990077 ) -- ( 0.06406477184154608 , 0.6423728544761721 ) -- ( 0.31947279823770486 , 0.2274752592073914 ) -- cycle;
\draw[meshtrig] ( 0.06406477184154608 , 0.6423728544761721 ) -- ( -0.19510386774107463 , 0.2177653724852468 ) -- ( 0.31947279823770486 , 0.2274752592073914 ) -- cycle;
\draw[meshtrig] ( -1.05 , -0.6300000000000003 ) -- ( -0.46573529073056247 , -0.6835123756350594 ) -- ( -0.6971923744228338 , -0.36057467183152875 ) -- cycle;
\draw[meshtrig] ( -0.46573529073056247 , -0.6835123756350594 ) -- ( -0.295081547266559 , -0.26766310993246095 ) -- ( -0.6971923744228338 , -0.36057467183152875 ) -- cycle;
\draw[meshtrig] ( 0.8212531032817161 , -0.668206736929906 ) -- ( 1.1578746161954439 , -0.44698643381097225 ) -- ( 0.6829318808011956 , -0.1985703945312043 ) -- cycle;
\draw[meshtrig] ( -0.650560958288211 , 0.02488343392998438 ) -- ( -0.6971923744228338 , -0.36057467183152875 ) -- ( -0.295081547266559 , -0.26766310993246095 ) -- cycle;
\draw[meshtrig] ( -0.32685022560532856 , 0.6762126992658455 ) -- ( -0.662003809555919 , 0.4617132239118765 ) -- ( -0.19510386774107463 , 0.2177653724852468 ) -- cycle;
\draw[meshtrig] ( 0.42125922650351183 , -0.6311516272278445 ) -- ( 0.8212531032817161 , -0.668206736929906 ) -- ( 0.6829318808011956 , -0.1985703945312043 ) -- cycle;
\draw[meshtrig] ( 1.1792544350003562 , 0.4754657340609835 ) -- ( 0.8537612867784096 , 0.7056698463200347 ) -- ( 0.7698780756308289 , 0.2979251021457717 ) -- cycle;
\draw[meshtrig] ( 0.8537612867784096 , 0.7056698463200347 ) -- ( 0.4772704753189263 , 0.6599589103990077 ) -- ( 0.7698780756308289 , 0.2979251021457717 ) -- cycle;
\draw[meshtrig] ( 0.4772704753189263 , 0.6599589103990077 ) -- ( 0.31947279823770486 , 0.2274752592073914 ) -- ( 0.7698780756308289 , 0.2979251021457717 ) -- cycle;
\draw[meshtrig] ( 0.42125922650351183 , -0.6311516272278445 ) -- ( 0.6829318808011956 , -0.1985703945312043 ) -- ( 0.16139908952035206 , -0.2112520920264157 ) -- cycle;
\draw[meshtrig] ( -0.662003809555919 , 0.4617132239118765 ) -- ( -0.650560958288211 , 0.02488343392998438 ) -- ( -0.19510386774107463 , 0.2177653724852468 ) -- cycle;
\draw[meshtrig] ( -0.650560958288211 , 0.02488343392998438 ) -- ( -0.295081547266559 , -0.26766310993246095 ) -- ( -0.19510386774107463 , 0.2177653724852468 ) -- cycle;
\draw[meshtrig] ( 0.16139908952035206 , -0.2112520920264157 ) -- ( -0.19510386774107463 , 0.2177653724852468 ) -- ( -0.295081547266559 , -0.26766310993246095 ) -- cycle;
\draw[meshtrig] ( 0.7698780756308289 , 0.2979251021457717 ) -- ( 0.31947279823770486 , 0.2274752592073914 ) -- ( 0.6829318808011956 , -0.1985703945312043 ) -- cycle;
\draw[meshtrig] ( 0.16139908952035206 , -0.2112520920264157 ) -- ( 0.31947279823770486 , 0.2274752592073914 ) -- ( -0.19510386774107463 , 0.2177653724852468 ) -- cycle;
\draw[meshtrig] ( 0.16139908952035206 , -0.2112520920264157 ) -- ( 0.6829318808011956 , -0.1985703945312043 ) -- ( 0.31947279823770486 , 0.2274752592073914 ) -- cycle;
\draw[Gamma_hbeg] ( -0.5325604384457461 , -0.8457754536867813 ) -- ( -0.2616833384212279 , -0.9483211465731678 );
\draw[Gamma_hbeg] ( -0.15683456406001578 , -0.986023442956576 ) -- ( 0.20069037679128188 , -0.9753830463086911 );
\draw[Gamma_hbeg] ( -0.15683456406001578 , -0.986023442956576 ) -- ( -0.2616833384212279 , -0.9483211465731678 );
\draw[Gamma_hbeg] ( 0.2927067577173778 , -0.9455521618221736 ) -- ( 0.549068096994426 , -0.8354165127540144 );
\draw[Gamma_hbeg] ( 0.2927067577173778 , -0.9455521618221736 ) -- ( 0.20069037679128188 , -0.9753830463086911 );
\draw[Gamma_hbeg] ( 0.549068096994426 , -0.8354165127540144 ) -- ( 0.7426508347836078 , -0.6609250862569299 );
\draw[Gamma_hbeg] ( 0.5644174942428252 , 0.8249130343206474 ) -- ( 0.3184413703054212 , 0.9325410708969099 );
\draw[Gamma_hbeg] ( 0.5644174942428252 , 0.8249130343206474 ) -- ( 0.715415619164979 , 0.6888728628220462 );
\draw[Gamma_hbeg] ( 0.2157176305351797 , 0.9748425134648014 ) -- ( -0.1441933352082881 , 0.9855107198567988 );
\draw[Gamma_hbeg] ( 0.2157176305351797 , 0.9748425134648014 ) -- ( 0.3184413703054212 , 0.9325410708969099 );
\draw[Gamma_hbeg] ( -0.21358143271716906 , 0.9712193405108395 ) -- ( -0.4811439631725012 , 0.8752112083711706 );
\draw[Gamma_hbeg] ( -0.21358143271716906 , 0.9712193405108395 ) -- ( -0.1441933352082881 , 0.9855107198567988 );
\draw[Gamma_hbeg] ( -0.8413141516556111 , 0.5394860478065603 ) -- ( -0.9452667919012407 , 0.2779458564554428 );
\draw[Gamma_hbeg] ( -0.9820314426361334 , 0.17850056047097587 ) -- ( -0.9821389909978679 , -0.17009547099836247 );
\draw[Gamma_hbeg] ( -0.9820314426361334 , 0.17850056047097587 ) -- ( -0.9452667919012407 , 0.2779458564554428 );
\draw[Gamma_hbeg] ( -0.5325604384457461 , -0.8457754536867813 ) -- ( -0.7195974847326421 , -0.6602613237240756 );
\draw[Gamma_hbeg] ( -0.6407332766726979 , 0.7377159099680539 ) -- ( -0.8413141516556111 , 0.5394860478065603 );
\draw[Gamma_hbeg] ( -0.9621407849762262 , -0.24749741079983237 ) -- ( -0.8694260179331358 , -0.49210274140910976 );
\draw[Gamma_hbeg] ( -0.9621407849762262 , -0.24749741079983237 ) -- ( -0.9821389909978679 , -0.17009547099836247 );
\draw[Gamma_hbeg] ( 0.9748776606619196 , 0.1426190169344961 ) -- ( 0.9300064979469389 , 0.36737049300657154 );
\draw[Gamma_hbeg] ( 0.9414894469590955 , -0.3338074278073215 ) -- ( 0.9902738417962321 , -0.054516369764659414 );
\draw[Gamma_hbeg] ( 0.9748776606619196 , 0.1426190169344961 ) -- ( 0.9902738417962321 , -0.054516369764659414 );
\draw[Gamma_hbeg] ( -0.6407332766726979 , 0.7377159099680539 ) -- ( -0.4811439631725012 , 0.8752112083711706 );
\draw[Gamma_hbeg] ( -0.7195974847326421 , -0.6602613237240756 ) -- ( -0.8694260179331358 , -0.49210274140910976 );
\draw[Gamma_hbeg] ( 0.7977618601807144 , -0.5884478879947267 ) -- ( 0.9414894469590955 , -0.3338074278073215 );
\draw[Gamma_hbeg] ( 0.7426508347836078 , -0.6609250862569299 ) -- ( 0.7977618601807144 , -0.5884478879947267 );
\draw[Gamma_hbeg] ( 0.9300064979469389 , 0.36737049300657154 ) -- ( 0.8217285286966978 , 0.5499630180937727 );
\draw[Gamma_hbeg] ( 0.715415619164979 , 0.6888728628220462 ) -- ( 0.8217285286966978 , 0.5499630180937727 );
\draw[Gamma_hend] ( -0.5243399164161836 , -0.8258145968339992 ) -- ( -0.27486605922334223 , -0.931213248051638 );
\draw[Gamma_hend] ( -0.1572533792636008 , -0.9870345976440426 ) -- ( 0.19847864364471224 , -0.9720361785162407 );
\draw[Gamma_hend] ( -0.1572533792636008 , -0.9870345976440426 ) -- ( -0.27486605922334223 , -0.931213248051638 );
\draw[Gamma_hend] ( 0.2892702296533673 , -0.9539568721375452 ) -- ( 0.5569980738913057 , -0.8480902482846124 );
\draw[Gamma_hend] ( 0.2892702296533673 , -0.9539568721375452 ) -- ( 0.19847864364471224 , -0.9720361785162407 );
\draw[Gamma_hend] ( 0.5569980738913057 , -0.8480902482846124 ) -- ( 0.8124300937690206 , -0.6673893804542803 );
\draw[Gamma_hend] ( 0.5727304558956006 , 0.8406480224310864 ) -- ( 0.3113443790501913 , 0.9447209116711861 );
\draw[Gamma_hend] ( 0.5727304558956006 , 0.8406480224310864 ) -- ( 0.7741800909395888 , 0.6960076425350201 );
\draw[Gamma_hend] ( 0.21488830975341394 , 0.9730243875166262 ) -- ( -0.14559691438800168 , 0.9878233364960367 );
\draw[Gamma_hend] ( 0.21488830975341394 , 0.9730243875166262 ) -- ( 0.3113443790501913 , 0.9447209116711861 );
\draw[Gamma_hend] ( -0.21610195893659165 , 0.9646546729161561 ) -- ( -0.4736048612554368 , 0.8654877415270561 );
\draw[Gamma_hend] ( -0.21610195893659165 , 0.9646546729161561 ) -- ( -0.14559691438800168 , 0.9878233364960367 );
\draw[Gamma_hend] ( -0.7690190346446344 , 0.508129265270201 ) -- ( -0.8053397706434955 , 0.36872376227163134 );
\draw[Gamma_hend] ( -0.8671152181345365 , 0.12524363101735514 ) -- ( -0.8672047634653507 , -0.10251032478121684 );
\draw[Gamma_hend] ( -0.8671152181345365 , 0.12524363101735514 ) -- ( -0.8053397706434955 , 0.36872376227163134 );
\draw[Gamma_hend] ( -0.5243399164161836 , -0.8258145968339992 ) -- ( -0.6499225882934572 , -0.6666427963193476 );
\draw[Gamma_hend] ( -0.6464323114401795 , 0.6637662392010715 ) -- ( -0.7690190346446344 , 0.508129265270201 );
\draw[Gamma_hend] ( -0.8321978029619526 , -0.30295574121885555 ) -- ( -0.7983497415833583 , -0.43782458047309175 );
\draw[Gamma_hend] ( -0.8321978029619526 , -0.30295574121885555 ) -- ( -0.8672047634653507 , -0.10251032478121684 );
\draw[Gamma_hend] ( 1.120540358307192 , 0.03226609437043336 ) -- ( 1.0176821477484226 , 0.4053941597696086 );
\draw[Gamma_hend] ( 1.0299019573471413 , -0.38005107786578207 ) -- ( 1.1232051755180021 , 0.007789778013563897 );
\draw[Gamma_hend] ( 1.120540358307192 , 0.03226609437043336 ) -- ( 1.1232051755180021 , 0.007789778013563897 );
\draw[Gamma_hend] ( -0.6464323114401795 , 0.6637662392010715 ) -- ( -0.4736048612554368 , 0.8654877415270561 );
\draw[Gamma_hend] ( -0.6499225882934572 , -0.6666427963193476 ) -- ( -0.7983497415833583 , -0.43782458047309175 );
\draw[Gamma_hend] ( 0.819115876769187 , -0.6609502998899349 ) -- ( 1.0299019573471413 , -0.38005107786578207 );
\draw[Gamma_hend] ( 0.8124300937690208 , -0.6673893804542803 ) -- ( 0.819115876769187 , -0.6609502998899349 );
\draw[Gamma_hend] ( 1.0176821477484226 , 0.4053941597696086 ) -- ( 0.8388952870050037 , 0.6334082674942362 );
\draw[Gamma_hend] ( 0.7741800909395888 , 0.6960076425350201 ) -- ( 0.8388952870050037 , 0.6334082674942362 );

%% file: num_exp/nze2_table_k1.tex
DG($k\!=\!1$): 
& 3.57K&4.24K &  $1.3\!\cdot\!10^{-1}$
& 3.3K&3.31K &  $3.4\!\cdot\!10^{-2}$
& 2.79K&3.14K &  $3.3\!\cdot\!10^{-2}$
\\ 
CG($k\!=\!1$): 
& 1.63K&1.72K &  $1.7\!\cdot\!10^{-1}$
& 978&1.6K &  $3.3\!\cdot\!10^{-2}$
& 728&867 &  $3.3\!\cdot\!10^{-2}$
\\ 

%% file: num_exp/nze2_table_k3.tex
DG($k\!=\!3$): 
& 353K&442K &  $6.6\!\cdot\!10^{-4}$
& 293K&327K &  $2.6\!\cdot\!10^{-4}$
& 267K&300K &  $2.7\!\cdot\!10^{-4}$
\\ 
CG($k\!=\!3$): 
& 221K&268K &  $1.6\!\cdot\!10^{-3}$
& 172K&190K &  $2.9\!\cdot\!10^{-4}$
& 150K&171K &  $2.7\!\cdot\!10^{-4}$
\\ 
GCC($k\!=\!3$): 
& 145K&151K &  $6.1\!\cdot\!10^{-3}$
& 84.5K&92.2K &  $3.1\!\cdot\!10^{-4}$
& 62.9K&76.1K &  $2.8\!\cdot\!10^{-4}$
\\ 

%% file: num_exp/nze2_table_k5.tex
DG($k\!=\!5$): 
& 4.39M&5.58M &  $5.4\!\cdot\!10^{-6}$
& 3.60M&4.5M &  $2.8\!\cdot\!10^{-7}$
& 3.26M&3.69M &  $2.2\!\cdot\!10^{-7}$
\\ 
CG($k\!=\!5$): 
& 3.71M&4.43M &  $1.7\!\cdot\!10^{-5}$
& 2.69M&2.99M &  $1.2\!\cdot\!10^{-6}$
& 2.28M&2.63M &  $4.3\!\cdot\!10^{-7}$
\\ 